\documentclass[a4paper,12pt]{article}
\usepackage[all]{xy}
\usepackage[T1]{fontenc}
\usepackage[dvips]{graphicx}
\usepackage{amsmath,amssymb,amsfonts,amsthm,latexsym,mathrsfs,textcomp,verbatim, enumerate}
\xyoption{web}
\begin{document}

\title{Topological Galois Theory}
\author{Olivia Caramello}
\bgroup           % fake a titlepage 
\let\footnoterule\relax  % no rule above thanks footnotes 
\date{January 2, 2013}
\maketitle

\begin{abstract}
We introduce an abstract topos-theoretic framework for building Galois-type theories in a variety of different mathematical contexts; such theories are obtained from representations of certain atomic two-valued toposes as toposes of continuous actions of a topological group. Our framework subsumes in particular Grothendieck's Galois theory and allows to build Galois-type equivalences in new contexts, such as for example graph theory and finite group theory.        
\end{abstract} 
\vspace{1cm}
\egroup

%MACROS-----------------------------------------------------------------------------------------------------------------------

%	European dates ``19 April 1990'' not ``April 19, 1990''
\def\Monthnameof#1{\ifcase#1\or
   January\or February\or March\or April\or May\or June\or
   July\or August\or September\or October\or November\or December\fi}
\def\today{\number\day~\Monthnameof\month~\number\year}

%===========================================================================
%	END OF PROOF BOX
%
%
%  The complexity of the macro necessary to get a little box on the
%  right-hand-side at the end of a proof is amazing.  It really does
%  have to be this long!  Otherwise you're liable to get it at the
%  beginning of the next line, or even on the next page.
%
\def\pushright#1{{%        set up
   \parfillskip=0pt            % so \par doesnt push \square to left
   \widowpenalty=10000         % so we dont break the page before \square
   \displaywidowpenalty=10000  % ditto
   \finalhyphendemerits=0      % TeXbook exercise 14.32
  %
  %                 horizontal
   \leavevmode                 % \nobreak means lines not pages
   \unskip                     % remove previous space or glue
   \nobreak                    % don't break lines
   \hfil                       % ragged right if we spill over
   \penalty50                  % discouragement to do so
   \hskip.2em                  % ensure some space
   \null                       % anchor following \hfill
   \hfill                      % push \square to right
   {#1}                        % the end-of-proof mark (or whatever)
  %
  %                   vertical
   \par}}                      % build paragraph

% prefer proofs with statements, also space after
\def\qed{\pushright{$\square$}\penalty-700 \smallskip}

\newtheorem{theorem}{Theorem}[section]

\newtheorem{proposition}[theorem]{Proposition}

\newtheorem{scholium}[theorem]{Scholium}

\newtheorem{lemma}[theorem]{Lemma}

\newtheorem{corollary}[theorem]{Corollary}

\newtheorem{conjecture}[theorem]{Conjecture}

\newenvironment{proofs}%
 {\begin{trivlist}\item[]{\bf Proof }}%
 {\qed\end{trivlist}}

  \newtheorem{rmk}[theorem]{Remark}
\newenvironment{remark}{\begin{rmk}\em}{\end{rmk}}

  \newtheorem{rmks}[theorem]{Remarks}
\newenvironment{remarks}{\begin{rmks}\em}{\end{rmks}}

  \newtheorem{defn}[theorem]{Definition}
\newenvironment{definition}{\begin{defn}\em}{\end{defn}}

  \newtheorem{eg}[theorem]{Example}
\newenvironment{example}{\begin{eg}\em}{\end{eg}}

  \newtheorem{egs}[theorem]{Examples}
\newenvironment{examples}{\begin{egs}\em}{\end{egs}}

%%%%%%%%%%%%%%%%%%%%%%%%%%%%%%%%%%%%%%%%%%%%%%%%%%%%%%%%%%%%%%%%%%%%%%

%  change some single-character symbols to be more appropriate for logic
\mathcode`\<="4268  % < = \langle
\mathcode`\>="5269  % > = \rangle
\mathcode`\.="313A  % make . a binary  relation
\mathchardef\semicolon="603B % the original
\mathchardef\gt="313E
\mathchardef\lt="313C

\newcommand{\cod}% codomain
 {{\rm cod}}

\newcommand{\comp}% composition
 {\circ}

\newcommand{\Cont}% category of continuous G-sets
 {{\bf Cont}}

\newcommand{\dom}% domain
 {{\rm dom}}

\newcommand{\empstg}% empty string
 {[\,]}

\newcommand{\epi}% epimorphism
 {\twoheadrightarrow}

\newcommand{\hy}% hyphen (in math mode)
 {\mbox{-}}

\newcommand{\im}% image
 {{\rm im}}

\newcommand{\imp}% implication
 {\!\Rightarrow\!}

\newcommand{\Ind}[1]% ind-completion of #1
 {{\rm Ind}\hy #1}

\newcommand{\mono}% monomorphism 
 {\rightarrowtail}

\newcommand{\ob}% class of objects
 {{\rm ob}}

\newcommand{\op}% opposite category
 {^{\rm op}}

\newcommand{\Set}% category of sets
 {{\bf Set }}

\newcommand{\Sh}% category of sheaves
 {{\bf Sh}}

\newcommand{\sh}% category of sheaves
 {{\bf sh}}

\newcommand{\Sub}% subobject lattice
 {{\rm Sub}}

%\newenvironment{changemargin}[2]{
% \begin{list}{}{  \setlength{\topsep}{-100pt}  \setlength{\leftmargin}{#1}  \setlength{\rightmargin}{#2}  \setlength{\listparindent}{\parindent} \setlength{\itemindent}{\parindent}  \setlength{\parsep}{\parskip} }\item[]}{\end{list}}

%\begin{changemargin}{-2.3cm}{-2.3cm} 
\tableofcontents

\newpage
\section{Introduction}\label{intro}

The present work provides a general framework, based on Topos Theory, for building Galois-type theories in a variety of different mathematical contexts. 

Let us start by briefly recalling the classical (infinite) Galois theory and its categorical interpretation. 

Let $F \subseteq L$ be a Galois extension, not necessarily finite-dimensional. Then the group $Aut_{F}(L)$ of automorphisms of $L$ which fix $F$ can be naturally made into a topological group by endowing it with the so-called \emph{Krull topology}, that is the topology in which the subgroups of $Aut_{F}(L)$ consisting of the automorphisms which fix a given finite family of elements of $L$ form a basis of open neighbourhoods of the identity.  

The fundamental theorem of infinite Galois provides an order-reversing bijective correspondence between the intermediate extensions $K$ and the closed subgroups of the Galois group $Aut_{F}(L)$, which restricts to a bijective correspondence between the \emph{finite extensions} and the \emph{open subgroups} of $Aut_{F}(L)$. The correspondence assigns to any intermediate extension $K$ the subgroup of $Aut_{F}(L)$ consisting of the automorphisms which fix $K$, and conversely associates to any closed subgroup $C$ of $Aut_{F}(L)$ the field extension of $F$ consisting of all the elements of $L$ which are fixed by every automorphism in $C$. In fact, the correspondence between finite extensions and open subgroups can be seen as the `kernel' of the extended correspondence between the intermediate extensions and the closed subgroups, as the latter can be immediately obtained from it by observing that, on the one hand, every intermediate extension can be expressed as a union of finite extensions, and on the other, every closed subgroup of $Aut_{F}(L)$ is equal to the intersection of all the open subgroups containing it. The advantage of restricting our attention to this `kernel' of the correspondence is that, as it was already observed by Grothendieck in \cite{grothendieckg}, it admits a natural categorical `externalization', that is a formulation as a categorical duality between the category ${\cal L}_{F}^{L}$ of intermediate extensions and field homomorphisms between them and the category $\Cont_{t}(Aut_{F}(K))$ of transitive (non-empty) continuous actions of the topological group $Aut_{F}(L)$ over a discrete (finite) set.

The point of depart of the present work is the observation that the above-mentioned equivalence ${{\cal L}_{F}^{L}}^{\textrm{op}}\simeq \Cont_{t}(Aut_{F}(K))$ naturally extends to an equivalence of toposes 
\[
\Sh({{\cal L}_{F}^{L}}^{\textrm{op}}, J_{at})\simeq \Cont(Aut_{F}(K)),
\]
where $\Cont(Aut_{F}(K))$ is the topos of continuous actions (over discrete sets) of the topological group $Aut_{F}(K)$ and $J_{at}$ is the atomic topology on ${{\cal L}_{F}^{L}}^{\textrm{op}}$ (recall that the \emph{atomic topology} on a category satisfying the property that any pair of arrows with common codomain can be completed to a commutative square is the Grothendieck topology whose covering sieves are exactly the non-empty ones). Indeed, $\Cont_{t}(Aut_{F}(K))$ is a $J_{can}$-dense subcategory of $\Cont(Aut_{F}(K))$, where $J_{can}$ is the canonical topology on $\Cont(Aut_{F}(K))$, and the Grothendieck topology induced by $J_{can}$ on $\Cont_{t}(Aut_{F}(K))$ coincides with the atomic topology on $\Cont_{t}(Aut_{F}(K))$ (notice that the objects of the subcategory $\Cont_{t}(Aut_{F}(K))$ are precisely the atoms of the topos $\Cont(Aut_{F}(K))$), whence the Comparison Lemma yields an equivalence 

\[
\Sh(\Cont_{t}(Aut_{F}(K)), J_{at})\simeq \Cont(Aut_{F}(K))
\]            

This remark suggests that a natural context for building analogues of classical Galois theory in different mathematical domains might be provided by equivalences of toposes of the form
\[
\Sh({\cal C}^{\textrm{op}}, J_{at})\simeq \Cont(Aut_{\cal D}(u)),
\]  
where $\cal C$ is a category satisfying the \emph{amalgamation property} (i.e. the property that every pair or arrows with common domain can be completed to a commutative square), $J_{at}$ is the atomic topology on it, $\cal D$ is a category in which $\cal C$ embeds, $u$ is an object of $\cal D$ and $Aut_{\cal D}(u)$ is the group of automorphisms of $u$ in $\cal D$ endowed with a topology in which the collection of subgroups of the form $\{f:u\cong u \textrm{ | } f\circ \chi=\chi\}$ for arrows $\chi:c\to u$ in $\cal C$ forms a basis of open neighbourhoods of the identity. Notice that, by Diaconescu's theorem, such an equivalence must be induced by a functor $F:{\cal C}^{\textrm{op}}\to \Cont(Aut_{\cal D}(u))$ - necessarily taking values in the subcategory $\Cont_{t}(Aut_{\cal D}(u))$ of $\Cont(Aut_{\cal D}(u))$ on the transitive (non-empty) continuous $Aut_{\cal D}(u)$-actions (since the images of representable functors on $\cal C$ via the associated sheaf functor $a_{J_{at}}:[{\cal C}, \Set]\to \Sh({\cal C}^{\textrm{op}}, J_{at})$ are atoms of the topos $\Sh({\cal C}^{\textrm{op}}, J_{at})$). In light of classical Galois theory, it is natural to require this functor to be the one sending any object of $\cal C$ to the set $Hom_{\cal D}(c, u)$, with the obvious action of $Aut_{\cal D}(u)$; and any arrow $f:d\to c$ of $\cal C$ to the $Aut_{\cal D}(u)$-equivariant map $-\circ f:Hom_{\cal D}(c, u) \to Hom_{\cal D}(d, u)$.     

Some conditions on $\cal C$, $\cal D$ and $u$ are necessary in order for such an equivalence to exist:

\begin{enumerate}[(i)]
\item $\cal C$ must be non-empty and satisfy the \emph{joint embedding property}, i.e. the property that for any pair of objects $a,b\in {\cal C}$ there exists an object $c\in {\cal C}$ and arrows $a\to c$ and $b\to c$ in $\cal C$; indeed, any topos of the form $\Cont(G)$ for a topological group $G$ is atomic and two valued, and it is easy to see that a topos of the form $\Sh({\cal C}^{\textrm{op}}, J_{at})$, for a non-empty category $\cal C$ satisfying the amalgamation property, is (atomic and) two-valued if and only if $\cal C$ satisfies the joint embedding property;  

\item The object $u$ must satisfy the property that any object $c$ of $\cal C$ admits an arrow to $u$ in $\cal D$ ($u$ is $\cal C$-universal) and the action of $Aut_{\cal D}(u)$ on each $Hom_{\cal D}(c, u)$ is transitive ($u$ is $\cal C$-ultrahomogeneous). 

\end{enumerate}  

If ${\cal C}\hookrightarrow {\cal D}$ is the canonical embedding of $\cal C$ into its ind-completion $\Ind{\cal C}$, then our main representation theorem proved in section \ref{repres} assures that the two above-mentioned conditions are also sufficient, while Theorem \ref{catd} provides sufficient conditions, valid for any embedding ${\cal C}\hookrightarrow {\cal D}$, for such an equivalence to exist in the case the topology $J_{at}$ is subcanonical.   

In fact, besides our criteria, which are widely applicable in practice thanks in particular to the effectiveness and generality of Fra\"iss\'e's construction in Model Theory and its generalizations (cf. \cite{Hodges} and \cite{OC2}), there are other natural means for obtaining such equivalences of toposes: exploiting the theory of special models for atomic complete theories, Grothendieck's Galois Theory or the representation theory of Grothendieck topos (cf. section \ref{four}).  

Given an equivalence 
\[
\Sh({\cal C}^{\textrm{op}}, J_{at})\simeq \Cont(Aut_{\cal D}(u)),
\] 
induced by a functor $F:{\cal C}^{\textrm{op}}\to \Cont(Aut_{\cal D}(u))$ as above, the functor $F$ is full and faithful if and only if the composite of the Yoneda embedding ${\cal C}^{\textrm{op}}\to [{\cal C}, \Set]$ with the associated sheaf functor $a_{J_{at}}:[{\cal C}^{\textrm{op}}, \Set]\to \Sh({\cal C}^{\textrm{op}}, J_{at})$ is full and faithful, i.e. if and only if the topology $J_{at}$ is subcanonical. This condition can be formulated concretely as the condition that all the arrows of $\cal C$ are strict monomorphisms, and under this hypotheses $F$ yields an essentially full and faithful functor from the opposite of the slice category ${\cal C}\slash u$ to the poset category of open subgroups of the group $Aut_{\cal D}(d)$; in other words, the following Galois-type property holds: for any arrows $\chi:c\to u$ and $\xi:d\to u$ from objects $c$ and $d$ of $\cal C$ to $u$, if all the automorphisms of $u$ which fix $\chi$ also fix $\xi$ then $\xi$ factors uniquely through $\chi$ in $\cal D$. We shall speak in this context of a \emph{concrete Galois theory}. Examples of such theories abound (cf. section \ref{examples}); in particular, there exists a concrete Galois theory for finite groups. 

Under the assumption that $F:{\cal C}^{\textrm{op}}\to \Cont_{t}(Aut_{\cal D}(u))$ is full and faithful, it is natural to wonder under which conditions $F$ is also essentially surjective and hence yields a Galois-type correspondence between the isomorphism classes of objects of the slice category ${\cal C}\slash u$ and the open subgroups of the topological group $Aut_{\cal D}(u)$. We shall speak in this case of a \emph{standard Galois theory}. By considering the invariant notion of atom in terms of the two toposes related by the equivalence, we arrive at a purely elementary necessary and sufficient condition on the category $\cal C$ for the associated functor $F$ to be essentially surjective, which can be interpreted as a form of `elimination of imaginaries' for $\cal C$ in the sense of Model Theory; also, we describe for a general $\cal C$ its completion ${\cal C}_{at}$ relatively to `imaginaries', for which a standard Galois equivalence holds. 

The scope of applicability of this framework is much broader than that of the Galois theory of Grothendieck of \cite{grothendieckg} (as well as of its infinitary generalization established in \cite{noohi}). Indeed, the latter corresponds to choosing categories $\cal C$ which are pretoposes (finitary in the case of \cite{grothendieckg} or infinitary in the case of \cite{noohi}) and the automorphisms groups arising from such constructions are always topologically prodiscrete (and in the finitary context profinite).   

Also, interpreting Galois-type equivalences in terms of different representations of a given topos paves the way for a systematic use of topos-theoretic invariants for transferring notions and results across the two sides, according to the methodologies `toposes as bridges' of \cite{OC10}; such an analysis is carried out in the present paper for several invariants, leading to insights of various nature on the given Galois theories (notably including the above-mentioned theorems).           

From a logical point of view, the `Galois-type' theories (i.e., the geometric theories classified by a topos of continuous actions of a topological group) are maximal with respect to the natural (Heyting algebra) ordering between geometric theories over a given signature defined in \cite{OC6} (${\mathbb T} \leq {\mathbb T}'$ if and only if ${\mathbb T}'$ is a quotient of $\mathbb T$, i.e. every geometric sequent which is provable in $\mathbb T$ is also provable in ${\mathbb T}$'); in order words, they do not have any proper quotients over their signature. Conversely, by the results in \cite{OC5} and \cite{OC6} (cf. also \cite{blasce}), any maximal theory which is finitary or written over a countable signature and has a model in $\Set$ is a Galois-type theory. The question thus naturally arises of whether there is a canonical way for extending an arbitrary geometric theory to a Galois-type theory. The results of \cite{OC2}, combined with those of the present paper, show that the answer to this question is positive for a very large class of theories. Specifically, for any theory of presheaf type $\mathbb T$ such that its category of finitely presentable models satisfies the amalgamation property, if there exists a homogeneous ${\mathbb T}$-model in $\Set$ and the theory $\mathbb T$ is finitary or written over a countable signature, there is a completion of the theory of homogeneous $\mathbb T$-models which is consistent (i.e., has a model in $\Set$) and hence is a Galois-type theory extending $\mathbb T$. 

This paper can be considered as a companion to \cite{dubuc3}, where a localic Galois theory for atomic connected toposes is developed. Nonetheless, it should be noted that there are important differences between the localic and the topological setting, the main ones being the following:

\begin{enumerate}[(i)]
\item Whilst any localic group gives naturally rise to a topological group, there are topological groups which do not arise in this way;

\item Whilst every atomic connected topos can be represented (by the main result in \cite{dubuc3}) as the topos of continuous actions of the localic automorphism group of \emph{any} of its points, an analogous representation in terms of the topological automorphism group of the point holds only if the point satisfies some special properties, as it is shown in this paper.   
\end{enumerate}

\section{Topological groups and their toposes of continuous actions}\label{topgroups}

Let us recall that a topological group is a group $G$ with a topology such that the group operation and the inverse operation are continuous with respect to it; for basic background on topological groups we refer the reader to \cite{DD}.

\subsection{Algebraic bases}

The following well-known result allows to make a given group into a topological group starting from a collection of subsets of the group satisfying particular properties:

\begin{lemma}
Let $G$ be a group and $\cal B$ be a collection of subsets $N$ of $G$ containing the neutral element $e$. Then there exists a topology $\tau$ on $G$ having $\cal B$ as a neighbourhood basis of $e$ and making $(G, \tau)$ into a topological group if and only if all the following conditions are satisfied:
\begin{enumerate}[(i)]
\item For any $N, M\in {\cal B}$ there exists $P\in {\cal B}$ such that $P\subseteq N\cap M$;
\item For any $N\in {\cal B}$ there exists $M\in {\cal B}$ such that $M^{2}\subseteq N$;
\item For any $N\in {\cal B}$ there exists $M\in {\cal B}$ such that $M\subseteq N^{-1}$;
\item For any $N\in {\cal B}$ and any $a\in G$ there exists $M\in {\cal B}$ such that $M\subseteq aNa^{-1}$.
\end{enumerate}
\end{lemma}

A couple of remarks:
\begin{itemize}
\item If the topology $\tau$ in the statement of the lemma exists then it is necessarily unique; indeed, if $\cal B$ is a basis of neighbourhoods of $e$ then for any $a\in G$, the collection of sets of the form $aN$ for a set $N$ in $\cal B$ forms a basis of neighbourhoods of $a$, whence the open sets in $\tau$ are exactly the unions of sets of the form $aN$ for $N\in {\cal B}$. 
 
We shall call $\tau$ the \emph{topology on $G$ generated by ${\cal B}$} and denote it by $\tau_{\cal B}^{G}$; we shall denote the resulting topological group by $G_{\cal B}$. 

\item If all the subsets in the family $\cal B$ are subgroups of $G$ then conditions $(ii)$ and $(iii)$ in the statement of the lemma are automatically satisfied. 

This motivates the following definition: we say that a collection $\cal B$ of subgroups of $G$ is an \emph{algebraic base} for $G$ if it is a basis of neighbourhood of $e$, any intersection of subgroups in $\cal B$ contains a subgroup in $\cal B$, and any conjungated of a subgroup in $\cal B$ lies again in $\cal B$.  
\end{itemize}
 
Let us denote by $\mathbf{GTop}$ the category of topological groups and continuous group homomorphisms between them. We can construct a category $\mathbf{GTop_{b}}$ of `groups paired with algebraic bases' as follows: the objects of $\mathbf{GTop_{b}}$ are pairs $(G, {\cal B})$ consisting of a group $G$ and an algebraic base $\cal B$ for it, while the arrows $(G, {\cal B})\to (G', {\cal B}')$ in $\mathbf{GTop_{b}}$ are the group homomorphisms $f:G\to G'$ such that for any $V\in {\cal B}'$, $f^{-1}(V)\in {\cal B}$. We have a functor $F:\mathbf{GTop_{b}} \to \mathbf{GTop}$ sending to any object $(G, {\cal B})$ of $\mathbf{GTop_{b}}$ the topological group $(G, \tau_{\cal B}^{G})$ and acting accordingly on arrows. On the other hand, any topological group $G$ has a canonical algebraic base $C_{G}$, namely the one consisting of all the open subgroups of it; this allows one to define a functor $G:\mathbf{GTop} \to \mathbf{GTop_{b}}$ sending $G$ to $(G, C_{G})$ and acting on arrows in the obvious way. It is easily verified that $G$ is left adjoint to $F$ and $F\circ G\cong 1_{\mathbf{GTop}}$, which allows us to regard $\mathbf{GTop}$ as a full subcategory of $\mathbf{GTop_{b}}$.      
 
For any topological group $G$, one can construct a topos $\Cont(G)$, whose objects are the left continuous actions $G\times X\to X$ (where $X$ is endowed with the discrete topology and $G\times X$ with the product topology) and whose arrows are the $G$-equivariant maps between them. Recall that a left action $\alpha:G\times X\to X$ is continuous if and only if for every $x\in X$ the isotropy subgroup $I_{x}:=\{g\in G \textrm{ | } \alpha(g, x)=x\}$ is open in $G$. The topos $\Cont(G)$ is atomic; in fact, its atoms are precisely its \emph{transitive} non-empty continuous actions. Notice that a non-empty transitive action $\alpha:G\times X\to X$ can be identified with the canonical action $G \times G\slash I_{x} \to G\slash I_{x}$ on the set $G\slash U$ of left cosets $gI_{x}$ of the isotropy group $I_{x}$ of $\alpha$ at any point $x\in X$; conversely, for any open subgroup $U$ of $G$, the canonical action of $G$ on the set $G\slash U$ makes $G\slash U$ into a non-empty transitive $G$-set. 

Notice that for any open subgroups $U$ and $V$ of $G$, the arrows $G\slash U \to G\slash V$ in $\Cont(G)$ can be identified with their action on the trivial coset, that is with the $V$-cosets of the form $aV$ where $a$ is an element of $G$ such that $U\subseteq a^{-1}Va$. Such an arrow is an isomorphism if and only if $U=a^{-1}Va$; indeed, it is invertible if and only if there exists $b\in G$ such that $ba\in U$, $ab\in V$ and $V\subseteq b^{-1}Ub$; but $ba\in U$ implies that $a^{-1}$ and $b$ are equivalent modulo $U$ and hence $V\subseteq b^{-1}Ub$ is equivalent to $V\subseteq {a^{-1}}^{-1} U a^{-1}= aUa^{-1}$.          
We denote by $\Cont_{t}(G)$ the full subcategory of $\Cont(G)$ on the $G$-sets of the form $G\slash U$.

Notice that two continuous $G$-sets $X$ and $Y$ are isomorphic in $\Cont(G)$ if and only if the the sets of isotropy subgroups at elements of $X$ and of isotropy subgroups at elements of $Y$ contain exactly the same subgroups as elements, if the latter are considered up to conjugation (by some element of $G$); in particular, any two transitive continuous $G$-sets are isomorphic if and only if any two isotropy subgroups of them are conjugate to each other.

\subsection{Algebraic bases and dense subcategories of actions}
  
There is a natural link between algebraic bases for a topological group and dense subcategories of the associated topos of continuous actions.
  
\begin{remark}\label{algbases}
For any algebraic base for $G$, the $G$-sets of the form $G\slash U$ for $U\in {\cal B}$ define a dense full subcategory of $\Cont_{t}(G)$ (in the sense that for any object of $\Cont_{t}(G)$ there exists an arrow from a $G$-set of this form to it) which is closed under isomorphisms. Conversely, any dense full subcategory of $\Cont_{t}(G)$ which is closed under isomorphisms gives rise to an algebraic base for $G$, namely the base consisting of the open subgroups $U$ of $G$ such that $G\slash U$ lies in the subcategory. The algebraic bases for $G$ can be thus identified with the dense full subcategories of $\Cont_{t}(G)$ which are closed under isomorphisms.    
\end{remark}

\begin{proposition}\label{propAPJEP}
For any algebraic base $B$ of a group $G$, the full subcategory $\Cont_{{\cal B}}(G)$ of $\Cont_{t}(G)$ on the objects of the form $G\slash U$ for $U\in {\cal B}$ satisfies the dual of the amalgamation property and the dual of the joint embedding property.
\end{proposition}

\begin{proofs}
To prove that $\Cont_{{\cal B}}(G)$ satisfies the dual of the amalgamation property it suffices to recall for any atomic topos $\cal E$ and any separating set of atoms of $\cal E$, the full subcategory $\cal A$ of $\cal E$ on this set of atoms satisfies the dual of the amalgamation property (this can be proved by observing that the pullback of two arrows between atoms in $\cal E$ cannot be zero, since the domains of the two arrows are epimorphic images of it, and hence there is an arrow to it from an atom in the separating set). 
  
On the other hand, $\Cont_{{\cal B}}(G)$ satisfies the dual of the joint embedding property, as for any two open subgroups $U$ and $V$ of $G$ belonging to $\cal B$, we have arrows $G\slash W \to G\slash U$ and $G\slash W \to G\slash V$ for any $W$ in $B$ such that $W\subseteq U\cap V$.  
\end{proofs}

Now, since the subcategory $\Cont_{{\cal B}}(G)$ is dense in $\Cont_{t}(G)$ and hence in the topos $\Cont(G)$, Grothendieck's Comparison Lemma yields an equivalence
\[
\Sh(\Cont_{{\cal B}}(G), J_{at}) \simeq \Cont(G),
\]
where $J_{at}$ is the atomic topology on $\Cont_{{\cal B}}(G)$.  
  
\subsection{Complete topological groups}\label{completegroups}  
  
We notice that the topos $\Cont(G)$ has a canonical point $p_{G}$, namely the geometric morphism $\Set \to \Cont(G)$ whose inverse image is the forgetful functor $\Cont(G)\to \Set$. Let us denote by $Aut(p_{G})$ the group of automorphisms of $p_{G}$ in the category of points of $\Cont(G)$; then we have a canonical map $\xi_{G}:G\to Aut(p_{G})$, sending any element $g\in G$ to the automorphism of $p_{G}$ which acts at each component as multiplication by the element $g$ (this is indeed an automorphism because the naturality conditions hold as the maps in $\Cont(G)$ are $G$-equivariant). 

As shown in \cite{noohi}, for any topological group $G$, the group $Aut(p_{G})$ can be intrinsically be endowed with a pro-discrete topology (that is a topology which is a projective limit of discrete topologies) in which the open subgroups are those subgroups of the form $U_{(X, x)}$ for a continuous $G$-sets $X$ and an element $x\in X$, where $U_{(X, x)}$ denotes the set of automorphisms $\alpha:p_{G}\cong p_{G}$ such that $\alpha(X)(x)=x$, and the canonical map $\xi_{G}:G \to Aut(p_{G})$ is continuous with respect to this topology.

It is natural to characterize the topological groups $G$ for which the map $\xi_{G}$ is a bijection (equivalently, a homeomorphism). Following the terminology of \cite{noohi}, we shall call such groups \emph{complete}, and we shall refer to the topological group $Aut(p_{G})$ as to the \emph{completion} of $G$. For any complete group $G$ with an algebraic base $\cal B$, we can alternatively describe the topology on $Aut(p_{G})$ induced by the topology on $G$ via the bijection $\xi_{G}$ as follows: a basis of open neighbourhoods of the identity is given by the sets of the form $\{\alpha:p_{G}\cong p_{G} \textrm{ | } \alpha(G\slash U)(eU)=eU\}$ for $U\in {\cal B}$. 

We have the following characterization of complete groups.

\begin{proposition}
Let $G$ be a topological group with an algebraic base $\cal B$. Then $G$ is complete if and only if for any assignment $U\to a_{U}$ of an element $a_{U}\in G\slash U$ to any subset $U\in {\cal B}$ such that for any $G$-equivariant map $m:G\slash U \to G\slash V$, $m(a_{U})=a_{V} m(eU)$ (equivalently, such that for any inclusion $U\subseteq V$ of open subgroups in $\cal B$, $a_{U}\equiv a_{V}$ modulo $V$ and for any open subgroup $U\in {\cal B}$ and any element $g\in G$ $a_{gUg^{-1}}\equiv a_{U}$ modulo $U$) there exists a unique $g\in G$ such that $a_{U}=gU$ for all $U\in {\cal B}$.   
\end{proposition}

\begin{proofs}
Since the category $\Cont_{t}(G)$ is dense in $\Cont_{t}(G)$, the automorphisms of $p_{G}$ correspond exactly to the isomorphisms of the flat functor $F:\Cont_{t}(G)\to \Set$ corresponding to $p_{G}$, which in fact can be identified with the forgetful functor. It thus remains to show that giving an automorphisms of $F$ corresponds to giving, for each $U\in {\cal B}$, an element $a_{U}\in G\slash U$ satisfying the hypotheses of the proposition. One direction is obvious; to prove the other one it suffices to observe that for any automorphism $\chi$ of $F$, $\chi_{U}([g])=g[\chi_{U}(e)]$ for any $U\in {\cal B}$; this follows as an immediate consequence of the following two remarks:

$(1)$ for any $U, V\in {\cal B}$ and any $G$-equivariant map $\gamma:G\slash V \to G\slash U$, we have $\chi_{U}(\gamma(eV))=\chi_{U}(eU) \gamma(eV)$. Indeed, by the naturality of $\chi$, we have $\chi_{U}(\gamma(eV))=\gamma(\chi_{V}(eV))$; but there exists $W\in {\cal B}$ such that $W\subseteq U\cap V$ whence there exists $g\in \chi_{W}(eW)$ such that $gU=\chi_{U}(eU)$ and $gV=\chi_{V}(eV)$, and by the equivariance of $\gamma$ we have $\gamma(\chi_{V}(eV))=\gamma(gV)=\gamma(g(eV))=g\gamma(eV)=\chi_{U}(eU)\gamma(eV)$, as required.    

$(2)$ The map $G\slash g^{-1}Ug \to G\slash U$ sending $e(g^{-1}Ug)$ to $gU$ is well-defined and equivariant.
     
\end{proofs}

We note that the uniqueness requirement in the statement of the proposition holds if and only if $G$ is \emph{nearly discrete}, that is if and only if the intersection of all the open subgroups of $G$ is equal to the singleton $\{e\}$.

There are two classes of complete groups: discrete groups (obviously) and profinite groups (by Grothendieck's Galois Theory \cite{grothendieckg}). We shall obtain an alternative characterization of complete groups in section \ref{functorialization} below. 

For any group $G$ and algebraic base $\cal B$ for $G$, the collection of subsets of the form ${\cal I}_{U, x}:=\{\alpha:p_{G}\cong p_{G} \textrm{ | } \alpha(G\slash U)(x)=x\}$ for $x\in G\slash U$ and $U\in {\cal B}$ forms an algebraic base for the group $Aut(p_{G})$ of automorphisms of $p_{G}$, and, if we consider $Aut(p_{G})$ endowed with the resulting topology, the canonical map $\xi_{G}:G\to Aut(p_{G})$ becomes a homomorphism of topological groups which induces a Morita-equivalence $\Cont(\xi_{G}):\Cont(G)\simeq \Cont(Aut(p_{G}))$ between them (cf. section \ref{functorialization} below).  

\section{A representation theorem}\label{repres}

In this section we establish our main representation theorem, which generalizes the representation theorem for coherent Boolean classifying toposes established by A. Blass and A. \v{S}\v{c}edrov in \cite{blasce}. We shall establish our theorem by using logical techniques, and later translate it into categorical language.

\subsection{Logical statement}

Recall from section D3.4 \cite{El} that, given a geometric theory $\mathbb T$ over a signature $\Sigma$, a geometric formula-in-context $\phi(\vec{x})$ over $\Sigma$ is said to be \emph{$\mathbb T$-complete} if the sequent $(\phi \vdash_{\vec{x}} \bot)$ is not provable in $\mathbb T$, but for any for any geometric formula $\psi(\vec{x})$ over $\Sigma$ in the same context, either $(\phi \vdash_{\vec{x}} \psi)$ is provable in $\mathbb T$ or $(\phi\wedge \psi \vdash_{\vec{x}} \bot)$ is provable in $\mathbb T$.  

\begin{theorem}\label{main}
Let $\mathbb T$ be an atomic and complete theory, ${\cal I}:=\{\phi_{i}(\vec{x_{i}}) \textrm{ | } i\in I\}$ be a set of ${\mathbb T}$-complete formulae such that for any $\mathbb T$-complete formula $\phi(\vec{x})$ there exists a  $\mathbb T$-provably functional formula (equivalently, $\mathbb T$-provably functional cover) from $\phi_{i}(\vec{x_{i}})$ to $\phi(\vec{x})$. Let $M$ be a model of $\mathbb T$ in $\Set$ in which each of the $\phi_{i}(\vec{x_{i}})$ is realized and such that for any $i\in I$ and any $\vec{a}, \vec{b}\in [[\vec{x_{i}}. \phi_{i}]]_{M}$ there exists an automorphism $f$ of $M$ such that $f(\vec{a})=\vec{b}$. Then, if we denote by $Aut_{\cal I}(M)$ the group of (${\mathbb T}$-model) automorphisms of $M$, we have that the sets of the form ${\cal I}_{\vec{a}}:=\{f:M\cong M \textrm{ | } f(\vec{a})=\vec{a}\}$, where $\vec{a}\in [[\vec{x_{i}}.\phi_{i}]]_{M}$ for some $i\in I$, form an algebraic base for $Aut_{\cal I}(M)$ and, if we endow $Aut_{\cal I}(M)$ with the resulting topology, we have an equivalence
\[
\Sh({\cal C}_{\mathbb T}, J_{\mathbb T})\simeq \Cont(Aut_{\cal I}(M))
\]
between the classifying topos of $\mathbb T$ and the topos of continuous $Aut_{\cal I}(M)$-sets (where $({\cal C}_{\mathbb T}, J_{\mathbb T})$ is the geometric syntactic site of $\mathbb T$), which restricts to an equivalence
\[
{\cal C}^{\cal I}_{\mathbb T} \simeq \Cont_{\cal I}(Aut_{\cal I}(M)),
\]
where ${\cal C}^{\cal I}_{\mathbb T}$ is the full subcategory of ${\cal C}_{\mathbb T}$ on the formulae of the form $\phi_{i}(\vec{x_{i}})$ and $\Cont_{\cal I}(Aut_{\cal I}(M))$ is the full subcategory of $\Cont(Aut_{\cal I}(M))$ on the $Aut_{\cal I}(M)$-sets isomorphic to one of the form $Aut_{\cal I}(M)\slash {\cal I}_{\vec{a}}$, which sends any formula $\phi_{i}(\vec{x_{i}})$ in ${\cal C}^{\cal I}_{\mathbb T}$ to the set $[[\phi_{i}(\vec{x_{i}})]]_{M}$ with the obvious $Aut_{\cal I}(M)$-action and a $\mathbb T$-provably functional formula $\theta$ from $\phi_{i}$ to $\phi_{j}$ to the $Aut_{\cal I}(M)$-equivariant map $[[\phi_{i}(\vec{x_{i}})]]_{M} \to [[\phi_{j}(\vec{x_{j}})]]_{M}$ whose graph is the interpretation $[[\theta]]_{M}$.   
\end{theorem}

\begin{proofs}
First, let us prove that the sets of the form ${\cal I}_{\vec{a}}$ form an algebraic base for the group $Aut_{\cal I}(M)$. For any $\vec{a}\in [[\phi_{i}(\vec{x_{i}})]]_{M}$ and $\vec{b}\in [[\phi_{j}(\vec{x_{j}})]]_{M}$, by the atomicity of $\mathbb T$ there exists a $\mathbb T$-complete formula $\chi(\vec{x_{i}}, \vec{x_{j}})$ such that $M \vDash \chi(\vec{a}, \vec{b})$. By our hypothesis, there exists a $\mathbb T$-provably functional formula $\theta(\vec{x_{k}}, \vec{x_{i}}, \vec{x_{j}})$ from a formula of the form $\phi_{k}(\vec{x_{k}})$ to $\chi(\vec{x_{i}}, \vec{x_{j}})$, from which it follows, $\chi$ being $\mathbb T$-complete, that there exists $\vec{c}\in [[\phi_{k}(\vec{x_{k}})]]_{M}$ such that $[[\theta]]_{M}(\vec{c})=(\vec{a}, \vec{b})$. Then ${\cal I}_{\vec{c}}\subseteq {\cal I}_{\vec{a}}\cap {\cal I}_{\vec{b}}$, since for any automorphism $f:M\cong M$, $f([[\theta]]_{M}(\vec{c}))=[[\theta]]_{M}(f(\vec{c}))$. To verify that every conjugate of a subgroup of the form ${\cal I}_{\vec{a}}$ is again a subgroup of this form, we observe that for any automorphism $h:M\cong M$, $h{\cal I}_{\vec{a}}h^{-1}=\{hfh^{-1} \textrm{ | } f(\vec{a})=\vec{a}\}=\{hfh^{-1} \textrm{ | } hfh^{-1}(h(\vec{a}))=h(\vec{a})\}$ and hence by symmetry $h{\cal I}_{\vec{a}}h^{-1}=\{s:M\cong M \textrm{ | } s(h(\vec{a}))=h(\vec{a})\}={\cal I}_{h(\vec{a})}$ (note that, $h$ being an automorphism, $h(\vec{a})\in [[\phi(\vec{x})]]_{M}$ if $\vec{a}\in [[\phi(\vec{x})]]_{M}$). 

From these remarks it follows in particular that the generated topology on $Aut_{\cal I}(M)$ coincides with the \emph{topology of pointwise convergence} on $Aut_{\cal I}(M)$, that is the topology whose basis of open neighbourhoods of the identity is given by the sets of the form $\{f:M\cong M \textrm{ | } f(\vec{a})=\vec{a}\}$ for \emph{any} $\vec{a}\in M$; indeed, as we have seen, any such set contains one of the form ${\cal I}_{\vec{c}}$.

Further, we note that our hypotheses imply that $M$ satisfies the more general property that for any $\vec{a}, \vec{b}\in M$ which satisfy exactly the same geometric (equivalently, the same $\mathbb T$-complete) formulae over the signature of $\mathbb T$, there exists an automorphism $f:M\cong M$ of $M$ such that $f({\vec{a}})=\vec{b}$. Indeed, if $\vec{a}, \vec{b}\in [[\chi(\vec{x})]]_{M}$ for a $\mathbb T$-complete formula $\chi(\vec{x})$ then there exists a $\mathbb T$-provably functional formula $\theta(\vec{x_{k}}, \vec{x})$ from a formula of the form $\phi_{k}(\vec{x_{k}})$ to $\chi(\vec{x})$, from which it follows, $\chi$ being $\mathbb T$-complete, that there exists $\vec{c}\in [[\phi_{k}(\vec{x_{k}})]]_{M}$ such that $[[\theta]]_{M}(\vec{c})=(\vec{a})$ and $\vec{d}\in [[\phi_{k}(\vec{x_{k}})]]_{M}$ such that $[[\theta]]_{M}(\vec{d})=(\vec{b})$; any automorphism $f:M\cong M$ such that $f(\vec{c})=\vec{d}$ will thus satisfy $f(\vec{a})=\vec{b}$. By using similar arguments, one can also prove that every $\mathbb T$-complete (equivalently, geometric) formula is satisfied in $M$. 

Now, since the sets ${\cal I}_{\vec{a}}$ form an algebraic base for the group $Aut_{\cal I}(M)$, by the results of section \ref{topgroups} we have an equivalence 
\[
\Cont(Aut_{\cal I}(M))\simeq \Sh(\Cont_{\cal I}(Aut_{\cal I}(M)), J_{at}),
\]
where $\Cont_{\cal I}(Aut_{\cal I}(M))$ is the full subcategory of $\Cont(Aut_{\cal I}(M))$ on the $Aut_{\cal I}(M)$-sets of the form $Aut_{\cal I}(M)\slash {\cal I}_{\vec{a}}$; therefore, showing that we have an equivalence $\Sh({\cal C}_{\mathbb T}, J_{\mathbb T})\simeq \Cont(Aut_{\cal I}(M))$ amounts precisely to proving that we have an equivalence 
\[
{\cal C}^{\cal I}_{\mathbb T} \simeq \Cont_{\cal I}(Aut_{\cal I}(M)),
\]
where ${\cal C}^{\cal I}_{\mathbb T}$ is the full subcategory of ${\cal C}_{\mathbb T}$ on the formulae of the form $\phi_{i}(\vec{x_{i}})$. To construct such an equivalence, we argue as follows. We have a functor $F:{\cal C}^{{\cal I}}_{\mathbb T} \to \Cont_{{\cal I}}(Aut_{\cal I}(M))$ sending any $\mathbb T$-complete formula $\phi(\vec{x})$ in ${\cal C}^{{\cal I}}_{\mathbb T}$ to the set $[[\phi(\vec{x})]]_{M}$, equipped with the obvious action by $Aut_{\cal I}(M)$ and any $\mathbb T$-provably functional formula $\theta$ from a formula $\phi(\vec{x})$ in $\cal I$ to a formula $\psi(\vec{y})$ in $\cal I$ to the $Aut_{\cal I}(M)$-equivariant map $[[\phi(\vec{x})]]_{M} \to [[\psi(\vec{y})]]_{M}$ whose graph coincides with the interpretation $[[\theta]]_{M}$ (this functor takes indeed values in $\Cont_{{\cal I}}(Aut_{\cal I}(M))$ by our hypotheses on $M$). Clearly, $F$ is essentially surjective by definition of the topology on $Aut_{\cal I}(M)$. To prove that $F$ is faithful, we observe that for any two $\mathbb T$-provably functional formulae $\theta_{1}(\vec{x}, \vec{y})$ and $\theta_{2}(\vec{x}, \vec{y})$ from a formula $\phi(\vec{x})$ in $\cal I$ to a formula $\psi(\vec{y})$ in $\cal I$, if $[[\theta_{1}(\vec{x}, \vec{y})]]_{M}=[[\theta_{2}(\vec{x}, \vec{y})]]_{M}$ then the sequent $(\phi(\vec{x}) \vdash_{\vec{x}} \exists \vec{y}(\theta_{1}(\vec{x}, \vec{y})\wedge \theta_{2}(\vec{x}, \vec{y})))$ is provable in $\mathbb T$ (since the conjunction of the antecedent and the consequent is satisfied in $M$, $\mathbb T$ is complete and $\phi(\vec{x})$ is $\mathbb T$-complete), from which it follows that $\theta_{1}$ and $\theta_{2}$ are provably equivalent in $\mathbb T$, that is they are equal as arrows $\phi(\vec{x})\to \psi(\vec{y})$ in ${\cal C}^{\cal I}_{\mathbb T}$. The fact that $F$ is full can be proved as follows. Let $\gamma:[[\phi(\vec{x})]]_{M}\to [[\psi(\vec{y})]]_{M}$ be a $Aut_{\cal I}(M)$-equivariant map, where $\phi(\vec{x})$ and $\psi(\vec{y})$ are formulae in $\cal I$. Take $(\vec{a}, \vec{b})$ belonging to the graph of this map (notice that this is possible since $\phi(\vec{x})$ is satisfied in $M$ by our hypotheses); by the atomicity of $\mathbb T$, there exists a $\mathbb T$-complete formula $\theta({\vec{x}}, \vec{y})$ such that $(\vec{a}, \vec{b})\in [[\theta(\vec{x}, \vec{y})]]_{M}$. To conclude the proof of fullness, it suffices to show that $\theta(\vec{x}, \vec{y})$ is $\mathbb T$-provably functional from $\phi(\vec{x})$ to $\psi(\vec{y})$, since $\gamma$ is equivariant and the action of $Aut_{\cal I}(M)$ on $[[\phi(\vec{x})]]_{M}$ is transitive. The sequent $(\theta \vdash_{\vec{x}, \vec{y}} \phi(\vec{x}) \wedge \psi(\vec{y}))$ is provable in $\mathbb T$ since it is satisfied in $M$, $\mathbb T$ is complete and $\theta$ is $\mathbb T$-complete; a similar argument proves that the sequent $(\phi(\vec{x}) \vdash_{\vec{x}} \exists \vec{y}\theta(\vec{x}, \vec{y}))$ is provable in $\mathbb T$. It remains to show that the sequent $(\theta \wedge \theta[\vec{z}\slash \vec{y}] \vdash_{\vec{x}, \vec{y}, \vec{z}} \vec{y}=\vec{z})$ is provable in $\mathbb T$. Since $\mathbb T$ is complete, it suffices to verify that this sequent is satisfied in $M$; suppose therefore that $(\vec{a}, \vec{b}), (\vec{a}, \vec{c})\in [[\theta]]_{M}$; since $\theta$ is $\mathbb T$-complete then, as we observed above, there exists an automorphism $f:M\cong M$ such that $f$ sends $(\vec{a}, \vec{b})$ to $(\vec{a}, \vec{b'})$, in other words such that $f(\vec{a})=\vec{a}$ and $f(\vec{b})=\vec{b'}$; but $\gamma$ is equivariant, which implies that $\vec{b'}=f(\vec{b})=f(\gamma(\vec{a}))=\gamma(f(\vec{a}))=\gamma(\vec{a})=\vec{b}$, as required. This completes the proof.     
\end{proofs}

The proof of the theorem motivates the following definition: given an atomic (and complete) theory $\mathbb T$, a set-based model $M$ of $\mathbb T$ is said to be \emph{special} if every $\mathbb T$-complete formula is realized in $M$ and for any $\vec{a}, \vec{b}\in M$ which satisfy exactly the same geometric (equivalently, the same $\mathbb T$-complete) formulae over the signature of $\mathbb T$, there exists an automorphism $f:M\cong M$ of $M$ sending $\vec{a}$ to $\vec{b}$. The theorem thus says that for any special model $M$ of an atomic complete theory $\mathbb T$, we have a representation $\Sh({\cal C}_{\mathbb T}, J_{\mathbb T})\simeq \Cont(Aut(M))$ of its classifying topos (where $Aut(M)$ is endowed with the topology of pointwise convergence).  

\begin{remarks}\label{rem}
\begin{enumerate}[(a)]
\item We showed in \cite{OC5} that if $\mathbb T$ is a complete atomic theory then $\mathbb T$ is countably categorical (i.e., any two countable models of it are isomorphic), and its unique (up to isomorphism) countable model $M$, if it exists, satisfies the hypotheses of the theorem. We thus obtain, as a corollary of the theorem, that the classifying topos of $\mathbb T$ can be represented as the topos $\Cont(Aut(M))$ of continuous $Aut(M)$-sets, where $Aut(M)$ is the group of automorphisms of $M$ endowed with the topology of pointwise convergence.

\item The group $Aut_{\cal I}(M)$ is complete (in the sense of section \ref{completegroups}). Indeed, by the equivalence of the theorem, the topos $\Cont(Aut_{\cal I}(M))$ is the classifying topos for $\mathbb T$, with the canonical point 
\[
p_{Aut_{\cal I}(M)}:\Set \to \Cont(Aut_{\cal I}(M))
\]
corresponding exactly to the $\mathbb T$-model $M$; the canonical map $Aut_{\cal I}(M) \to Aut(p_{Aut_{\cal I}(M)})$ is thus a bijection.

\item Given two Morita-equivalent theories $\mathbb T$ and ${\mathbb T}'$, $\mathbb T$ is atomic and complete (resp., is atomic and complete and has a special model) if and only if ${\mathbb T}'$ is atomic and complete (resp., is atomic and complete and has a special model). Indeed, the property of atomicity and completeness of a theory can be expressed as the topos-theoretic invariant property of its classifying topos to be atomic and two-valued, and hence it is stable under Morita-equivalence. On the other hand, by Theorem \ref{main} and Remark \ref{rem}(b), if an atomic complete theory $\mathbb T$ has a special model then there exists a point $p$ of its classifying topos ${\cal E}_{\mathbb T}:=\Sh({\cal C}_{\mathbb T}, J_{\mathbb T})$ such that the latter can be represented as the topos of continuous $Aut(p)$-sets where $Aut(p)$ is the group of automorphisms of $p$ endowed with the topology in which the subgroups of the form $U_{e, x}:\{\alpha:p\cong p \textrm{ | } \alpha(e)(x)=x\}$ for $e\in {\cal E}_{\mathbb T}$ and $x\in p^{\ast}(e)$ form a basis of open neighbourhoods of the identity; also, the equivalence ${\cal E}_{\mathbb T}\simeq \Cont(Aut(p))$ can be described as sending any object $e$ of ${\cal E}_{\mathbb T}$ to the canonical continuous action of $Aut(p)$ on the set $p^{\ast}(e)$. If $\mathbb T$ is Morita-equivalent to a geometric theory ${\mathbb T}'$ then we have an equivalence of classifying toposes ${\cal E}_{\mathbb T}\simeq {\cal E}_{{\mathbb T}'}$ which sends the point $p$ of ${\cal E}_{\mathbb T}\simeq {\cal E}_{{\mathbb T}'}$ to a point $q$ of ${\cal E}_{{\mathbb T}'}$; hence, denoted by $Aut(q)$ is the group of automorphisms of $q$, endowed with the topology in which the subgroups of the form $U'_{e', y}:\{\beta:q\cong q \textrm{ | } \beta(e')(y)=y\}$ for $e'\in {\cal E}_{{\mathbb T}'}$ and $y\in q^{\ast}(e')$ form a basis of open neighbourhoods of the identity, we have a homeomorphism $\tau:Aut(p)\cong Aut(q)$ such that the equivalence of classifying toposes ${\cal E}_{\mathbb T}\simeq {\cal E}_{{\mathbb T}'}$, composed with the equivalences ${\cal E}_{\mathbb T}\simeq \Cont(Aut(p))$  and $\Cont(\tau):\Cont(Aut(p))\simeq \Cont(Aut(q))$, yields an equivalence ${\cal E}_{{\mathbb T}'}\simeq \Cont(Aut(q))$ which sends any object $e'$ of ${\cal E}_{{\mathbb T}'}$ to the canonical continuous action of $Aut(q)$ on the set $q^{\ast}(e')$. If we denote by $N$ the ${\mathbb T}'$-model in $\Set$ corresponding to the point $q$ of the classifying topos ${\cal E}_{{\mathbb T}'}$ of ${\mathbb T}'$, we thus obtain an equivalence $\Sh({\cal C}_{{\mathbb T}'}, J_{{\mathbb T}'})\simeq \Cont(Aut(N))$ (where $Aut(N)$ is endowed with the topology of pointwise convergence) which sends any ${\mathbb T}'$-complete formula $\phi(\vec{x})$ to the canonical action of $Aut(N)$ on its interpretation $[[\phi(\vec{x})]]_{N}$ in $N$. As the notion of atom is a topos-theoretic invariant, this action is necessarily non-empty and transitive; in other words, the model $N$ is special for ${\mathbb T}'$, as required.      
\end{enumerate} 
\end{remarks}

\subsection{Categorical formulation}

Thanks to the results obtained by the author on the class of theories of presheaf type (cf. \cite{OC2}), we can recast our main representation theorem in purely categorical language.

To this end, we have to recall from \cite{OC2} some natural categorical notions.

\begin{definition}
Let $\cal C$ be a small category.
\begin{itemize}
\item $\cal C$ is said to satisfy the \emph{amalgamation property} (AP) if for every objects $a,b,c\in {\cal C}$ and morphisms $f:a\rightarrow b$, $g:a\rightarrow c$ in $\cal C$ there exists an object $d\in \cal C$ and morphisms $f':b\rightarrow d$, $g':c\rightarrow d$ in $\cal C$ such that $f'\circ f=g'\circ g$:
\[  
\xymatrix {
a \ar[d]_{g} \ar[r]^{f} & b  \ar@{-->}[d]^{f'} \\
c \ar@{-->}[r]_{g'} & d } 
\] 

\item $\cal C$ is said to satisfy the \emph{joint embedding property} (JEP) if for every pair of objects $a,b\in {\cal C}$ there exists an object $c\in \cal C$ and morphisms $f:a\rightarrow c$, $g:b\rightarrow c$ in $\cal C$:
\[  
\xymatrix {
 & a  \ar@{-->}[d]^{f} \\
b \ar@{-->}[r]_{g} & c } 
\] 

\item Given a full embedding of categories ${\cal C}\hookrightarrow {\cal D}$, an object $u$ of $\cal D$ is said to be \emph{$\cal C$-homogeneous} if for every objects $a,b \in {\cal C}$ and arrows $j:a\rightarrow b$ in ${\cal C}$ and $\chi:a\rightarrow u$ in ${\cal D}$ there exists an arrow $\tilde{\chi}:b\rightarrow u$ such that $\tilde{\chi}\circ j=\chi$:
\[  
\xymatrix {
a \ar[d]_{j} \ar[r]^{\chi} & u \\
b \ar@{-->}[ur]_{\tilde{\chi}}; &  } 
\] 
\item Given a full embedding of categories ${\cal C}\hookrightarrow {\cal D}$, an object $u$ of $\cal D$ is said to be \emph{$\cal C$-ultrahomogeneous} if for every objects $a,b \in {\cal C}$ and arrows $j:a\rightarrow b$ in ${\cal C}$ and $\chi_{1}:a\rightarrow u$, $\chi_{2}:b\rightarrow u$ in $\cal D$ there exists an isomorphism $\check{j}:u\rightarrow u$ such that $\check{j}\circ \chi_{1}=\chi_{2}\circ j$:
\[  
\xymatrix {
a \ar[d]_{j} \ar[r]^{\chi_{1}} & u \ar@{-->}[d]^{\check{j}}\\
b \ar[r]_{\chi_{2}} & u } 
\] 
\item Given a full embedding of categories ${\cal C}\hookrightarrow {\cal D}$, an object $u$ of $\cal D$ is said to be \emph{$\cal C$-universal} if it is $\cal C$-cofinal, that is for every $a\in {\cal C}$ there exists an arrow $\chi:a\rightarrow u$ in $\cal D$:
\[  
\xymatrix {
a \ar@{-->}[r]^{\chi} & u.  } 
\]     
\end{itemize}
\end{definition}

\begin{remarks}\label{ultrahom}
\begin{enumerate}[(a)]
\item In the above definition of $\cal C$-ultrahomogeneous object, one can suppose $j=1_{a}$ without loss of generality; 

\item Any object which is both $\cal C$-ultrahomogeneous and $\cal C$-universal is $\cal C$-homogeneous; in particular, any $\cal C$-ultrahomogeneous and $\cal C$-universal object of the ind-completion $\Ind{\cal C}$ of $\cal C$ can be identified with an object of the subcategory $\mathbf{Flat}_{J_{at}}({\cal C}^{\textrm{op}}, \Set)$ of $\Ind{\cal C}$. 
\end{enumerate}

\end{remarks}

\begin{theorem}\label{maincategorical}
Let $\cal C$ be a small non-empty category satisfying AP and JEP, and let $u$ be a ${\cal C}$-universal and ${\cal C}$-ultrahomogeneous object in $\Ind{\cal C}$. Then the collection ${\cal I}_{\cal C}$ of sets of the form ${\cal I}_{\chi}:=\{f:u\cong u \textrm{ | } f\circ \chi=\chi\}$, for an arrow $\chi:c\to u$ from an object $c$ of $\cal C$ to $u$ defines an algebraic base for the group of automorphisms of $u$ in $\Ind{\cal C}$, and, denoted by $Aut_{\cal C}(u)$ the resulting topological group, we have an equivalence of toposes
\[
\Sh({\cal C}^{\textrm{op}}, J_{at})\simeq \Cont(Aut_{\cal C}(u))
\] 
induced by the functor $F:{\cal C}^{\textrm{op}}\to \Cont_{t}(Aut_{\cal C}(u))$ which sends any object $c$ of $\cal C$ to the set $Hom_{\Ind{\cal C}}(c, u)$ (equipped with the obvious action by $Aut_{\cal C}(u)$) and any arrow $f:c\to d$ in $\cal C$ to the $Aut_{\cal C}(u)$-equivariant map $-\circ f:Hom_{\Ind{\cal C}}(d, u)\to Hom_{\Ind{\cal C}}(c, u)$.
\end{theorem}

\begin{proofs}
Let $\mathbb T$ be the geometric theory of flat functors on ${\cal C}^{\textrm{op}}$, and ${\mathbb T}'$ be its quotient axiomatizing the $J_{at}$-continuous flat functors on ${\cal C}^{\textrm{op}}$. Then ${\mathbb T}'$ is classified by the topos $\Sh({\cal C}^{\textrm{op}}, J_{at})$, which is atomic (since $\cal C$ satisfies the amalgamation property) and two-valued (since $\cal C$ satisfies the joint embedding property, cf. \cite{OC2}). Therefore $\mathbb T$ is atomic and complete. Every object of $\cal C$ can be seen as a finitely presentable $\mathbb T$-model, and, as proved in \cite{OCS}, the formulae which present these objects are ${\mathbb T}'$-complete and hence satisfy the hypotheses of Theorem \ref{main} since the objects of the form $l(c)$, where $l$ is the composite of the Yoneda embedding ${\cal C}^{\textrm{op}}\to [{\cal C}, \Set]$ with the associated sheaf functor $a_{J_{at}}:[{\cal C}, \Set]\to \Sh({\cal C}^{\textrm{op}}, J_{at})$, define a separating set for the topos $\Sh({\cal C}^{\textrm{op}}, J_{at})$ and hence are dense in the family of its atoms. To conclude, it remains to observe that, for any $c\in {\cal C}$, if $\phi_{c}$ is the geometric formula which presents the object $c$ then the interpretation of $\phi_{c}$ in $u$, regarded as a $\mathbb T$-model, corresponds precisely to the set of arrows from $c$ to $u$ in $\Ind{\cal C}$. 
\end{proofs}

\begin{remarks}
\begin{enumerate}[(a)]
\item Theorem \ref{main} can be deduced from its categorical counterpart (Theorem \ref{maincategorical}) by taking $\cal C$ equal to the category ${\cal C}_{\mathbb T}^{\cal I}$ and $u$ to be object of $\Ind{\cal C}=\mathbf{Flat}({\cal C}_{\mathbb T}^{\cal I}, \Set)$ corresponding to the object $M$ of the category ${\mathbb T}\textrm{-mod}(\Set)$ under the Morita-equivalence ${\mathbb T}\textrm{-mod}(\Set)\simeq \mathbf{Flat}_{J_{at}}({\cal C}_{\mathbb T}^{\cal I}, \Set)$; notice that the category ${\cal C}_{\mathbb T}^{\cal I}$ satisfies the amalgamation property (by our denseness hypothesis on $\cal I$) and the joint embedding property (since it is complete and classified by the atomic topos $\Sh({\cal C}_{\mathbb T}^{\cal I}, J_{at})$, cf. \cite{OC2}) and $u$, regarded as a flat functor $F:{\cal C}_{\mathbb T}^{\cal I} \to \Set$ satisfies the property that for any $\phi(\vec{x})$ in $\cal I$, $F(\phi(\vec{x}))=[[\phi(\vec{x})]]_{M}$, which implies that $Hom_{\Ind{\cal C}}(\phi(\vec{x}), u)$ is isomorphic to $[[\phi(\vec{x})]]_{M}$.
  
\item All the hypotheses of Theorem \ref{maincategorical} are necessary. Indeed, the amalgamation property on $\cal C$ is necessary for defining the atomic topology on the opposite of it, while the joint embedding property is necessary because it is equivalent, under the assumption that $\cal C$ is non-empty, to the condition that the topos $\Sh({\cal C}^{\textrm{op}}, J_{at})$ is two-valued, which is always satisfied if the latter is equivalent to a topos of continuous actions of a topological group; the condition that $u$ should be $\cal C$-universal and $\cal C$-ultrahomogeneous is also necessary since it is equivalent to the requirement that the functor $F$ should take values in the subcategory of non-empty transitive actions, which is always the case if $F$ induces an equivalence of toposes $\Sh({\cal C}^{\textrm{op}}, J_{at})\simeq \Cont(Aut_{\cal C}(u))$ (the composite of the Yoneda embedding ${\cal C}^{\textrm{op}}\to [{\cal C}, \Set]$ with the associated sheaf functor $a_{J_{at}}:[{\cal C}, \Set]\to \Sh({\cal C}^{\textrm{op}}, J_{at})$ takes values in the full subcategory of $\Sh({\cal C}^{\textrm{op}}, J_{at})$ on the atoms, and any equivalence sends atoms to atoms). It is also worth to note that if all arrows in $\cal C$ are monomorphisms, the existence of a $\cal C$-universal and $\cal C$-homogeneous object in $\Ind{\cal C}$ implies the fact that $\cal C$ satisfies the amalgamation property. 
\end{enumerate}
\end{remarks}

A particularly natural context in which Theorems \ref{main} and \ref{maincategorical} can be applied is that of theories of presheaf type. 

\begin{corollary}\label{corollarypresheaf}
Let $\mathbb T$ be a theory of presheaf type such that its category $\textrm{f.p.}{\mathbb T}\textrm{-mod}(\Set)$ of finitely presentable models satisfies AP and JEP, and let $M$ be a $\textrm{f.p.}{\mathbb T}\textrm{-mod}(\Set)$-universal and $\textrm{f.p.}{\mathbb T}\textrm{-mod}(\Set)$-ultrahomogeneous model of $\mathbb T$. Then we have an equivalence of toposes 
\[
\Sh({\textrm{f.p.}{\mathbb T}\textrm{-mod}(\Set)}^{\textrm{op}}, J_{at})\simeq \Cont(Aut(M)),
\]
where $Aut(M)$ is endowed with the topology of pointwise convergence.

Let $\phi(\vec{x})$ and $\psi(\vec{y})$ be formulae presenting respectively $\mathbb T$-models $M_{\phi}$ and $M_{\psi}$, and let $\vec{a}$ and $\vec{b}$ elements of $[[\phi(\vec{x})]]_{M}$ and $[[\psi(\vec{y})]]_{M}$. If every automorphism of $f$ which fixes $\vec{b}$ fixes $\vec{a}$ then there exists a unique ${\mathbb T}'$-provably functional formula $\theta(\vec{x}, \vec{y})$ from $\phi(\vec{x})$ to $\psi(\vec{y})$ such that $[[\theta]]_{M}(\vec{a})=\vec{b}$. If moreover all the arrows in $\textrm{f.p.}{\mathbb T}\textrm{-mod}(\Set)$ are strict monomorphisms, the formula $\theta$ can be taken to be provably functional in $\mathbb T$; in other words it induces a $\mathbb T$-model homomorphism $z:M_{\psi}\to M_{\phi}$ such that, denoted by $a:M_{\phi}\to M$ the $\mathbb T$-model homomorphism corresponding to the element $\vec{a}$ and by $b:M_{\psi}\to M$ the $\mathbb T$-model homomorphism corresponding to the element $\vec{b}$, we have $b\circ z=a$.     
\end{corollary}\qed

\begin{remarks}
\begin{enumerate}[(a)]
\item The topos $\Sh({\textrm{f.p.}{\mathbb T}\textrm{-mod}(\Set)}^{\textrm{op}}, J_{at})$ classifies (in every Grothendieck topos) the homogeneous $\mathbb T$-models (in the sense of \cite{OC2}); in particular, it classifies in $\Set$ the ${\textrm{f.p.}{\mathbb T}\textrm{-mod}(\Set)}$-homogeneous objects of the category ${\mathbb T}\textrm{-mod}(\Set)$;

\item For any theory of presheaf type $\mathbb T$ such that its category $\textrm{f.p.}{\mathbb T}\textrm{-mod}(\Set)$ of finitely presentable models satisfies AP, the joint embedding property on $\textrm{f.p.}{\mathbb T}\textrm{-mod}(\Set)$ can always be achieved at the cost of replacing $\mathbb T$ with any of its `connected components' (that is, with any of the quotients of $\mathbb T$ axiomatized by the geometric sequents which are valid in all the models belonging to a given connected component of the category $\textrm{f.p.}{\mathbb T}\textrm{-mod}(\Set)$). In fact, the completions of the geometric theory of homogeneous $\mathbb T$-models correspond exactly to the connected components of the category $\textrm{f.p.}{\mathbb T}\textrm{-mod}(\Set)$ (cf. \cite{OC2}).  
\end{enumerate}
\end{remarks}

\subsection{Prodiscreteness}

In this section we investigate some conditions on $({\cal C}, u)$ under which the topological group $Aut_{\cal C}(u)$ is prodiscrete. Recall that a topological group is prodiscrete if it is an inverse limit of discrete groups in the category of topological groups.

\begin{definition}
Given a pair $({\cal C}, u)$ satisfying the hypotheses of Theorem \ref{maincategorical}, a pair $(c, f)$ where $c\in {\cal C}$ and $f$ is an arrow $c\to u$ in $\Ind{\cal C}$ is said to be a $({\cal C}, u)$\emph{-Galois object} if for any automorphism $\xi:u\cong u$ in $\Ind{\cal C}$ there exists exactly one automorphism $s:c\cong c$ such that $\xi \circ f=f\circ s$:
\[  
\xymatrix {
c \ar@{-->}[d]_{s} \ar[r]^{f} & u \ar[d]^{\xi}\\
c \ar[r]_{f} & u } 
\] 
\end{definition}

Given $({\cal C}, u)$, let $\cal F$ be the category whose objects are the pairs $(c, f)$, where $c$ is an object of $\cal C$ and $f:c\to u$ is an arrow in $\Ind{\cal C}$, and whose arrows $(c, f)\to (d, g)$ are the arrows $l:c\to d$ in $\cal C$ such that $g\circ l=f$ in $\Ind{\cal C}$. Let ${\cal F}_{g}$ be the subcategory of $\cal F$ whose objects are the $({\cal C}, u)$-Galois objects and whose arrows $(c, f)\to (d, g)$ are the arrows $l:c\to d$ in $\cal C$ such that for any automorphism $\chi$ on $d$ in $\cal C$ there exists exactly one automorphism $t:c\cong c$ of $c$ in $\cal C$ such that $\chi\circ l=l \circ t$. Notice that the category ${\cal F}_{g}$ is closed under isomorphisms in ${\cal F}$. Note also that if all the arrows in $\cal C$ are monic and for any $({\cal C}, u)$-Galois object $(c, f)$ the arrow $f:c\to u$ is monic in $\Ind{\cal C}$ then the category ${\cal F}_{g}$ coincides with the full subcategory of $\cal F$ on the $({\cal C}, u)$-Galois objects (since $u$ is $\cal C$-ultrahomogeneous); of course, these conditions are satisfied if any arrow of $\Ind{\cal C}$ is monic.

Recall that a subcategory $\cal C$ of a category $\cal D$ is said to be \emph{cofinal} in $\cal D$ if for every $d\in {\cal D}$ the comma category $d\slash {\cal C}$ is non-empty and connected. If $\cal C$ is cofinal in $\cal D$ then for any functor $F:{\cal D}\to \Set$, the colimit of $F$ is isomorphic to the colimit of $F\circ i$, where $i$ is the inclusion ${\cal C}\hookrightarrow D$. Given a family $\cal F$ of objects of $\cal D$ we say that $\cal F$ is cofinal in $\cal D$ if every object of $\cal D$ admits an arrow to an object in $\cal F$.   

\begin{lemma}\label{cofinality}
Let $\cal C$ be a full subcategory of a filtered category $\cal D$. Then, if for any object $d$ of $\cal D$ there exists an arrow $d\to c$ in $\cal C$ to an object $c$ of $\cal D$, $\cal C$ is filtered and cofinal in $\cal D$.
\end{lemma}

\begin{proofs}
To prove that $\cal C$ is filtered we observe that $\cal C$ is non-empty since $\cal D$ is non-empty and for any object $d$ of $\cal D$ there exists an arrow $d\to c$ in $\cal C$ to an object $c$ of $\cal C$. The fact that $\cal C$ satisfies the joint embedding property and the weak coequalizer property follows at once from the fact that $\cal C$ is full in $\cal D$ as these properties are satisfied by $\cal D$ by our hypotheses. It remains to prove that for any $d\in {\cal D}$ the comma category $d\slash {\cal C}$ is connected. We shall prove that it satisfies the joint embedding property, from which our thesis will clearly follow. Given any two objects $f:d\to c$ and $g:d\to c'$ in $d\slash {\cal C}$, since $\cal D$ satisfies the joint embedding property, there exist an object $d'$ and two arrows $h:c\to d'$ and $k:c'\to d'$ in $\cal D$; now, the fact that $\cal D$ satisfies the weak coequalizer property and the fact that there exists an arrow $d'\to c''$ where $c''$ is an object of $\cal C$ imply that we can suppose without loss of generality that $h\circ f=k\circ g$ and that $d'\in {\cal C}$; the arrow $h\circ f=k\circ g:d\to d''$ thus defines an object of $(d\slash {\cal C})$ (as $\cal C$ is full in $\cal D$) and the arrows $h$ and $k$ define respectively arrows $f\to h\circ f$ and $g\to g\circ k$ in $d\slash {\cal C}$.      
\end{proofs}

\begin{theorem}
Let $({\cal C}, u)$ be a pair satisfying the hypotheses of Theorem \ref{maincategorical}. Then, if every arrow in $\Ind{\cal C}$ is monic and the $({\cal C}, u)$-Galois objects are cofinal in $\cal F$, the topological group $Aut_{\cal C}(u)$ is prodiscrete; specifically, it can be expressed as the pro-limit of the diagram $D:{{\cal F}_{g}}^{\textrm{op}}\to \mathbf{TopGr}$ sending any $(c, f)$ in ${\cal F}_{g}$ to the discrete automorphism group $D((c, f))=Aut_{\cal C}(c)$ and any arrow $l:(c, f)\to (d, g)$ in ${\cal F}_{g}$ to the group homomorphism $D(l):Aut_{\cal C}(d)\to Aut_{\cal C}(c)$ sending any automorphism $\chi$ of $d$ to the unique automorphism $t$ of $c$ such that $\chi\circ l=l\circ t$, the colimit arrows (for $(c, f)\in {\cal F}_{g}$) being precisely the maps $Aut_{\cal C}(u)\to Aut(c)$ sending to any automorphism $\xi$ of $u$ the unique automorphism $s$ of $c$ such that $\xi \circ f=f\circ s$.  
\end{theorem}  

\begin{proofs}
Recall from \cite{DD} that the (projective) limit $L_{D}$ of a cofiltered diagram $D:{\cal I}\to \mathbf{TopGr}$ of topological groups can be described as follows: $L_{D}$ is the subgroup of the direct product $\prod_{i\in {\cal I}}D(i)$ consisting of the strings $(x_{i})_{i\in {\cal I}}$ such that for any arrow $f:i\to j$ in $\cal I$, $D(f)(x_{i})=x_{j}$, with the topology of a subspace of the product $\prod_{i\in {\cal I}}D(i)$ of the discrete topological spaces $D_{i}$. In fact, this topology can alternatively be described as the coarsest topology making all the natural projection maps $L_{D}\to D(i)$ (for $i\in {\cal I}$) continuous.

Since all the arrows in $\Ind{\cal C}$ are monic, the category ${\cal F}_{g}$ coincides with the full subcategory of $\cal F$ on the $({\cal C}, u)$-Galois objects. Therefore, by our hypotheses and Lemma \ref{cofinality}, the category ${\cal F}_{g}$ is filtered and hence we can consider the prolimit $L_{D}$ of the diagram $D:{{\cal F}_{g}}^{\textrm{op}}\to \mathbf{TopGr}$ defined in the statement of the theorem. 

We can establish an isomorphism between $Aut_{\cal C}(u)$ and $L_{D}$ as follows. We assign to any element $\xi\in Aut_{\cal C}(u)$ the element $(x_{(c, f)})_{(c, f)\in {\cal F}_{g}}$ of $L(D)$, where $x_{(c, f)}:c\cong c$ is the unique automorphism $s$ of $c$ (that is, element $s$ of $Aut_{\cal C}(c)$) such that $\xi \circ f=f\circ s$; this is well-defined since for any arrow $l:(c, f)\to (d, g)$ in ${\cal F}_{g}$, $x_{(c, f)}$ coincides with $D(l)((d, g))$ since both arrows are factorizations of $\xi$ across $g$. To show that this map $a:Aut_{\cal C}(u) \to L_{D}$ is an isomorphism of groups, we shall construct an inverse $b:L_{D}\to Aut_{\cal C}(u)$. We define $b$ as follows: given $(x_{(c, f)})_{(c, f)\in {\cal F}_{g}}$ in $L_{D}$, we define $b((x_{(c, f)})_{(c, f)\in {\cal F}_{g}})$ equal to the unique arrow $h:u\to u$ such that for any $(c, f)$ in ${\cal F}_{g}$, $h\circ f=f\circ x_{(c, f)}$, using the fact that, by Lemma \ref{cofinality}, $u$ is the colimit of the forgetful functor ${\cal F}_{g} \to \Ind{\cal C}$ and the arrows $f\circ x_{(c, f)}:c\to u$ (for $(c, f)\in {\cal F}_{g}$) define a cocone on $u$ over it; $h$ is an automorphism of $u$ because, symmetrically, using the inverses $x_{(c, f)}^{-1}$ of the $x_{(c, f)}$, one can obtain an arrow $k:u\to u$ which is easily verified to be a two-sided inverse of $h$. It is now clear that the arrows $a$ and $b$ are inverse to each other, and that the colimiting cone to $D$ arising from the composition of the canonical maps $L(D)\to D((c, f))$ (for $(c, f)\in {\cal F}_{g}$) with the arrow $a$ is precisely given by the maps $Aut_{\cal C}(u)\to Aut(c)$ sending to any automorphism $\xi$ of $u$ the unique automorphism $s$ of $c$ such that $\xi \circ f=f\circ s$ (for $(c, f)\in {\cal F}_{g}$). This shows that the map $a$ is an isomorphism of groups, so it remains to verify that $a:Aut_{\cal C}(u) \to L_{D}$ is actually an homeomorphism; for this, it suffices to check that the topology on $Aut_{\cal C}(u)$ is the coarsest one making all the canonical maps $Aut_{\cal C}(u)\to Aut_{\cal C}(c)$ continuous. For any object $(c, f)$ of ${\cal F}_{g}$, the canonical map $a_{(c, f)}:=Aut_{\cal C}(u)\to Aut_{\cal C}(c)$ corresponding to it sends any automorphism $\xi$ of $u$ to the unique automorphism $s:c\cong c$ such that $\xi \circ f=f\circ s$; therefore the inverse image $a_{(c,f)}^{-1}(1_{c})$ of the identity element of $Aut_{\cal C}(c)$ via this map is equal to the set ${\cal I}_{f}$ of all automorphisms $\xi$ of $u$ such that $\xi \circ f=f$. Now, if $\chi:G\to H$ is a surjective homomorphism of groups then for any $h\in H$ there exists $g\in G$ such that $\chi^{-1}(\{h\})=g\chi^{-1}(\{e\})$ (where $e$ is the neutral element of $H$); indeed, if $\chi(g)=h$ (notice that such a $g$ exists since $\chi$ is surjective) then for any $x\in G$, $x\in \chi^{-1}(\{h\})$, i.e. $\chi(x)=h$ if and only if $\chi(x)\chi(g^{-1})=\chi(xg^{-1})=e$, that is if and only if $x\in g\chi^{-1}(\{e\})$. This remark, applied to our situation (notice that $a_{(c, f)}$ is surjective as $u$ is $\cal C$-ultrahomogeneous) shows that the inverse image of any singleton set in $Aut_{\cal C}(u)$ under $a_{(c, f)}$ is the translation of an open set of the form ${\cal I}_{f}$ by an element of $Aut_{\cal C}(u)$ and hence it is open in $Aut_{\cal C}(u)$. Hence, the topology on $Aut_{\cal C}(c)$ being discrete, we can conclude that $a_{(c, f)}$ is continuous. Now, the fact that the topology on $Aut_{\cal C}(u)$ is the coarsest to make all these maps continuous follows from the fact that the ${\cal I}_{f}$ (for $(c, f)$ in ${\cal F}_{c}$) form an algebraic base for $Aut_{\cal C}(u)$ (by our hypothesis, the $({\cal C}, u)$-Galois objects are cofinal in $\cal F$ and hence Remark \ref{algbases} applies).      
\end{proofs}

\subsection{Functorialization}\label{functorialization}

We can functorialize the Morita-equivalence established in Theorem \ref{maincategorical} by using the notion of geometric morphism as a `bridge' across the two different representations.

We observe that, for any categories $\cal C$ and ${\cal C}'$ satisfying the amalgamation property, any morphism of sites $F:({\cal C}^{\textrm{op}}, J_{at}) \to ({{\cal C}'}^{\textrm{op}}, J_{at})$ and hence induces a geometric morphism $\Sh(F):\Sh({{\cal C}'}^{\textrm{op}}, J_{at}) \to \Sh({\cal C}^{\textrm{op}}, J_{at})$, which in turn corresponds, via Diaconsecu' equivalence, to a morphism $\tilde{F}:\mathbf{Flat}_{J_{at}}({{\cal C}'}^{\textrm{op}}, \Set)\to \mathbf{Flat}_{J_{at}}({\cal C}^{\textrm{op}}, \Set)$ which can be identified with $-\circ F$. Recall that a morphism of sites is a flat-covering cover-preserving functor (cf. \cite{Shulman}, for instance); in particular, any cartesian cover-preserving morphism between sites whose underlying categories are cartesian is a morphism of sites. 

Let us define $\cal G$ to be the category whose objects are the pairs $({\cal C}, u)$, where $\cal C$ is a small category satisfying AP and JEP and $u$ is a $\cal C$-ultrahomogeneous and $\cal C$-universal object of $\Ind{\cal C}$, and whose arrows $({\cal C}, u)\to ({\cal C}', u')$ are the morphisms of sites $f:({\cal C}^{\textrm{op}}, J_{at}) \to ({{\cal C}'}^{\textrm{op}}, J_{at})$ such that $\tilde{F}(u')=u$ (notice that this is well-defined by Remark \ref{ultrahom}(b), as $u$ is $\cal C$-ultrahomogeneous and $\cal C$-universal). Then we have a functor $A:{\cal G}^{\textrm{op}}\to \mathbf{GTop_{b}}$ sending any pair $({\cal C}, u)$ to the object $(Aut_{\cal C}(u), {\cal I}_{\cal C})$ of $\mathbf{GTop_{b}}$ and any arrow $F:({\cal C}, u)\to ({\cal C}', u')$ to the arrow $\tilde{F}:Aut_{{\cal C}'}(u') \to Aut_{\cal C}(u)$ in $\mathbf{GTop_{b}}$:

\[  
\xymatrix@1@=3pt@M=3pt {
 & & & & & \Sh({\cal C}^{\textrm{op}}, J_{at})  & \simeq & \Cont(Aut_{\cal C}(u))  & & & & &\\
 & & & & & & &  & \ar@/^10pt/@{--}[dddrrr] &  & \\
 & & & & & \Sh({{\cal C}'}^{\textrm{op}}, J_{at}) \ar[uu]  & \simeq & \Cont(Aut_{{\cal C}'}(u')) \ar[uu]  & & & & & \\
 ({\cal C}, u) \ar[dd]^{F} &  & & & & & & & & & & & (Aut_{\cal C}(u), {\cal I}_{\cal C}) \\
  & \ar@/^10pt/@{--} [uuurrr]  & & & & & & & & & & &\\
 ({\cal C}', u')  & & & & & & & & & & & & (Aut_{{\cal C}'}(u), {\cal I}_{{\cal C}'}) \ar[uu]^{\tilde{F}} 
 }
\]   

This is well-defined since for any arrow $\chi:c\to u$ in $\Ind{\cal C}$ from an object $c$ of $\cal C$ to $u$, $\tilde{F}^{-1}({\cal I}_{\chi})={\cal I}_{\xi}$ for some arrow $\xi:F(c)\to u'$ in $\Ind{{\cal C}'}$. Indeed, if we denote by $\tilde{u}:{\cal C}^{\textrm{op}}\to \Set$ and by $\tilde{u'}:{{\cal C}'}^{\textrm{op}}\to \Set$ the flat functors corresponding to the objects $u$ and $u'$ respectively of $\Ind{\cal C}$ and of $\Ind{{\cal C}'}$, we have that $\tilde{u'}\circ F=\tilde{u}$ and hence the arrows $\chi:c\to u$ in $\Ind{\cal C}$ (i.e., by Yoneda, the elements of $\tilde{u}(c)=\tilde{u'}(F(c))$) correspond exactly to the arrows $\xi:F(c)\to u'$ in $\Ind{{\cal C}'}$ (i.e., by Yoneda, the elements of $\tilde{u'}(F(c))$); we thus have $\tilde{F}^{-1}({\cal I}_{\chi})={\cal I}_{\xi}$, where $\xi$ is the arrow $F(c) \to u'$ in $\Ind{{\cal C}'}$ associated to $\chi:c\to u$ via this correspondence. 
 
In the converse direction, we can define a functor $B:\mathbf{GTop_{b}} \to {\cal G}^{\textrm{op}}$ sending any object $(G, {\cal B})$ of $\mathbf{GTop_{b}}$ to the pair $({\Cont_{\cal B}(G)}^{\textrm{op}}, p_{G_{{\cal B}}})$ (where $p_{G_{{\cal B}}}$, namely the canonical point of the topos $\Cont(G_{\cal B})$, is regarded as an object of $\Ind{{\Cont_{\cal B}(G)}^{\textrm{op}}}$ in the canonical way) and any arrow $f:(G, {\cal B})\to (G', {\cal B}')$ in $\mathbf{GTop_{b}}$ to the arrow ${\Cont(f)^{\ast}}^{\textrm{op}}|:{\Cont_{{\cal B}'}(G')}^{\textrm{op}}\to {\Cont_{\cal B}(G)}^{\textrm{op}}$ (notice that this restriction is indeed well-defined since by our hypotheses the inverse image under $f$ of any open subgroup of $G'$ belonging to ${\cal B}'$ is an open subgroup of $G$ belonging to $\cal B$). To prove that this functor is well-defined we observe that by Proposition \ref{propAPJEP} for any object $(G, {\cal B})$ of $\mathbf{GTop_{b}}$ the category $B((G, {\cal B}))$ satisfies the amalgamation and joint embedding properties; on the other hand, $p_{G_{{\cal B}}}$ is a ${\Cont_{\cal B}(G)}^{\textrm{op}}$-ultrahomogeneous and ${\Cont_{\cal B}(G)}^{\textrm{op}}$-universal object, since for any object $c$ of ${\Cont_{\cal B}(G)}^{\textrm{op}}$,
\[
Hom_{\Ind{{\Cont_{\cal B}(G)}^{\textrm{op}}}}(c, p_{G_{{\cal B}}})\cong p_{G_{{\cal B}}}(c)\cong c,
\]
which is a non-empty $Aut({p_{G_{{\cal B}}}})$-transitive set (since it is a transitive $G_{\cal B}$-set and for any element $g$ of $G_{\cal B}$, the action of $g$ on $c$ coincides with the component at $c$ of the action on $p_{G_{{\cal B}}}$ of the image of $g$ under the canonical map $G_{\cal B}\to Aut({p_{G_{{\cal B}}}})$.       
The fact that for any arrow $f:(G, {\cal B})\to (G', {\cal B}')$ in $\mathbf{GTop_{b}}$, $B(f):B((G', {\cal B}')) \to B((G, {\cal B}))$ is an arrow in $\cal G$ is immediate from the fact that $B(f)$ is the restriction of the inverse image of a geometric morphism to subcanonical sites.     

We can visualize this as follows:

\[  
\xymatrix@1@=3pt@M=3pt {
 & & & & & \Sh(\Cont_{{\cal B}'}(G'), J_{at})  & \simeq & \Cont(G_{{\cal B}'})  & & & & &\\
 & & & & & & &  & \ar@/^10pt/@{--}[dddrrr] &  & \\
 & & & & & \Sh(\Cont_{\cal B}(G), J_{at}) \ar[uu]  & \simeq & \Cont(G_{\cal B}) \ar[uu]  & & & & & \\
 ({\Cont_{{\cal B}'}(G')}^{\textrm{op}}, p_{G'_{{\cal B}'}}) \ar[dd]^{{\Cont(f)^{\ast}}^{\textrm{op}}|} &  & & & & & & & & & & & (G', {\cal B}') \\
  & \ar@/^10pt/@{--} [uuurrr]  & & & & & & & & & & &\\
 ({\Cont_{\cal B}(G)}^{\textrm{op}}, p_{G_{{\cal B}}})  & & & & & & & & & & & & (G, {\cal B}) \ar[uu]^{f} 
 }
\]

We have the following result. 

\begin{theorem}\label{adjunction}
The functors 
\[
A:{\cal G}^{\textrm{op}}\to \mathbf{GTop_{b}}
\]
and
\[
B:\mathbf{GTop_{b}} \to {\cal G}^{\textrm{op}}
\]
defined above are adjoint to each other ($A$ on the right and $B$ on the left) and restrict to a duality between the full subcategory ${\cal G}^{\textrm{op}}_{sm}$ of ${\cal G}^{\textrm{op}}$ on the objects $({\cal C}, u)$ such that every morphism in $\cal C$ is a strict monomorphism and the full subcategory $\mathbf{GTop_{b}}^{c}$ of $\mathbf{GTop_{b}}$ on the objects $(G, {\cal B})$ such that the topological group $G_{\cal B}$ is complete.    
\end{theorem}

\begin{proofs}
The counit $\epsilon$ of the adjunction is given, for any $({\cal C}, u)$ in $\cal G$, by $\epsilon({\cal C}, u)=F^{\textrm{op}}:=({\cal C}, u)\to ({\Cont_{{\cal I}_{\cal C}}(Aut_{\cal C}(u))}^{\textrm{op}}, p_{{Aut_{\cal C}(u)}_{{\cal I}_{\cal C}}})$, regarded as an arrow in $\cal G$, where $F$ is the functor defined in the statement of Theorem \ref{maincategorical}, while the unit $\eta$ is given, for any $(G, {\cal B})\in \mathbf{GTop_{b}}$, by $\eta((G, {\cal B}))=\xi_{G_{\cal B}}:=(G, {\cal B})\to (Aut_{{{\Cont_{\cal B}(G)}^{\textrm{op}}}}(p_{G_{{\cal B}}}), {\cal I}_{{\Cont_{\cal B}(G)}^{\textrm{op}}})$ (cf. section \ref{topgroups} for the definition of the canonical map $\xi$). The naturality of the unit and counit, and the triangular identities, are easily verified.  

Now, by Theorem \ref{subcanonical}, $\epsilon({\cal C}, u)$ is an isomorphism in $\cal G$ if and only if every arrow of $\cal C$ is a strct monomorphism, while $\eta((G, {\cal B}))$ is an isomorphism in $\mathbf{GTop_{b}}$ if and only if $G_{\cal B}$ is complete (cf. section \ref{completegroups}).   
\end{proofs}

\begin{remark}
Up to Morita-equivalence, the functors $A$ and $B$ defining the adjunction of Theorem \ref{adjunction} are inverse to each other. Indeed, for any $({\cal C}, u)$ in $\cal G$, $\Sh(B(A({\cal C}, u)), J_{at})\simeq \Sh({\cal C}^{\textrm{op}}, J_{at})$, while for any $(G, {\cal B})\in \mathbf{GTop_{b}}$, $\Cont(G_{\cal B})\simeq \Cont(G'_{{\cal B}'})$, where  $A(B((G, {\cal B}))=(G', {\cal B}')$.  
\end{remark}

The following result is an immediate consequence of Theorem \ref{adjunction}.

\begin{corollary}\label{atomint}
Let $({\cal C}, u)$ and $({\cal C}', u')$ be objects of ${\cal G}^{\textrm{op}}_{sm}$. Then a continuous group homomorphism $h:Aut_{{\cal C}'}(u')\to Aut_{\cal C}(u)$ is induced by a (unique) functor $F:{\cal C}\to {\cal C}'$ such that $\tilde{F}(u')=u$ if and only if the inverse image under $h$ of any open subgroup of the form ${\cal I}_{\chi}$ (where $\chi:c\to u$ is an arrow in $\Ind{{\cal C}}$) is of the form ${\cal I}_{\chi'}$ (where $\chi':c'\to u'$ is an arrow in $\Ind{{\cal C}'}$).  
\end{corollary}\qed

To apply Corollary \ref{atomint} in a logical context, we recall from \cite{OCthesis} that there are various natural notions of interpretations between theories; for instance, it is natural to define an interpretation between geometric theories $\mathbb T$ and ${\mathbb T}'$ as a geometric functor ${\cal C}_{\mathbb T} \to {\cal C}_{{\mathbb T}'}$, while for coherent theories there two additional natural notions of interpretations directly inspired by classical Model Theory: we can define an interpretation of a coherent theory $\mathbb T$ into a coherent theory ${\mathbb T}'$ as a coherent functor ${\cal C}^{coh}_{\mathbb T} \to {\cal C}^{coh}_{{\mathbb T}'}$, where ${\cal C}^{coh}_{\mathbb T}$ and ${\cal C}^{coh}_{\mathbb T}$ are respectively the coherent syntactic categories of $\mathbb T$ and of ${\mathbb T}'$, or alternatively as a coherent functor ${\cal P}_{\mathbb T}\to {\cal P}_{{\mathbb T}'}$, where ${\cal P}_{\mathbb T}$ and ${\cal P}_{{\mathbb T}'}$ are respectively the pretopos completions of ${\cal C}^{coh}_{\mathbb T}$ and ${\cal C}^{coh}_{\mathbb T}$ (that is, the categories of model-theoretic coherent imaginaries of $\mathbb T$ and $\mathbb T$ - notice incidentally that if the theories in question are Boolean then any first-order formula over their signature is provably equivalent in the theory to a coherent formula so that these categories coincide with the usual first-order syntactic categories or first-order categories of imaginaries arising in classical model theory). It is natural to define an \emph{atomic interpretation} of an atomic complete theory $\mathbb T$ into an atomic complete theory ${\mathbb T}'$ to be a morphism of sites $({\cal C}^{c}_{\mathbb T}, J_{at})\to ({\cal C}^{c}_{{\mathbb T}'}, J_{at})$, where ${\cal C}^{c}_{\mathbb T}$ (resp. ${\cal C}^{c}_{{\mathbb T}'}$) is the full subcategory of the geometric syntactic category ${\cal C}_{\mathbb T}$ of $\mathbb T$ (resp. ${\cal C}_{{\mathbb T}'}$ of ${\mathbb T}'$) on the $\mathbb T$-complete (resp. the ${\mathbb T}'$-complete) formulae.

Corollary \ref{atomint} can be notably applied to the pairs of the form $({\cal C}^{c}_{\mathbb T}, M)$, where $\mathbb T$ is a theory satisfying the hypotheses of Theorem \ref{main} with respect to the model $M$ (cf. Remark \ref{completeformulae}).

\begin{proposition}\label{int}
Let ${\mathbb T}$ and ${\mathbb T}'$ be atomic complete theories with special models $M$ and $M'$ respectively. Then a continuous group homomorphism $h:Aut(M')\to Aut(M)$ (where the groups $Aut(M)$ and $Aut(M')$ are endowed with the topology of pointwise convergence) is induced by an atomic interpretation of ${\mathbb T}$ into ${\mathbb T}'$ if and only if for any string $\vec{a}$ of elements of $M$ there exists a string $\vec{b}$ of elements of $M'$ such that $\{f:M'\cong M' \textrm{ | } f(\vec{b})=\vec{b}\}=\{f:M'\cong M' \textrm{ | } h(f)(\vec{a})=\vec{a}\}$.
\end{proposition}

In view of Corollary \ref{atomint}, it is natural to wonder whether we can explicitly characterize, given any two atomic complete theories $\mathbb T$ and ${\mathbb T}'$ satisfying the hypotheses of Theorem \ref{main}, the continuous homomorphisms $Aut(M')\to Aut(M)$ which are induced by an interpretation of $\mathbb T$ into ${\mathbb T}'$. We shall give an answer to these questions in section \ref{coherence} and \ref{universal} below. 

As we was in section \ref{topgroups}, the category $\mathbf{GTop}$ can be identified with a full subcategory of the category $\mathbf{GTop_{b}}$, by choosing the canonical algebraic base associated to any topological group. It is thus natural to wonder whether it is possible to characterize the objects $({\cal C}, u)$ of the category ${\cal G}$ which correspond to such objects under the adjunction of Theorem \ref{adjunction}. To this end, we remark that the objects of $\mathbf{GTop_{b}}$ of the form $(G, C_{G})$ can be characterized as the objects $(G, {\cal B})$ such that the category $\Cont_{{\cal B}}(G)$ coincides with the full subcategory of $\Cont(G_{\cal B})$ on its atoms. Therefore, the pairs $({\cal C}, u)$ of the form $B((G, C_{G}))$ for some topological group $G$ satisfy the property that every atom of $\Sh({\cal C}^{\textrm{op}}, J_{at})$ has, up to isomorphism, the form $l(c)$ for some object $c$ of $\cal C$ (where $l:{\cal C}^{\textrm{op}}\to \Sh({\cal C}^{\textrm{op}}, J_{at})$ is the composite of the Yoneda embedding $y:{\cal C}^{\textrm{op}}\to [{\cal C}, \Set]$ with the associated sheaf functor $a_{J_{at}}:[{\cal C}, \Set] \to \Sh({\cal C}^{\textrm{op}}, J_{at})$). Conversely, if $\cal C$ satisfies this condition then $A({\cal C}, u)=(Aut_{\cal C}(u), {\cal I}_{\cal C})$ is of the form $(G, C_{G})$, since the category $\Cont_{{\cal I}_{\cal C}}(Aut_{\cal C}(u))$ coincides with the full subcategory $\Cont_{{\cal I}_{\cal C}}(Aut_{\cal C}(u))$ of $\Cont(Aut_{\cal C}(u))$ on its atoms. An alternative characterization of the objects $({\cal C}, u)$ of $\cal G$ such that $A({\cal C}, u)$ is of the form $(G, C_{G})$ is the following: $A({\cal C}, u)$ is of the form $(G, C_{G})$ if and only if every open subgroup of $Aut_{\cal C}(u)$ is of the form ${\cal I}_{\chi}$ for some arrow $\chi:c\to u$.  

Summarizing, we have the following result.

\begin{proposition}
For any object $({\cal C}, u)$ of $\cal G$, the following conditions are equivalent:

\begin{enumerate}[(i)]
\item $A(({\cal C}, u))$ is, up to isomorphism in $\mathbf{GTop_{b}}$, of the form $(G, C_{G})$ for some topological group $G$;

\item Every open subgroup of $Aut_{\cal C}(u)$ is of the form ${\cal I}_{\chi}$ for some $\chi:c\to u$ in $\Ind{\cal C}$;

\item Every atom of the topos $\Sh({\cal C}^{\textrm{op}}, J_{at})$ is, up to isomorphism, of the form $l(c)$ for some object $c\in {\cal C}$.
\end{enumerate}
\end{proposition}    

This motivates the following definition: we shall say that a category $\cal C$ is \emph{atomically complete} if is satisfies AP, the atomic topology on ${\cal C}^{\textrm{op}}$ is subcanonical and every atom of the topos $\Sh({\cal C}^{\textrm{op}}, J_{at})$ is, up to isomorphism, of the form $l(c)$ for some object $c\in {\cal C}$.    

\begin{proposition}\label{dualatom}
Let $\cal C$ be a small category satisfying the dual of the amalgamation property and such that the atomic topology on it is subcanonical. Then every atom of $\Sh({\cal C}, J_{at})$ is, up to isomorphism, of the form $l(c)$ (for some $c\in {\cal C}$) if and only if for every object $c$ and equivalence relation $R$ on ${\cal C}(-,c)$ in $[{\cal C}^{\textrm{op}}, \Set]$ there exists an arrow $f:c\to d$ such that for any arrows $a,b: e\to c$, $f\circ a=f\circ b$ if and only if there exists an arrow $h$ such that $(a\circ h, b\circ h)\in R$.   
\end{proposition}

\begin{proofs}
Let $A$ be an atom of the topos $\Sh({\cal C}, J_{at})$; then, as the objects of the form $l(c)$ form a separating set for the topos $\Sh({\cal C}^{\textrm{op}}, J_{at})$, there exists an arrow (in fact, an epimorphism) $l(c)\to A$ for some object $c$ of $\cal C$. The object $A$ is therefore isomorphic to the quotient $l(c)\slash R$ of $l(c)$ by an equivalence relation $R$ on $l(c)$ in $\Sh({\cal C}^{\textrm{op}}, J_{at})$. As $J_{at}$ is subcanonical, the equivalence relations on $l(c)$ in $\Sh({\cal C}^{\textrm{op}}, J_{at})$ can be identified with the equivalence relation on ${\cal C}(-, c)$ in $[{\cal C}^{\textrm{op}}, \Set]$ which are $J_{at}$-closed as subobjects of ${\cal C}(-, c)\times {\cal C}(-, c)$. The condition that $A$ should be, up to isomorphism, of the form $l(d)$ thus amounts to the existence of an object $d\in {\cal C}$ and an arrow $f:c\to d$ in $\cal C$ such that $l(f)$ is isomorphic to the canonical projection $l(c)\to l(c)\slash R$ in $\Sh({\cal C}, J_{at})$. Notice that, since the associated sheaf functor $a_{J_{at}}:[{\cal C}^{\textrm{op}}, \Set]\to \Sh({\cal C}, J_{at})$ preserves coequalizers, this condition is equivalent to requiring the existence of an arrow $f:c\to d$ such that for any $(\xi, \chi)\in R$, $f\circ \xi=f\circ \chi$ and the canonical arrow $f\circ -:{\cal C}(-, c)\slash R\to {\cal C}(-, d)$ in $[{\cal C}^{\textrm{op}}, \Set]$ is sent by $a_{J_{at}}$ to an isomorphism.  
Now, for any elementary topos $\cal E$ and any local operator $j$ on $\cal E$, it is possible to express the condition for an arrow $h:A\to B$ in $\cal E$ to be sent by the associated sheaf functor $a_{j}:{\cal E}\to \sh_{j}({\cal E})$ to an isomorphism as the condition that the canonical monomorphism $m:A\to R$ from $A$ to the kernel pair $R$ of $f$ should be $j$-dense; indeed, denoted by $(r_{1}, r_{2}):R\to A$ the kernel pair of $f$, since $a_{j}$ preserves pullbacks and coequalizers we have that $a_{j}(f)$ is an isomorphism if and only if $a_{j}(r_{1})=a_{j}(r_{2})$, if and only if $a_{j}(m)$ is an isomorphism, i.e. $m$ is $j$-dense. Recall that a subobject $B\mono E$ in $[{\cal C}^{\textrm{op}}, \Set]$ is $J$-dense for a Grothendieck topology $J$ on $\cal C$ if and only if for every object $c\in {\cal C}$ and any element $x\in E(c)$, the sieve $\{f:d\to c \textrm{ | } E(f)(x)\in B(d) \}\in J(c)$. Applying this characterization to our particular case, we obtain that the canonical monomorphism ${\cal C}(-, c)\slash R \mono K$ in $[{\cal C}^{\textrm{op}}, \Set]$, where $K$ is the kernel pair of the arrow $f\circ -:{\cal C}(-, c)\slash R\to {\cal C}(-, d)$, is $J_{at}$-dense if and only if for every pair of arrows $(a,b)$ such that $f\circ a=f\circ b$ there exists an arrow $h$ such that $(a\circ h, b\circ h)\in R$. 
\end{proofs}

\begin{remark}\label{remeq}
An equivalence relation in a topos $[{\cal C}^{\textrm{op}}, \Set]$ on an object $E$ can be seen as a subfunctor ${\cal C}^{\textrm{op}} \to \Set$ of $E\times E$ sending to every object $c$ of $\cal C$ an equivalence relation of $E(c)$. Indeed, the concept of equivalence relation can be expressed by using geometric logic and the evaluation functors $ev_{c}:{\cal C}^{\textrm{op}}\to \Set$ (for $c\in {\cal C}$) are geometric and jointly conservative. 
\end{remark}

These results lead to the following duality theorem.

\begin{theorem}
The functors $A$ and $B$ defined above restrict to a duality between the full subcategory of $\cal G$ on the objects of the form $({\cal C}, u)$ for $\cal C$ atomically complete and the category of complete topological groups.
\end{theorem}\qed

\section{Concrete Galois theories}\label{concrete}

Now that we have established an equivalence of classifying toposes, namely that of Theorem \ref{maincategorical}, it is natural to consider it in conjunction with appropriate topos-theoretic invariants to construct `bridges' (in the sense of \cite{OC10}) for connecting the two sides to each other: 

\[  
\xymatrix {
 &   \Sh({\cal C}^{\textrm{op}}, J_{at})\simeq \Cont(Aut_{\cal C}(u)) \ar@/^12pt/@{--}[dr] &  \\
 {\cal C}^{\textrm{op}} \ar@/^12pt/@{--}[ur] & & Aut_{\cal C}(u) }
\]

\subsection{Strict monomorphisms and Galois-type equivalences}

By considering the invariant notion of arrow between two given objects in the context of the two representations above, we obtain the following `Galois-type' theorem.

\begin{theorem}\label{subcanonical}
Let $\cal C$ be a small category satisfying both the amalgamation and joint embedding properties, and let $u$ be a ${\cal C}$-universal and ${\cal C}$-ultrahomog-\\eneous object in $\Ind{\cal C}$. Then the following conditions are equivalent:
\begin{enumerate}[(i)]
\item Every arrow $f:d\to c$ in $\cal C$ is a strict monomorphism (in the sense that for any arrow $g:e\to c$ such that $h\circ g=k\circ g$ whenever $h\circ f=k\circ f$, $g$ factors uniquely through $f$);
\item The functor $F:{\cal C}^{\textrm{op}}\to \Cont_{{\cal I}}(Aut_{\cal C}(u))$ which sends any object $c$ of $\cal C$ to the set $Hom_{\Ind{\cal C}}(c, u)$ (endowed with the obvious action by $Aut_{\cal C}(u)$) and any arrow $f:c\to d$ in $\cal C$ to the $Aut_{\cal C}(u)$-equivariant map
\[
-\circ f:Hom_{\Ind{\cal C}}(d, u)\to Hom_{\Ind{\cal C}}(c, u)
\]
is full and faithful;
\item For any objects $c,d\in {\cal C}$ and any arrows $\chi:c\to u$ and $\xi:d\to u$ in $\Ind{\cal C}$, ${\cal I}_{\xi}\subseteq {\cal I}_{\chi}$ (that is, for any automorphism $f$ of $u$, $f\circ \xi=\xi$ implies $f\circ \chi=\chi$) if and only if there exists a unique arrow $f:c\to d$ in $\cal C$ such that $\chi=\xi \circ f$:
\[  
\xymatrix {
c \ar@{-->}[d]_{f} \ar[r]^{\chi} & u \\
d \ar[ur]_{\xi} &  } 
\] 
\end{enumerate}  
\end{theorem}

\begin{proofs}
Condition $(i)$ is exactly equivalent to the condition that $J_{at}$ should be subcanonical (cf. Example C2.1.12(b) \cite{El}); but this condition is in turn equivalent, by Lemma C2.2.15 \cite{El}, to the requirement that the composite $l:{\cal C}^{\textrm{op}}\to \Sh({\cal C}^{\textrm{op}}, J_{at})$ of the Yoneda embedding ${\cal C}^{\textrm{op}}\to [{\cal C}, \Set]$ with the associated sheaf functor $a_{J_{at}}:[{\cal C}, \Set]\to \Sh({\cal C}^{\textrm{op}}, J_{at})$ be full and faithful, that is to condition $(ii)$.

To prove the equivalence of property $(iii)$ with the condition for $l$ to be full and faithful, we observe that, given an equivariant map $\gamma:F(c)\to F(d)$, we can describe $\gamma$ as follows. Take $\xi\in F(c)$ and $\chi=\gamma(\xi)$; then ${\cal I}_{\xi}\subseteq {\cal I}_{\chi}$ and $\gamma$ can be identified with the map sending $g\circ \xi$ to $g\circ \chi$ for every automorphism $g$ of $u$. Conversely, if we have two arrows $\chi:c\to u$ and $\xi:d\to u$ in $\Ind{\cal C}$ such that ${\cal I}_{\xi}\subseteq {\cal I}_{\chi}$ then the map $F(c)\to F(d)$ sending $g\circ \xi$ to $g\circ \chi$ for every automorphism $g$ of $u$ is well-defined and $Aut_{\cal I}(u)$-equivariant. To conclude, it suffices to notice that such a map is equal to $F(f)$ for an arrow $f:d\to c$ in ${\cal C}$ if and only if $\xi \circ f=\chi$.      
\end{proofs}

The theorem motivates the following definition: given a small category $\cal C$ satisfying AP and JEP, and a $\cal C$-universal and $\cal C$-ultrahomogeneous object $u$ in $\Ind{\cal C}$, we say that the pair $({\cal C}, u)$ defines a \emph{concrete Galois theory} if every arrow of $\cal C$ is a strict monomorphism.  

\begin{remark}\label{completeformulae}
A canonical example of categories satistfying the conditions of Theorem \ref{subcanonical} is given by the subcategories of the syntactic categories of theories $\mathbb T$ satisfying the hypotheses of Theorem \ref{main} on the $\mathbb T$-complete formulae. We shall see other examples of categories satisfying this condition in section \ref{examples}. 
\end{remark}

Under the hypotheses that all the arrows in $\cal C$ are strict monomorphisms, we can construct alternative `Galois representations' of the topos $\Sh({\cal C}^{\textrm{op}}, J_{at})$, as follows.

\begin{theorem}\label{catd}
Suppose that ${\cal C}\hookrightarrow {\cal D}$ is a full embedding of a category $\cal C$ into a category $\cal D$ containing a $\cal C$-ultrahomogeneous and $\cal C$-universal object $u$. Suppose that the following two properties are satisfied:
\begin{enumerate}[(a)]
\item For any objects $a, b\in {\cal C}$ and arrows $\xi:a\to u$ and $\chi:b\to u$ in $\cal D$ there exists an object $c$ of $\cal C$, arrows $f:a\to c$ and $g:b\to c$ in $\cal C$ and an arrow $\epsilon:c\to u$ in $\cal D$ such that $\epsilon \circ f=\xi$ and $\epsilon \circ g= \chi$;  

\item For any arrows $f,g:a\to b$ in $\cal C$ and any arrow $\chi:b\to u$, $\chi\circ f=\chi\circ g$ implies $f=g$.
\end{enumerate}
Then the category $\cal C$ satisfies AP and JEP, the sets of the form ${\cal I}_{\chi}:=\{f:u\cong u \textrm{ | } f\circ \chi=\chi\}$ define an algebraic base for the group $Aut_{\cal D}(u)$ of automorphisms of $u$ in $\cal D$ and we have a functor $F:{\cal C}^{\textrm{op}}\to \Cont(Aut_{\cal D}(u))$ taking values in the subcategory of transitive $Aut_{\cal D}(u)$-actions which sends every object of $\cal C$ to the set $Hom_{\cal D}(c, u)$, endowed with the obvious action of $Aut_{\cal D}(u)$, and any arrow $f:d\to c$ of $\cal C$ to the $Aut_{\cal D}(u)$-equivariant map $Hom_{\cal D}(-, f):Hom_{\cal D}(c, u) \to Hom_{\cal D}(d, u)$. Moreover, the following conditions are equivalent:
\begin{enumerate}[(i)]
\item All the arrows of $\cal C$ are strict monomorphisms; 

\item The functor $F$ is full and faithful and induces an equivalence
\[
\Sh({\cal C}^{\textrm{op}}, J_{at})\simeq \Cont(Aut_{\cal D}(u)).
\] 
\end{enumerate}
\end{theorem}

\begin{proofs}
The fact that $\cal C$ satisfies the amalgamation property follows from the $\cal C$-universality of $u$ and properties $(a)$ and $(b)$, while the fact that $\cal C$ satisfies the joint embedding property follows from property $(a)$. The fact that the sets ${\cal I}_{\chi}$ form an algebraic base for $Aut_{\cal D}(u)$ also follows from property $(a)$. 

The fact that $(ii)$ implies $(i)$ follows as in the proof of Theorem \ref{subcanonical}. To prove that $(i)$ implies $(ii)$ it suffices to prove that the functor $F$ is full and faithful, since if this is the case we have a full and faithful functor ${\cal C}^{\textrm{op}}\to \Cont(Aut_{\cal D}(u))$ whose image can be identified with the dense subcategory of $\Cont(Aut_{\cal D}(u))$ corresponding to the algebraic base ${\cal I}_{\cal C}$ as in Remark \ref{algbases}, and hence the toposes of sheaves $\Sh({\cal C}^{\textrm{op}}, J_{at})$ and $\Cont(Aut_{\cal D}(u))$ on the two categories with respect to the atomic topology are equivalent. Let us thus suppose that $\xi:a\to u$ and $\chi:b\to u$ are arrows in $\cal D$ from objects of $\cal C$ such that ${\cal I}_{\chi}\subseteq {\cal I}_{\xi}$. By property $(a)$ there exist arrows $f:a\to c$ and $g:b\to c$ in $\cal C$ and an arrow $\epsilon:c\to u$ in $\cal D$ such that $\epsilon \circ f=\xi$ and $\epsilon \circ g= \chi$. By property $(b)$, $\chi$ factors (uniquely) through $\xi$ if and only if $g$ factors (uniquely) through $f$; for this, since all the arrows of $\cal C$ are strict monomorphisms, it is equivalent to verify that for any pair of arrows $h,k:b\to e$, if $h\circ g=k\circ g$ then $h\circ f=k\circ f$. Let us therefore suppose that $h\circ g=k\circ g$. Since $u$ is $\cal C$-homogeneous, there exist arrows $\alpha:e\to u$ and $\beta:e\to u$ in $\cal D$ such that $\alpha \circ h=\epsilon$ and $\beta \circ k=\epsilon$. Now, since $u$ is ${\cal C}$-ultrahomogeneous, there exists an automorphism $j:u\cong u$ of $u$ in $\cal D$ such that $\beta=j\circ \alpha$. We have that $j\circ \chi=j\circ \epsilon \circ g=j\circ \alpha \circ h\circ g=\beta \circ h \circ g= \beta\circ k\circ g=\epsilon \circ g=\chi$; therefore, as ${\cal I}_{\chi}\subseteq {\cal I}_{\xi}$, we have that $j\circ \xi=\xi$; but $j\circ \xi=j\circ \alpha \circ h\circ f=\beta \circ h\circ f$, and $\xi=\beta \circ k\circ f$; hence $\beta \circ (h\circ f)=\beta \circ (k\circ f)$, from which it follows, by condition $(b)$, that $h\circ f=k\circ f$, as required.   
\end{proofs}

\begin{remarks}\label{necessarycondition}

\begin{enumerate}[(a)]
\item Under the hypothesis that all the arrows of $\cal C$ are strict mono-\\morphisms, condition $(b)$ in the statement of the theorem is necessary; indeed, the fact that the functor $F$ induces an equivalence of toposes implies that for any arrows $\chi:c\to u$ and $\xi:d\to u$ there exists at most one factorization of $\chi$ through $\xi$. Now, if for two arrows $f,g:a\to b$ in $\cal C$ and an arrow $\chi:b\to u$ we have $\chi\circ f=\chi\circ g$, $f$ and $g$ are two factorizations of $\chi\circ f=\chi\circ g$ through $\chi$ and hence $f=g$. 

\item In the case all the arrows of $\cal C$ are strict monomorphisms, Theorem \ref{maincategorical} can be deduced as the particular case of Theorem \ref{catd} when ${\cal C}\hookrightarrow$ is the embedding ${\cal C}\hookrightarrow \Ind{\cal C}$.
\end{enumerate}
\end{remarks}

Let us now discuss the relationship between a category $\cal D$ with an object $u$ as in the hypotheses of Theorem \ref{catd} and the ind-completion $\Ind{\cal C}$ of the category $\cal C$.

In \cite{El} (section C4.2), the following criterion for a full embedding ${\cal C}\hookrightarrow {\cal D}$ of categories to be of the form ${\cal C}\hookrightarrow \Ind{\cal C}$ is given: the category $\cal D$ has all (small) filtered colimits, every object of $\cal C$ is finitely presentable in $\cal D$ (in the sense that the corresponding representable hom functor preserves small filtered colimits) and every object of $\cal D$ can be expressed as a (small) filtered colimit of objects of $\cal C$. The embedding ${\cal C}\hookrightarrow \Ind{\cal C}$ can also be characterized as the filtered-colimit completion of the category $\cal C$.

Notice that if ${\cal C}\hookrightarrow {\cal D}$ is isomorphic to ${\cal C}\hookrightarrow \Ind{\cal C}$ then we have an equivalence ${\cal D}\simeq \Ind{\cal C}\simeq \mathbf{Flat}({\cal C}^{\textrm{op}}, \Set)$ sending any object $d$ of $\cal D$ to the flat functor $Hom_{\cal D}(-, d):{\cal C}^{\textrm{op}}\to \Set$, and hence the three properties in the following proposition (expressing the fact that the functor $Hom_{\cal D}(-, d):{\cal C}^{\textrm{op}}\to \Set$ is flat, i.e. that its category of elements is filtered) hold:

\begin{proposition}\label{propind}
Let ${\cal C}\hookrightarrow {\cal D}$ be a full embedding isomorphic to ${\cal C}\hookrightarrow \Ind{\cal C}$. Then the following three properties hold:
\begin{enumerate}[(i)]

\item Every object $d$ of $\cal D$ is the colimit of the canonical diagram on ${\cal C}\slash d$ in $\cal D$; 

\item For any object $d$ of $\cal D$ and any arrows $\xi:a\to d$ and $\chi:b\to d$ in $\cal D$ from objects $a$ and $b$ of $\cal C$, there exists an object $c$ of $\cal C$, arrows $f:a\to c$ and $g:b\to c$ of $\cal C$ and an arrow $\epsilon:c\to d$ in $\cal D$ such that $\epsilon \circ f=\xi$ and $\epsilon \circ g= \chi$;

\item For any arrows $f,g:a\to b$ in $\cal C$, arrows $\xi:a\to u$ in $\cal D$ and $\chi:b\to u$ such that $\chi \circ f=\xi$ and $\chi \circ g=\xi$, there exists an arrow $h:b\to c$ in $\cal C$ and an arrow $\epsilon:c\to u$ in $\cal D$ such that $\epsilon \circ h =\chi$ and $h\circ f=h\circ g$. 
\end{enumerate}
\end{proposition}\qed

\begin{remark}
From the theorem it follows that the two conditions in the statement of Theorem \ref{catd} are automatically satisfied if the embedding ${\cal C}\hookrightarrow {\cal D}$ is (isomorphic to) the canonical embedding ${\cal C}\hookrightarrow \Ind{\cal C}$.
\end{remark}

The following proposition provides some convenient criteria for identifying ind-completions.

\begin{proposition}

\begin{enumerate}[(i)]
\item Let $\mathbb T$ be a theory of presheaf type. Then the category ${\mathbb T}\textrm{-mod}(\Set)$ can be identified as the ind-completion of the category $\textrm{f.p.}{\mathbb T}\textrm{-mod}(\Set)$; in particular, the embedding $\textrm{f.p.}{\mathbb T}\textrm{-mod}(\Set)\hookrightarrow {\mathbb T}\textrm{-mod}(\Set)$ satisfies the three conditions in the statement of Proposition \ref{propind} above.

\item Let $\mathbb T$ be a geometric theory over a signature which does not contain any relation symbols, and let $\cal C$ be a set of finitely generated $\mathbb T$-models such that every $\mathbb T$-model can be expressed as a directed union of models in $\cal C$. Then we have ${\mathbb T}\textrm{-mod}_{i}(\Set)\simeq \Ind{\tilde{\cal C}}$, where ${\mathbb T}\textrm{-mod}_{i}(\Set)$ is the category of $\mathbb T$-models in $\Set$ and injective homomorphisms between them and $\tilde{\cal C}$ is the category of $\mathbb T$-models in $\cal C$ and injective homomorphisms between them.  

\item Let $\mathbb T$ be a finitary first-order theory and ${\mathbb T}\textrm{-mod}_{e}(\Set)$ be the category of $\mathbb T$-models in $\Set$ and elementary embeddings between them. If $\cal C$ is a set of finitely generated $\mathbb T$-models in $\Set$ such that every $\mathbb T$-model can be expressed as a filtered colimit in ${\mathbb T}\textrm{-mod}_{e}(\Set)$ of models in $\cal C$. Then we have ${\mathbb T}\textrm{-mod}_{e}(\Set)\simeq \Ind{\tilde{\cal C}}$, where $\tilde{\cal C}$ is the category of models in $\cal C$ and elementary embeddings between them.     
\end{enumerate}

\end{proposition}

\begin{proofs}
$(i)$ If $\mathbb T$ is a theory of presheaf type then $\mathbb T$ is classified by the presheaf topos $[\textrm{f.p.}{\mathbb T}\textrm{-mod}(\Set), \Set]$, and hence by Diaconescu's equivalence we have that ${\mathbb T}\textrm{-mod}(\Set)\simeq \mathbf{Flat}(\textrm{f.p.}{\mathbb T}\textrm{-mod}(\Set)^{\textrm{op}}, \Set)$. 

$(ii)$ The theory ${\mathbb T}'$ obtained from $\mathbb T$ by adding a binary predicate and the coherent sequents asserting that such predicate is complemented to the equality relation is clearly geometric and hence its category of models, which coincides with the category of $\mathbb T$-models in $\Set$ and injective homomorphisms between them, has filtered colimits (cf. Lemma D2.4.9 \cite{El}). Now, if the signature of $\mathbb T$ does not contain any relation symbols, any finitely generated $\mathbb T$-model $A$ is finitely presentable; indeed, for any $\mathbb T$-model $M$ equal to the directed colimit $colim(M_{i})$ of a diagram $M_{i}$ of models, any homomorphism $f:A\to M$ factors, as a function, through a canonical embedding $M_{i}\hookrightarrow M$ (indeed, $A$ being finitely generated, it suffices to choose $i$ so that the image under $f$ of any of the generators of $A$ belongs to $M_{i}$), and such factorization must be a $\mathbb T$-model homomorphism since the signature of $\mathbb T$ does not contain any relation symbol and the arrows $M_{i}\to M$ are injective. The above-mentioned characterization for ind-completions thus applies. 

$(iii)$ As observed at p. 887 \cite{El}, the category ${\mathbb T}\textrm{-mod}_{e}(\Set)$ has all small filtered colimits, the proof of condition $(ii)$ shows that any factorization $A\to M_{i}$ of a $\mathbb T$-model homomorphism $f:A\to M$ from a finitely generated $\mathbb T$-model $A$ to a model $M\cong colim(M_{i})$ is a $\mathbb T$-model homomorphism; indeed, the fact that it respects function symbols follows from the injectivity of the morphisms in ${\mathbb T}\textrm{-mod}_{e}(\Set)$, while the fact that it respects relation symbols follows from the fact that such maps are elementary embeddings and hence in particular they reflect the relations over the signature of $\mathbb T$.   
\end{proofs}
 
Given a category $\cal D$ with an object $u$ as in the hypotheses of Theorem \ref{catd}, we can naturally obtain an object $\tilde{u}$ of $\Ind{\cal C}$, as follows. The hypotheses of Theorem \ref{catd} ensure that the category ${\cal C}\slash u$ is filtered; we can thus consider the colimit $\tilde{u}$ in $\Ind{\cal C}$ of the canonical functor ${\cal C}\slash u \to {\cal C}$. Let us show that $\tilde{u}$ is a $\cal C$-universal and $\cal C$-ultrahomogeneous object. Let us denote by $J_{(c, \xi)}:c\to \tilde{u}$ the canonical colimit arrows in $\Ind{\cal C}$. The fact that $u$ is $\cal C$-universal implies the existence, for any object $c$ of $\cal C$, of an arrow $\xi:c \to u$ in ${\cal D}$; hence for every object $c$ of $\cal C$ we have an arrow $c\to \tilde{u}$ in $\Ind{\cal C}$, namely $J_{(c, \xi)}$. To prove that $\tilde{u}$ is $\cal C$-ultrahomogeneous, we first observe that any arrow $r:c\to \tilde{u}$ in $\Ind{\cal C}$ from an object $c$ of $\cal C$ to $\tilde{u}$ is of the form $J_{(c, \xi)}$ for some object $(c, \xi)\in {\cal C}\slash u$. As $c$ is finitely presentable in $\Ind{\cal C}$, we have a factorization of $r$ as $J_{(d, \chi)}\circ f$ for some arrow $f:c\to d$ in $\cal C$; but $(c, \xi \circ f)$ is an object of ${\cal C}\slash u$, and $f$ is an arrow $(c, \xi \circ f)\to (d, \chi)$ in ${\cal C}\slash u$, whence $r=J_{(d, \chi)}\circ f=J_{(c, \xi \circ f)}$. Now, given two arrows $J_{(c, \xi_{1})}, J_{(c, \xi_{2})}:c\to \tilde{u}$ in $\Ind{\cal C}$, by the $\cal C$-ultrahomogeneity of $u$ we have an automorphism $j:u\cong u$ of $u$ in $\cal D$ such that $j\circ \xi_{1}=\xi_{2}$. By using the universal property of the colimit $\tilde{u}$, we can obtain an automorphism $\tilde{j}:\tilde{u}\cong \tilde{u}$ of $\tilde{u}$ in $\Ind{\cal C}$ such that $\tilde{j}\circ J_{(c, \xi_{1})}=J_{(c, \xi_{2})}$, as follows. The arrows $J_{(d, j\circ \chi)}: d\to \tilde{u}$ (for $(d, \chi)\in {\cal C}\slash u$) clearly form a cocone in $\Ind{\cal C}$ on ${\cal C}\slash u$ with vertex $\tilde{u}$; we set $\tilde{j}:\tilde{u}\to \tilde{u}$ equal to the unique arrow $s$ in $\Ind{\cal C}$ such that $s\circ J_{(d, \chi)}=J_{(d, j\circ \chi)}$ for any object $(d, \chi)$ of ${\cal C}\slash u$. The arrow $\tilde{j}$ is an automorphism because $\tilde{j^{-1}}$ is an inverse for it, and it clearly satisfies the required property. This completes the proof that $\tilde{u}$ is $\cal C$-universal and $\cal C$-ultrahomogeneous. Therefore the object $\tilde{u}$ satisfies the hypotheses of Theorem \ref{maincategorical}, and we have an equivalence
\[
\Sh({\cal C}^{\textrm{op}}, \Set)\simeq \Cont(Aut_{\cal C}(\tilde{u})),
\]     
which, combined with the equivalence of Theorem \ref{catd}, exhibits the topological group $Aut_{\cal C}(\tilde{u})$ as the completion (in the sense of section \ref{completegroups}) of the topological group $Aut_{\cal D}(u)$.

Summarizing, we have the following result.

\begin{theorem}\label{maincategoricalD}
Let $({\cal C}\hookrightarrow {\cal D}, u)$ be a pair satisfying all the hypotheses of Theorem \ref{catd}. Then there is a $\cal C$-universal and $\cal C$-ultrahomogeneous object $\tilde{u}$ in $\Ind{\cal C}$ and equivalences of toposes 
\[
\Sh({\cal C}^{\textrm{op}}, \Set)\simeq \Cont(Aut_{\cal D}(u)) \simeq \Cont(Aut_{\cal C}(\tilde{u})),
\]
\end{theorem}\qed

Notice that this theorem shows in particular that there can be many different Morita-equivalent topological groups (i.e., topological groups having equivalent toposes of continuous actions).

\subsubsection{Regular and strict monomorphisms}

Let $\cal C$ be a small category. Recall that an arrow $f:d\to c$ in $\cal C$ is said to be a \emph{strict monomorphism} if for any arrow $g:e\to c$ such that $h\circ g=k\circ g$ whenever $h\circ f=k\circ f$, $g$ factors uniquely through $f$; in other words, $f$ is the limit of the (possibly large) diagram consisting of all the pairs $(h,k)$ which coequalize $f$. The dual concept of strict epimorphism was introduced by Grothendieck in \cite{grothendieck}, Exposé I, 10.2 and 10.3.

An arrow $f$ is said to be a \emph{regular monomorphism} in a category $\cal C$ if its cokernel pair exists and $f$ is the equalizer of it. Clearly, for any arrow $f$ of a category $\cal C$ having a cokernel pair in $\cal C$, $f$ is a regular monomorphism if and only if it is a strict monomorphism; this shows that the notion of strict monomorphism represents a natural generalization of the notion of regular monomorphism applicable to categories in which cokernel pairs of arrows do not necessarily exist.

The following result establishes relationships between strict monomorphisms in different categories.

\begin{proposition}\label{regstrict}

\begin{enumerate}[(i)]
\item Let ${\cal C}\hookrightarrow {\cal D}$ be an embedding of a small category $\cal C$ into a category $\cal D$ satisfying the following properties: 

\begin{enumerate}[(a)]
\item For any pair of composable arrows $g,h$ in $\cal D$ with domains and codomains in $\cal C$, if $g\circ h$ belongs to $\cal C$ then $h$ belongs to $\cal C$;  

\item For any objects $a, b$ of ${\cal C}$, $d$ of $\cal D$ and arrows $\xi:a\to d, \chi:b\to d$ in $\cal D$, there exists an object $c$ of $\cal C$, arrows $f:a\to c$ and $g:b\to c$ in $\cal C$ and an arrow $\epsilon:c\to d$ in $\cal D$ such that $\epsilon \circ f=\xi$ and $\epsilon \circ g= \chi$;  

\item For any arrows $f,g:a\to b$ in $\cal C$, object $d$ of $\cal D$ and arrow $\chi:b\to d$, $\chi\circ f=\chi\circ g$ implies $f=g$.
\end{enumerate}

Then any arrow of $\cal C$ which is a strict monomorphism in $\cal D$ is a strict monomorphism in $\cal C$. 

In particular (cf. point (ii) below), if $\cal C$ is a small category with all arrows monic then for any arrow $f$ in $\cal C$, $f$ is a strict monomorphism in $\cal C$ if and only if it is a strict monomorphisms in $\Ind{\cal C}$. 

\item Let ${\cal C}\hookrightarrow {\cal D}$ be a full embedding. If any object of $\cal D$ can be written as a colimit of objects of $\cal D$ (that is, of a diagram on $\cal D$ which factors through the embedding ${\cal C}\hookrightarrow {\cal D}$) then every strict monomorphism in $\cal C$ is a strict monomorphism in $\cal D$.  

\item If the domains $A$ of cokernel pairs of arrows in $\cal C$ exist in $\cal D$ and are limits in $\cal D$ of the canonical diagram on $A\slash {\cal C}$ then if $f$ is a regular monomorphism in $\cal D$, $f$ is a strict monomorphism in $\cal C$.   
\end{enumerate}

\end{proposition}
 
\begin{proofs}
$(i)$ This immediately follows from the fact that the two conditions in the statement of part $(i)$ of the proposition ensure that the limit in $\cal D$ over the diagram formed by the pairs of arrows which coequalize $f$ in $\cal D$, if it exists, coincide with the limit in $\cal C$ over the diagram formed by the pairs of arrows which coequalize $f$ in $\cal C$; in other words, if $f$ is a strict monomorphism in $\cal D$ then it is a strict monomorphism in $\cal C$. 

For any small category $\cal C$ whose arrows are all monic, the embedding ${\cal C}\hookrightarrow \Ind{\cal C}$ is full and satisfies property $(a)$ since for every object $d$ of $\Ind{\cal C}$ the category ${\cal C}\slash d$ satisfies the joint embedding property (it being filtered), and property $(b)$ since the category ${\cal C}\slash d$ satisfies the weak coequalizer property (it being filtered) and every morphism in $\cal C$ is monic. 

$(ii)$ Suppose that every object $d$ of $\cal D$ is a colimit of a diagram $D:{\cal I}\to {\cal D}$ which factors through the embedding ${\cal C}\hookrightarrow {\cal D}$, and let $f:a\to b$ be a strict monomorphism in $\cal C$. To prove that $f$ is a strict monomorphism in $\cal D$, suppose that $g:d\to b$ is an arrow in $\cal D$ such that any pair of arrows in $\cal D$ which coequalize $f$ coequalize $g$. Consider the composites $g_{i}:=D(i)\to b$ of the canonical colimit arrows $\xi_{i}:D(i)\to d$ with $g$; these arrows are in $\cal C$, as the embedding ${\cal C}\hookrightarrow {\cal D}$ is full, and clearly each of them satisfies the property of being coequalized by any pair of arrows in $\cal C$ which coequalizes $f$; therefore, as $f$ is a strict monomorphism, for every $i\in {\cal I}$ there exists a unique arrow $\gamma_{i}:D(i)\to a$ in $\cal C$ such that $f\circ \gamma_{i}=g_{i}$ (for all $i\in {\cal I}$). Now, by the uniqueness of the factorization in the definition of strict monomorphism, the arrows $\gamma_{i}:D(i)\to a$ (for $i\in {\cal I}$) define a conone on $D$ with vertex $a$ and hence they induce a unique arrow $\gamma:d\to a$ such that $\gamma \circ \xi_{i}=\gamma_{i}$ and therefore $f\circ \gamma=g$. The uniqueness of a factorization of $g$ through $f$ in $\cal D$ follows from the uniqueness of the factorization of each of the arrows $g_{i}$ through $f$.     

$(iii)$ Suppose that $f:d\to c$ is an arrow of $\cal C$ which is a regular monomorphism in $\cal D$ with cokernel pair  $\xi_{1}, \xi_{2}:c\to K_{f}$, and let us suppose that $g:e\to c$ is an arrow in $\cal C$ such that any pair of arrows in $\cal C$ which coequalizes $f$ coequalizes $g$; we want to prove that $g$ factors uniquely through $f$. Let us denote by $\chi_{1}, \chi_{2}:c\to K_{g}$ the cokernel pair of $g$. Then for any arrow $s:K_{f}\to a$ to an object $a$ of $\cal C$ there exists a unique arrow $t_{s}:K_{g}\to a$ such that $t_{s}\circ \chi_{1}=s\circ \xi_{1}$ and $t_{s}\circ \chi_{2}=s\circ \xi_{2}$. The arrows $t_{s}:K_{g}\to a$ form a cone with vertex $K_{g}$ on the diagram $D$ on $K_{f}\slash {\cal C}$ with vertex $K_{f}$ and legs $s:K_{f}\to a$ (indeed, if $z\circ s=s'$ for some arrows $s':K_{f}\to b$ in $\cal D$ and $z:a\to b$ in $\cal C$ then $z\circ t_{s}=t_{s'}$), and hence since by our hypothesis $D$ is limiting, there exists a unique arrow $y:K_{g}\to K_{f}$ such that for any $s:K_{f}\to a$ (with $a\in {\cal C}$), $s\circ y=t_{s}$; also, the universal property of the limit $K_{f}$ ensures that $y\circ \chi_{1}=\xi_{1}$ and $y\circ \chi_{2}=\xi_{2}$. Therefore, as $f$ is the equalizer of $\xi_{1}$ and $\xi_{2}$, $g$ factors uniquely through $f$. 
\end{proofs} 

The question now naturally arises of how one can obtain embeddings ${\cal C}\hookrightarrow {\cal D}$ satisfying the hypotheses of the proposition. In this respect, the following result is useful.

\begin{proposition}\label{srmono} 
\begin{enumerate}[(i)]

\item Let $\cal C$ a small category, $\Ind{\cal C}$ its ind-completion and $\cal A$ a subcategory of $\Ind{\cal C}$ whose objects are exactly the objects of $\Ind{\cal C}$ and which satisfies the following properties:

\begin{enumerate}[(a)]
\item The composite of two arrows $g\circ h$ belongs to $\cal A$ then $h$ belongs to $\cal A$,

\item $\cal A$ is closed under filtered colimits in $\Ind{\cal C}$, and

\item any object of $\Ind{\cal C}$ can be written as a filtered colimit of a diagram in $\cal C$ with values in $\cal A$. 
\end{enumerate}
Then, denoted by ${\cal C}_{\cal A}$ the full subcategory of $\cal A$ on the objects of $\cal C$, we have ${\cal A}\simeq \Ind{{\cal C}_{\cal A}}$.

\item Let $\cal C$ be a small category, ${\cal B}$ a full subcategory of $\cal C$ and $\cal D$ be a full subcategory of its ind-completion $\Ind{\cal C}$. If ${\cal B}\subseteq {\cal D}$, $\cal D$ is closed in $\Ind{\cal C}$ under filtered colimits and every object of $\cal D$ can be written as a filtered colimit of objects in $\cal B$ then ${\cal D}\simeq \Ind{{\cal B}}$. 
\end{enumerate}      
\end{proposition}

\begin{proofs}
$(i)$ We use the following classical criterion for an embedding to be isomorphic to the embedding of a small category into its ind-completion (cf. section C4.2 \cite{El}): for a full embedding ${\cal C}\hookrightarrow {\cal D}$ to be isomorphic to the embedding ${\cal C}\hookrightarrow \Ind{\cal C}$ it is necessary and sufficient that the category $\cal D$ has all (small) filtered colimits, every object of $\cal C$ is finitely presentable in $\cal D$ and every object of $\cal D$ can be expressed as a (small) filtered colimit of objects of $\cal C$. The first condition is satisfied by condition $(b)$, the third condition is satisfied by condition $(c)$. To verify the second condition, we use the explicit characterization of finitely presentable objects provided by Lemma \ref{fp}; it is immediately seen that condition $(a)$ ensures that every finitely presentable object in $\Ind{\cal C}$ is finitely presentable also an object of $\cal A$.\\         
$(ii)$ This result equally follows from the same criterion for ind-completions, using the explicit characterization of finitely presentable objects provided by Lemma \ref{fp} to prove that every object of $\cal B$ is finitely presentable in $\cal D$.    
\end{proofs}

\begin{remarks}\label{rem1}

\begin{enumerate}[(a)]

\item Condition $(i)$ in the proposition is satisfied by the subcategory $\cal A$ consisting of the objects of $\Ind{\cal C}$ and strict monomorphisms between them;

\item For any geometric theory $\mathbb T$ such that the category ${\mathbb T}\textrm{-mod}(\Set)$ can be expressed as $\Ind{\cal C}$ of a full subcategory $\cal C$ of finitely presentable $\mathbb T$-models, the subcategory ${\mathbb T}\textrm{-mod}_{i}(\Set)$ of ${\mathbb T}\textrm{-mod}(\Set)$ consisting of the $\mathbb T$-models in $\Set$ and injective $\mathbb T$-model homomorphisms between them satisfies condition $(i)$ (obviously) and condition $(ii)$ (by Proposition \ref{injective}); condition $(iii)$ holds if and only if every $\mathbb T$-model can be expressed as a directed union of models in $\cal C$. 
\end{enumerate}

\end{remarks}

\begin{proposition}\label{propinj}
Let $\mathbb T$ be a geometric theory whose category ${\mathbb T}\textrm{-mod}(\Set)$ can be expressed as $\Ind{\cal C}$, for a full subcategory $\cal C$ of finitely presentable $\mathbb T$-models (for instance, a theory of presheaf type $\mathbb T$) and such that every $\mathbb T$-model can be expressed as a directed union of models in $\cal C$. Then ${\mathbb T}\textrm{-mod}_{i}(\Set)\simeq \Ind{{\cal C}_{i}}$, where ${\cal C}_{i}$ is the category whose objects are the objects of $\cal C$ and whose arrows are the injective $\mathbb T$-model homomorphisms between them, and the category ${\cal C}_{sm}$ of objects of $\cal C$ and arrows between them which are strict monomorphisms in ${\mathbb T}\textrm{-mod}(\Set)$ has the property that every arrow in it is a strict monomorphism.  
\end{proposition}

\begin{proofs}
This immediately follows from Proposition \ref{srmono}(i) in view of Remark \ref{rem1}(b) and by Proposition \ref{regstrict}(i). 
\end{proofs}

\begin{remark}\label{application}
This proposition can be profitably applied in connection with theories $\mathbb T$ whose category ${\mathbb T}\textrm{-mod}(\Set)$ of set-based models possesses cokernel pairs of arrows (for instance, for algebraic theories $\mathbb T$) and for any arrow in ${\cal C}_{i}$ its cokernel pair lies in ${\mathbb T}\textrm{-mod}_{i}(\Set)$; indeed, in this case the subcategory ${\cal C}_{sm}$ coincides with the category whose objects are the objects of $\cal C$ and whose arrows are the injective $\mathbb T$-model homomorphisms between them which are regular monomorphisms in the category ${\mathbb T}\textrm{-mod}(\Set)$. In fact, in many algebraic categories (e.g., the category of vector spaces and linear maps, the category of groups, the category of Boolean algebras etc.), all monomorphisms are regular, and monomorphisms in categories of models of algebraic theories are precisely the injective model homomorphisms (cf. Proposition 3.4.1 \cite{borceux}); the problem of whether the cokernel pair of a monomorphism, provided that it exists, is again a pair of monomorphisms, is non-trivial in general, and a considerable amount of literature has been written on this subject (cf. for instance \cite{amalg1}, \cite{amalg2},\cite{amalg3},\cite{amalg4},\cite{amalg5},\cite{amalg6}). 
\end{remark}

Finally, we note that if $\mathbb T$ is a theory of presheaf type then in studying the property of an arrow in $\textrm{f.p.}{\mathbb T}\textrm{-mod}(\Set)$ to be a strict monomorphism one can exploit the syntactic characterization of finitely presentable $\mathbb T$-models as the finitely presented models of $\mathbb T$ (cf. \cite{OCS}), which provides a (natural) bijective correspondence between the arrows $M \to N$ in $\textrm{f.p.}{\mathbb T}\textrm{-mod}(\Set)$ from such a model $M$ presented by a formula $\phi(\vec{x})$ to any (finitely presentable) $\mathbb T$-model $N$ and the elements of the interpretation $[[\phi(\vec{x})]]_{N}$ of $\phi$ in $N$.

\subsection{Remarks on the general case}\label{sheafif}

We have seen in the last section that an equivalence of toposes
\[
\Sh({\cal C}^{\textrm{op}}, J_{at})\simeq \Cont(Aut_{\cal C}(u))
\] 
provided by Theorem \ref{maincategorical} restricts to a `Galois-type' categorical equivalence of ${\cal C}^{\textrm{op}}$ with a full subcategory of $\Cont(Aut_{\cal C}(u))$ if and only if the topology $J_{at}$ on the category $\cal C$ is subcanonical (equivalently, all the arrows in $\cal C$ are strict monomorphisms). Still, it is natural to ask what one can say about the relationships between the $Aut_{\cal C}(u)$-equivariant maps $Hom_{\Ind{\cal C}}(c, u)\to Hom_{\Ind{\cal C}}(d, u)$ (for objects $c, d$ of $\cal C$) in the general case. Clearly, the equivalence of Theorem \ref{maincategorical} sends, for each $c\in {\cal C}$, the object $a_{J_{at}}(yc)$ of the topos $\Sh({\cal C}^{\textrm{op}}, J_{at})$ to the $Aut_{\cal C}(u)$-set $Hom_{\Ind{\cal C}}(c, u)$, where $yc$ is the representable functor associated to $c$ and $a_{J_{at}}$ is the associated sheaf functor. Therefore the arrows $Hom_{\Ind{\cal C}}(c, u)\to Hom_{\Ind{\cal C}}(c, u)$ in $\Cont(Aut_{\cal C}(c))$ correspond to the arrows $a_{J_{at}}(yc)\to a_{J_{at}}(yd)$ in $\Sh({\cal C}^{\textrm{op}}, J_{at})$. To characterize such arrows in more explicit terms we shall delve in a general analysis valid for any site.         

Let $({\cal C}, J)$ be a site, $Y:{\cal C}\to [{\cal C}^{\textrm{op}}, \Set]$ be the Yoneda embedding and $a_{J}:[{\cal C}^{\textrm{op}}, \Set]\to \Sh({\cal C}, J)$ be the associated sheaf functor.  Let $C$ and $C'$ be objects of $[{\cal C}^{\textrm{op}}, \Set]$. Suppose that $r:a_{J}(C)\to a_{J}(C')$ is an arrow in $\Sh({\cal C}, J)$. Then, the unit of the adjunction between the inclusion $\Sh({\cal C}, J)\hookrightarrow [{\cal C}^{\textrm{op}}, \Set]$ and the associated sheaf functor $a_{J}$ provides an arrow $\eta_{C'}:C\to a_{J}(C')$ in $[{\cal C}^{\textrm{op}}, \Set]$, while the arrow $r$ corresponds, via this adjunction, to an arrow $r_{a}:C\to a_{J}(C')$ in $[{\cal C}^{\textrm{op}}, \Set]$. Let us consider the pullback  
\[  
\xymatrix {
A \ar[d]^{\alpha} \ar[r]^{\beta} & C' \ar[d]^{\eta_{C'}}\\
C \ar[r]_{r_{a}} & a_{J}(C'), } 
\] 
in the topos $[{\cal C}^{\textrm{op}}, \Set]$. 

Notice that $a_{J}(\alpha)$ is an isomorphism, while $a_{J}(\beta)$ is isomorphic to $r$. The arrow $r$ can thus be identified with the pair $(\alpha, \beta)$; conversely, any pair of arrows $(\alpha:A\to C, \beta:A\to C')$ in $[{\cal C}^{\textrm{op}}, \Set]$ with common domain with the property that $a_{J}(\alpha)$ is an isomorphism gives rise to an arrow $a_{J}(C)\to a_{J}(C')$ in $\Sh({\cal C}, J)$, namely $a_{J}(\beta)\circ a_{J}(\alpha)^{-1}$; in fact, this correspondence can be made into an reflection of the discrete category $Hom_{\Sh({\cal C}, J)}(a_{J}(C), a_{J}(C'))$ on the set of arrows $a_{J}(C)\to a_{J}(C')$ in $\Sh({\cal C}, J)$ into the category ${\cal K}_{C, C'}$ having as objects such pairs $(\alpha, \beta)$ and as arrows $(\alpha, \beta)\to (\alpha', \beta')$ the arrows $s:A\to A'$ in $[{\cal C}^{\textrm{op}}, \Set]$ such that $\alpha' \circ s=\alpha$ and $\beta' \circ s=\beta$ (where $dom(\alpha)=dom(\beta)=A$ and $dom(\alpha')=dom(\beta')=A'$).

Now, as observed in \cite{OCG}, we can give an explicit condition on a morphism $\gamma:F\to G$ in $[{\cal C}^{\textrm{op}}, \Set]$ for it to be sent by the associated sheaf functor $a_{J}$ to an isomorphism. Indeed, $a_{J}(\gamma)$ is an isomorphism if and only if its image $Im(\gamma)\mono G$ is $J$-dense and the equalizer of its kernel pair is $J$-dense; and the condition for a monomorphism $A\mono E$ in $[{\cal C}^{\textrm{op}}, \Set]$ to be $J$-dense can be explicitly described as follows: for any $x\in E(c)$, $\{f:d\to c \textrm{ | } E(f)(x)\in A(d)\}\in J(c)$.     

Notice that the objects of the category ${\cal K}_{C, C'}$ can be described in more elementary terms, that is in terms of discrete fibrations with codomain $\cal C$, by using Grothendieck's construction. Specifically, for any objects $A$ and $E$ of $[{\cal C}^{\textrm{op}}, \Set]$, the natural transformations $A\to E$ are in bijective correspondence with the functors $F:\int A \to \int E$ such that $\pi_{E}\circ F=\pi_{A}$, where $\int A$ and $\int E$ are respectively the category of elements of the presheaf $A$ and of the presheaf $E$ and $\pi_{A}$ and $\pi_{E}$ are the two canonical projections to $\cal C$. Conversely, any discrete fibration $P:{\cal B}\to {\cal C}$ with codomain $\cal C$ corresponds to a presheaf on $\cal C$ whose category of elements is isomorphic to $\cal B$, and any morphism of discrete fibrations of $\cal C$ corresponds to a natural transformation between the corresponding presheaves (whose corresponding functor between the associated categories of elements is isomorphic to the given morphism). The point of view of discrete fibrations is natural in this context because the property of a morphism $\gamma:A\to E$ to be sent by the associated sheaf functor $a_{J}$ to an isomorphism can be naturally expressed in terms of the associated morphism of fibrations. In particular, if $C$ and $C'$ are representable functors $Yc$ and $Yc'$, the discrete fibrations corresponding to them are respectively the canonical projections $p_{c}:{\cal C}\slash c \to {\cal C}$ and $p_{c'}:{\cal C}\slash c' \to {\cal C}$. The existence of a pair $(\alpha, \beta)$ satisfying the above conditions translates into the existence of a discrete fibration $r:{\cal I}\to {\cal C}$ and two morphisms of fibrations $a:r\to p_{c}$ and $b:r\to p_{c'}$ such that $a$ satisfies the (elementary condition equivalent to the) condition that its associated presheaf be sent by $a_{J}$ to an isomorphism.

\subsection{Atoms and transitive actions}

Let us now consider, in the context of an equivalence 
\[
\Sh({\cal C}^{\textrm{op}}, J_{at})\simeq \Cont(Aut_{\cal C}(u))
\] 
provided by theorem \ref{maincategorical}, the notion of atom of a topos. 

Let us first characterize the atoms of the topos $\Sh({\cal C}^{\textrm{op}}, J_{at})$. Let $A$ be an atom of this topos; then, since the $l(c)$ (for $c\in {\cal C}$) form a separating set for it, there is an arrow $e:l(c)\to A$ for some $c\in {\cal C}$; $A$ being an atom, $e$ is an epimorphism, so it is the coequalizer of its kernel pair $S\mono l(c)\times l(c)$. We can suppose, without loss of generality, $S$ to be of the form $a_{J_{at}}(R)$ where $R$ is a subobject of ${\cal C}(c, -)\times {\cal C}(c, -)$ in $[{\cal C}, \Set]$; in fact, we can suppose $R$, without loss of generality, to be an equivalence relation on ${\cal C}(c, -)$. Indeed, consider the language with just one binary relation symbol $K$ and the structure for its language in the topos $[{\cal C}, \Set]$ obtained by interpreting $K$ as $R$; then one can write down a geometric formula over this language whose interpretation $I_{Z}$ in this structure is precisely the equivalence relation on ${\cal C}(c, -)$ generated by $R$; since the associated sheaf functor $a_{J_{at}}$ is a geometric functor, $a_{J_{at}}(I_{Z})$ is isomorphic to the equivalence relation generated by $a_{J_{at}}(R)$, that is, to $a_{J_{at}}(R)$ itself, since $a_{J_{at}}(R)$ is an equivalence relation on $l(c)$.  

Notice that $R$ can be thought of as a set $R$ of pairs of arrows in $\cal C$ with common domain $c$ and common codomain such that for any $(f,g)\in R$ and any $h:cod(f)=cod(g)\to cod(h)$ in $\cal C$, $(h\circ f, h\circ g)\in R$. Let us denote by $(f,g)$ the set of pairs of arrows of the form $(h\circ f, h\circ g)$, for an arrow $h\in {\cal C}$; then $(f,g)$ is a subobject of ${\cal C}(c, -)\times {\cal C}(c, -)$ in $[{\cal C}, \Set]$, and $a_{J_{at}}(f,g)\mono l(c)\times l(c)$ can be identified with the monic part of the epi-mono factorization in the topos $\Sh({\cal C}^{\textrm{op}}, J_{at})$ of the arrow $<l(f), l(g)>:l(cod(f))=l(cod(g)) \to l(c)\times l(c)$. Therefore we can regard $a_{J_{at}}(R)$ as the union in $\Sub(l(c)\times l(c))$ of the epi-mono factorizations of all the arrows $<l(f), l(g)>$ for all the pairs $(f,g)\in R$. Notice that $a_{J_{at}}(R)$ corresponds under the equivalence $\Sh({\cal C}^{\textrm{op}}, J_{at})\simeq \Cont(Aut_{\cal C}(u))$ to the ($Aut_{\cal C}(u)$-equivariant) subset of 
$Hom_{\Ind{\cal C}}(c, u)\times Hom_{\Ind{\cal C}}(c, u)$ consisting of the collection $R_{t}$ of pairs of the form $(\chi \circ f, \chi\circ g)$ where $\chi\in Hom_{\Ind{\cal C}}(cod(f), u)$ and $(f,g)\in R$. The quotient in $\Cont(Aut_{\cal C}(u))$ by this equivalence relation on $Hom_{\Ind{\cal C}}(c, u)$ thus corresponds, under the equivalence $\Sh({\cal C}^{\textrm{op}}, J_{at})\simeq \Cont(Aut_{\cal C}(u))$, to the quotient in $\Sh({\cal C}^{\textrm{op}}, J_{at})$ of $l(c)$ by the relation $a_{J_{at}}(R)$ (notice that $R_{t}$ is indeed an equivalence relation since this concept is a topos-theoretic invariant).  

Now, let us characterize the atoms in $\Cont(Aut_{\cal C}(u))$. These are exactly, up to isomorphism, the $Aut_{\cal C}(u)$-sets of the form $Aut_{\cal C}(u)\slash U$, where $U$ is any open subgroup of $Aut_{\cal C}(u)$; notice that, since ${\cal I}_{\cal C}$ is an algebraic base for $Aut_{\cal C}(u)$ then any such $U$ contains an open subgroup of the form ${\cal I}_{\xi}$ for some $\xi:c\to u$, so that we have an epimorphism $Aut_{\cal C}(u)\slash {\cal I}_{\xi}\to Aut_{\cal C}(u)\slash U$. From the preceding remarks, we thus deduce that there exists an equivalence relation $R$ on ${\cal C}(-, c)$ in $[{\cal C}, \Set]$, with the property that $U$ is equal to the set $\{z \in Aut_{\cal C}(u) \textrm{ | } (\xi, z\circ \xi)\in R_{t}\}$. Notice that, by definition of $R_{t}$, $(\xi, z\circ \xi)\in R_{t}$ if and only if there exists $(f, g)\in R$ and $\chi:cod(f)\to u$ such that $\xi=\chi\circ f$ and $z\circ \xi =\chi\circ g$.  

Summarizing, we have the following result.          

\begin{theorem}\label{atoms}
Under the hypotheses of Theorem \ref{maincategorical}, for any subset $U$ of $Aut_{\cal C}(u)$, $U$ is an open subgroup of $Aut_{\cal C}(u)$ if and only if there exists and object $c$ of $\cal C$ and an equivalence relation $R$ on ${\cal C}(c, -)$ in $[{\cal C}, \Set]$ such that $U=\{z \in Aut_{\cal C}(u) \textrm{ | } z\circ \chi \circ f =\chi\circ g \textrm{ for some } (f,g)\in R \textrm{ and  } \chi:cod(f)\to u\}$. More precisely, for any $\xi:c\to u$, a subset $U\subseteq {\cal I}_{\xi}$ is an open subgroup of $Aut_{\cal C}(u)$ if and only if there exists an equivalence relation $R$ on ${\cal C}(c, -)$ in $[{\cal C}, \Set]$ such that $U=\{z \in Aut_{\cal C}(u) \textrm{ | } \xi=\chi\circ f \textrm{ and } z\circ \xi =\chi\circ g \textrm{ for some } (f,g)\in R \textrm{ and  } \chi:cod(f)\to u\}$.  
\end{theorem}\qed

Let us now suppose that the topology $J_{at}$ on $\cal C$ is subcanonical, i.e. that all the arrows in $\cal C$ are strict monomorphisms. It is natural to wonder whether we can give an elementary description of the dual ${{\cal C}_{at}}$ of the full subcategory of the topos $\Sh({\cal C}^{\textrm{op}}, J_{at})$ on its atoms, and of the canonical embedding $i_{at}:{\cal C}\hookrightarrow {\cal C}_{at}$. We shall call ${\cal C}_{at}$ the \emph{atomic completion} of the category $\cal C$, and we shall say that $\cal C$ is \emph{atomically complete} if the embedding $i_{at}$ is an equivalence. 

Let us denote by $\textbf{At}$ the category whose objects are the categories $\cal D$ with the amalgamation property and in which all arrows are strict monomorphisms, and whose morphisms ${\cal D}\to {\cal D}'$ are the morphisms of sites $({\cal D}, J_{at})\to ({\cal D}', J_{at})$; let us denote by $\textbf{AtC}$ the full subcategory of $\textbf{At}$ on the atomically complete categories. Then the assignment ${\cal C}\to {\cal C}_{at}$ can be made into a functor $C:\textbf{At} \to \textbf{AtC}$; indeed, for any arrow $f:{\cal C}\to {\cal D}$ in $\textbf{At}$, $\Sh(f^{op})^{\ast}:\Sh({\cal C}^{\op}, J_{at})\simeq \Sh({\cal C}_{at}^{\op}, J_{at}) \to \Sh({\cal D}^{\op}, J_{at})\simeq \Sh({\cal D}_{at}^{\op}, J_{at})$ restricts to an arrow ${\cal C}_{at} \to {\cal D}_{at}$ in $\textbf{AtC}$, since it preserves epimorphisms, every object of ${\cal C}_{at}$ is a quotient of an object of $\cal C$ and any quotient of an atom is an atom. In fact, $C$ defines a left adjoint to the canonical inclusion $\textbf{AtC} \hookrightarrow \textbf{At}$, making $\textbf{AtC}$ into a reflective subcategory of $\textbf{At}$.

To define explicitly an extension ${\cal C}_{e}$ of the category $\cal C$ in $\textbf{At}$ which is equivalent to ${\cal C}_{at}$, we observe that any atom of $\Sh({\cal C}^{\textrm{op}}, J_{at})$ is isomorphic to a quotient $yc\slash R$ for an equivalence relation $R$ on $yc$ in $\Sh({\cal C}^{\textrm{op}}, J_{at})$, where $y:{\cal C}^{\textrm{op}}\to [{\cal C}, \Set]$ is the Yoneda embedding. Notice that an equivalence relation $R$ on $yc$ in $\Sh({\cal C}^{\textrm{op}}, J_{at})$ is an equivalence relation $R$ on $yc$ in $[{\cal C}, \Set]$ which is $J_{at}$-closed as a subobject of $yc\times yc$. Concretely, an equivalence relation $R$ on $yc$ in $[{\cal C}, \Set]$ can be identified with a function which assigns to each object $e$ of ${\cal C}$ an equivalence relation $R_{e}$ on the set $Hom_{\cal C}(c, e)$ in such a way that for any arrow $h:e\to e'$ in $\cal C$ and any $(\chi, \xi)\in R_{e}$, $(h\circ \chi, h\circ \xi)\in R_{e'}$; $R$ is an equivalence relation in $\Sh({\cal C}^{\textrm{op}}, J_{at})$ if it is moreover $J_{at}$-closed, that is if it is equal to its $J_{at}$-closure $\overline{R}^{{J_{at}}}$, which is given by: for each $e\in {\cal C}$ $\overline{R}^{J_{at}}_{e}=\{(\chi, \xi) \textrm{ | for some $h:e\to e'$, } (h\circ \chi, h\circ \xi)\in R_{e'}\}$.       

This suggests that a natural choice for the objects of the category ${\cal C}_{e}$ is the pairs $(c, R)$, where $c$ is an object of $\cal C$ and $R$ is a $J_{at}$-closed equivalence relation on $yc$. In order to have an equivalence ${\cal C}_{e}\simeq {\cal C}_{at}$ given by the assignment $(c, R)\to a_{J_{at}}({yc\slash R})$, we need to be able to the define the arrows in ${\cal C}_{e}$ (as well as their composition) in such a way that the arrows $(c', R')\to (c, R)$ in ${\cal C}_{e}$ correspond bijectively with the arrows $a_{J_{at}}({yc\slash R}) \to a_{J_{at}}({yc'\slash R'})$ in the topos $\Sh({\cal C}^{\textrm{op}}, J_{at})$. To an arrow $z:a_{J_{at}}({yc\slash R}) \to a_{J_{at}}({yc'\slash R'})$ we can associate the following set of data: $z$ corresponds, by the universal property of $a_{J_{at}}$ as a left adjoint, to a unique arrow $\tilde{z}:yc\slash R \to a_{J_{at}}({yc'\slash R'})$ in $[{\cal C}, \Set]$; since $z$ is an epimorphism, it is isomorphic to the coequalizer of its kernel pair, which can be identified with an equivalence relation $P$ on $yc$ containing $R$, and we have a unique isomorphism $w:a_{J_{at}}(yc\slash P)\to a_{J_{at}}(yc'\slash R')$ such that $w\circ a_{J_{at}}(\pi_{R, P})=z$, where $\pi_{R, P}:yc\slash R \to yc\slash P$ is the canonical projection arrow. Notice that the equivalence relations on $yc\slash R$ are in bijection with the equivalence relations on $yc$ which contain $R$, via the correspondence which sends any equivalence relation on $yc\slash R$ to its pullback along the canonical arrow $\pi_{R}:yc \to yc\slash R$, and any equivalence relation $z:Z\mono yc\times yc$ on $yc$ which contains $R$ to the image $\tilde{Z}$ of the composite arrow $(\pi_{R}\times \pi_{R})\circ z$.   

As in section \ref{sheafif}, let us consider the pullback 
\[  
\xymatrix {
Q \ar[d]^{\alpha} \ar[r]^{\beta} & yc'\slash R' \ar[d]^{\eta_{yc'\slash R'}}\\
yc\slash P \ar[r]_{\tilde{w}} & a_{J_{at}}(yc'\slash R'), } 
\] 
in the topos $[{\cal C}, \Set]$, where $\tilde{w}$ is the composite of $w$ with the unit $yc\slash P \to a_{J_{at}}(yc\slash P)$ and $\eta_{yc'\slash R'}$ is the unit $yc'\slash R' \to a_{J_{at}}(yc'\slash R')$; clearly, both $a_{J_{at}}(\alpha)$ and $a_{J_{at}}(\beta)$ are isomorphisms, and we have that $a_{J_{at}}(w)=a_{J_{at}}(\beta)\circ {a_{J_{at}}(\alpha)}^{-1}$. Therefore the set of data consisting of the triple $(P, \alpha, \beta)$ completely determines our original arrow $z$, which can be recovered as $a_{J_{at}}(\pi_{R, P})\circ a_{J_{at}}(\beta)\circ {a_{J_{at}}(\alpha)}^{-1}$.    

Notice that, since the above square is a pullback, the arrow $<\alpha, \beta>:Q \to yc\slash P \times yc'\slash R'$ is a monomorphism. Now, for any subobject $<\delta, \gamma>: Q' \mono yc\slash P \times yc'\slash R'$ such that $a_{J_{at}}(\delta)$ and $a_{J_{at}}(\gamma)$ are isomorphisms, the composite $a_{J_{at}}(\gamma)\circ {a_{J_{at}}(\delta)}^{-1}$ yields an arrow $a_{J_{at}}(yc\slash P)\to a_{J_{at}}(yc'\slash R')$ in $\Sh({\cal C}^{\textrm{op}}, J_{at})$ and hence an arrow $a_{J_{at}}(\pi_{R, P})\circ a_{J_{at}}(\gamma)\circ {a_{J_{at}}(\delta)}^{-1}$. Therefore, we can identify the arrows $a_{J_{at}}({yc\slash R}) \to a_{J_{at}}({yc'\slash R'})$ in the topos $\Sh({\cal C}^{\textrm{op}}, J_{at})$ with the triples $(P, \alpha, \beta)$ consisting of a $J_{at}$-closed equivalence relation $P$ on $c$ containing $R$ and a subobject $<\alpha, \beta>: Q \to yc\slash P \times yc'\slash R'$ such that $a_{J_{at}}(\alpha)$ and $a_{J_{at}}(\beta)$ are isomorphisms, modulo the equivalence relation $\simeq$ given by: $(P, \alpha, \beta)\simeq (P, \delta, \gamma)$ if and only if $a_{J_{at}}(\beta)\circ {a_{J_{at}}(\alpha)}^{-1}=a_{J_{at}}(\gamma)\circ {a_{J_{at}}(\delta)}^{-1}$. This equivalence relation $\simeq$ admits an elementary description not involving the functor $a_{J_{at}}$. Indeed, let us consider the pullback     
\[  
\xymatrix {
B_{\beta, \gamma} \ar[d]^{p_{\gamma}} \ar[r]^{p_{\beta}} & Q \ar[d]^{\beta}\\
Q' \ar[r]_{\gamma} & yc'\slash R'} 
\] 
in the topos $[{\cal C}, \Set]$. This diagram is sent by the associated sheaf functor $a_{J_{at}}$ to a pullback in $\Sh({\cal C}^{\textrm{op}}, J_{at})$; therefore, both the arrows $a_{J_{at}}(p_{\beta})$ and $a_{J_{at}}(p_{\gamma})$ are isomorphisms (as they are pullbacks of isomorphisms) and the condition that $a_{J_{at}}(\beta)\circ {a_{J_{at}}(\alpha)}^{-1}=a_{J_{at}}(\gamma)\circ {a_{J_{at}}(\delta)}^{-1}$ can be formulated as the existence of a factorization of $<{a_{J_{at}}(\alpha)}^{-1}, {a_{J_{at}}(\delta)}^{-1}>$ through $<a_{J_{at}}(p_{\beta}), a_{J_{at}}(p_{\gamma})>$; but this is in turn equivalent to the condition that $a_{J_{at}}(\alpha \circ p_{\beta})=a_{J_{at}}(\delta \circ p_{\gamma})$, which admits an elementary formulation (since it is equivalent to the assertion that the equalizer of $\alpha \circ p_{\beta}$ and $\delta \circ p_{\gamma}$ is a $J_{at}$-dense monomorphism - cf. section \ref{sheafif}). 

We thus define an arrow $(c', R')\to (c, R)$ in ${\cal C}_{e}$ as a triple $(P, \alpha, \beta)$ consisting of a $J_{at}$-closed equivalence relation $P$ on $c$ containing $R$ and a subobject $<\alpha, \beta>: Q \mono yc\slash P \times yc'\slash R'$ such that $a_{J_{at}}(\alpha)$ and $a_{J_{at}}(\beta)$ are isomorphisms modulo the equivalence relation $\simeq$ defined above. 

It remains to understand how to define the composition of arrows in ${\cal C}_{e}$ in such a way to make the assignments $(c, R)\to a_{J_{at}}(yc\slash R)$ and $((P, \alpha, \beta):(c', R')\to (c, R)) \to a_{J_{at}}(\pi_{R, P})\circ a_{J_{at}}(\beta)\circ {a_{J_{at}}(\alpha)}^{-1}$ into a (full and faithful) functor $E:{\cal C}_{e}\to {\cal C}_{at}$. Given two (representatives of) arrows $(P, \alpha, \beta):(c', R')\to (c, R)$ and $(P', \alpha', \beta'):(c'', R'')\to (c', R')$ in ${\cal C}_{e}$ corresponding respectively to arrows $u:a_{J_{at}}(yc\slash R)\to a_{J_{at}}(yc'\slash R')$ and $v:a_{J_{at}}(yc'\slash R')\to a_{J_{at}}(yc''\slash R'')$, the kernel pair of the composite arrow $v\circ u$ can be identified with the pullback of the kernel pair of the arrow $v$ along the arrow $u\times u$; but the kernel pair of $v$ can be identified with the equivalence relation on $yc'\slash R'$ corresponding to $P'$, while $u$ is equal to the arrow $a_{J_{at}}(\pi_{R, P})\circ a_{J_{at}}(\beta)\circ {a_{J_{at}}(\alpha)}^{-1}$.

\[  
\xymatrix {
P \ar@<1ex>[dr] \ar[dr] & & & \\
R \ar[u] \ar[r] \ar@<1ex>[r] & yc \ar[ddr] \ar[r] & a_{J_{at}}(yc\slash R) \ar[dd]^{u} \ar[ddr]^{a_{J_{at}}(\pi_{R, P})} & \\
P' \ar[dr] \ar@<1ex>[dr]  & & & \\
R' \ar[u] \ar[r] \ar@<1ex>[r] & yc' \ar[ddr] \ar[r] & a_{J_{at}}(yc'\slash R') \ar[dd]^{v} \ar[ddr]^{a_{J_{at}}(\pi_{R', P'})} & a_{J_{at}}(yc\slash P) \ar[l]^{\cong} \\
& & & \\
& yc'' \ar[r] & a_{J_{at}}(yc'' \slash R'')  & a_{J_{at}}(yc'\slash P') \ar[l]^{\cong}} 
\]

Let us thus define $P''$ to be the pullback along the canonical arrow $yc \to yc\slash R$ of the equivalence relation on $yc\slash R$ generated by the arrow $h:H\to yc\slash R \times yc\slash R$ defined by the following diagram, where the upper and the lower square are pullbacks in $[{\cal C}, \Set]$:

\[  
\xymatrix {
H \ar[rr]^{h} \ar[d]^{r} & & yc\slash R \times yc\slash R \ar[d]^{\pi_{R, P}\times \pi_{R, P}}\\
L \ar[rr]^{<\alpha \circ \chi, \beta \circ \xi>    \quad} & & yc\slash P \times yc\slash P\\
L \ar[u]^{1_{L}} \ar[rr]^{<\chi, \xi>} \ar[d]^{s} & & Q \times Q \ar[u]_{\alpha\times \alpha} \ar[d]^{\beta\times \beta}\\ 
\tilde{P'} \ar[rr] & & yc'\slash R' \times yc'\slash R'.} 
\]     

The relation $P''$ will be the first component of a triple representing the composition of the two arrows $(P, \alpha, \beta):(c', R')\to (c, R)$ and $(P', \alpha', \beta'):(c'', R'')\to (c', R')$ in the category ${\cal C}_{e}$. To define the remaining two components, we observe that the arrow $m:a_{J_{at}}(yc\slash P'')\to a_{J_{at}}(yc''\slash R'')$ in $\Sh({\cal C}^{\textrm{op}}, J_{at})$ corresponding to the composition $v\circ u$ can be expressed as the composition of arrows $a_{J_{at}}(\beta')\circ {a_{J_{at}}(\alpha')}^{-1} \circ q$, where $q$ is the canonical arrow from the coequalizer of $a_{J_{at}}(h)$ to the coequalizer of $a_{J_{at}}(\tilde{P'})$. Now, if we denote by $b_{1}:yc\slash P''\to C$ (resp. by $b_{2}:C'\to C$, by $b_{3}:C'\to yc'\slash P'$) the canonical arrows between the coequalizers of $h$ and $<\alpha \circ \chi, \beta \circ \xi>$ (resp. of $<\chi, \xi>$ and $<\alpha \circ \chi, \beta \circ \xi>$, of $<\chi, \xi>$ and $\tilde{P'}$) then the arrow $q$ is equal to the composite $a_{J_{at}}(b_{3})\circ {a_{J_{at}}(b_{2})}^{-1}\circ a_{J_{at}}(b_{1})$. Notice that $a_{J_{at}}(b_{2})$ is an isomorphism since $a_{J_{at}}(\alpha)$ is an isomorphism, and $a_{J_{at}}(b_{3})$ is an isomorphism since $a_{J_{at}}(\beta)$ is an isomorphism; as $q$ is an isomorphism (by definition of $P''$), it follows that $a_{J_{at}}(b_{1})$ is also an isomorphism. Now, consider the pullbacks
\[  
\xymatrix {
M \ar[d]^{c} \ar[r]^{a} & C' \ar[d]^{b_{2}}\\
yc\slash P'' \ar[r]_{b_{1}} & C} 
\]
and
\[  
\xymatrix {
N \ar[d]^{b} \ar[r]^{d} & M \ar[d]^{a}\\
W \ar[r]_{y} \ar[d]^{x} & C' \ar[d]^{b_{3}}\\
Q' \ar[r]_{\alpha'} & yc'\slash P'} 
\]      
in the topos $[{\cal C}, \Set]$.

Notice that $a_{J_{at}}(c)$ is an isomorphism since $a_{J_{at}}(b_{2})$ is; similarly, $a_{J_{at}}(y)$ and $a_{J_{at}}(d)$ are isomorphisms as they are pullbacks of isomorphisms. Let us set $\alpha''_{p}=c\circ d$ and $\beta''_{p}=\beta' \circ x \circ b$, and verify that $m=a_{J_{at}}(\beta''_{p})\circ {a_{J_{at}}(\alpha''_{p})}^{-1}$. We have $m=a_{J_{at}}(\beta')\circ {a_{J_{at}}(\alpha')}^{-1}\circ q=a_{J_{at}}(\beta')\circ {a_{J_{at}}(\alpha')}^{-1} \circ a_{J_{at}}(b_{3})\circ {a_{J_{at}}(b_{2})}^{-1}\circ a_{J_{at}}(b_{1})=a_{J_{at}}(\beta')\circ a_{J_{at}}(x)\circ {a_{J_{at}}(y)}^{-1}\circ a_{J_{at}}(a)\circ {a_{J_{at}}(c)}^{-1}=a_{J_{at}}(\beta')\circ a_{J_{at}}(x)\circ a_{J_{at}}(b)\circ {a_{J_{at}}(d)}^{-1}\circ {a_{J_{at}}(c)}^{-1}=a_{J_{at}}(\beta''_{p})\circ {a_{J_{at}}(\alpha''_{p})}^{-1}$, as required. Now, we define $<\alpha'', \beta''>$ as the image of $<\alpha''_{p}, \beta''_{p}>$ in $[{\cal C}, \Set]$; clearly, $a_{J_{at}}(\beta'')\circ {a_{J_{at}}(\alpha'')}^{-1}=a_{J_{at}}(\beta''_{p})\circ {a_{J_{at}}(\alpha''_{p})}^{-1}=m$, and hence $(P'', \alpha'', \beta'')$ defines an arrow $(c'', R)\to (c, R)$ in ${\cal C}_{e}$ which is the composition of the arrows $(P, \alpha, \beta):(c', R')\to (c, R)$ and $(P', \alpha', \beta'):(c'', R'')\to (c', R')$ in ${\cal C}_{e}$.  

Note that the definition of the composition operation in the category ${\cal C}_{e}$ is completely elementary. We thus have a category ${\cal C}_{e}$ (the categorical axioms are easily verified at this point) and a functor ${\cal C}_{e}\to {\cal C}_{at}$ which is full and faithful and essentially surjective, i.e. a (weak) equivalence. The following theorem records this fact.

\begin{theorem}
Let $\cal C$ be a category in $\textbf{At}$. Then the category ${\cal C}_{e}$ defined above is (weakly) equivalent to the atomic completion ${\cal C}_{at}$ of $\cal C$.
\end{theorem}\qed 

Note that the definition of the category ${\cal C}_{e}$ is relatively complicated since, in order to be equivalent to ${\cal C}_{at}$, one has to characterize in `elementary terms' (i.e., in a way which involves only constructions in the presheaf topos $[{\cal C}, \Set]$, and not the associated sheaf functor $a_{J_{at}}$) \emph{all} the arrows $a_{J_{at}}(yc\slash R)\to a_{J_{at}}(yc'\slash R')$ in the topos $\Sh({\cal C}^{\textrm{op}}, J_{at})$. If instead one restricts to the arrows $a_{J_{at}}(yc\slash R)\to a_{J_{at}}(yc'\slash R')$ of the form $a_{J_{at}}(h)$ for an arrow $h:yc\slash R \to yc'\slash R'$ in $[{\cal C}, \Set]$, one obtains a subcategory ${\cal C}_{r}$ admitting a simpler description. Let us define the objects of ${\cal C}_{r}$ as the pairs $(c, R)$, where $c$ is an object of $\cal C$ and and $R$ is a $J_{at}$-closed equivalence relation on $yc$, and the arrows $(c', R')\to (c, R)$ in ${\cal C}_{r}$ to be the arrows $f:c\to c'$ in $\cal C$ such that the arrow $yf:yc\to yc'$ factorizes through the projection arrows $yc\to yc\slash R$ and $yc'\to yc'\slash R'$, modulo the equivalence relation which identifies two arrows $f, f':c\to c'$ when the images of the factorizations of $yf$ and $yf'$ through the projection arrows under the associated sheaf functor $a_{J_{at}}$ are equal (note that this condition admits an elementary formulation).   
Notice that to every object $c$ of $\cal C$ we can associate an object of ${\cal C}_{e}$, namely the pair $(c, R_{c})$, where $R_{c}$ is the empty relation on $yc$, and every arrow $h:c'\to c$ in $\cal C$ gives canonically rise to an arrow $(c', R_{c'})\to (c, R_{c})$ in ${\cal C}_{e}$, given by the triple $(P, \alpha, \beta)$, where $P$ is the kernel pair of $yh$ in $[{\cal C}, \Set]$, $\alpha$ is the identity on $yc\slash P$ and $\beta$ is equal to the factorization of $yh:yc \to yc'$ through the canonical projection $yc \to yc\slash P$. In fact, these assignments yield a functor ${\cal C}\to {\cal C}_{r}$ which, composed with the canonical embedding ${\cal C}_{r} \hookrightarrow {\cal C}_{e}$, yields a functor ${\cal C}\to {\cal C}_{e}$ which is isomorphic to the canonical functor ${\cal C}\to {\cal C}_{at}$.

The following result gives a criterion for recognizing atomically complete categories. 

\begin{theorem}\label{atomcompl}
Let $\cal C$ be a category in $\mathbf{At}$. Then the following conditions are equivalent:
\begin{enumerate}[(i)]
\item $\cal C$ is atomically complete;

\item For every object $c$ of $\cal C$ and any equivalence relation $R$ on $yc$ in $[{\cal C}, \Set]$, there exists an arrow $m:d\to c$ such that for any arrows $f,g:c\to e$ in $\cal C$, $f\circ m=g\circ m$ if and only if there exists $h:e\to e'$ such that $(h\circ f, h\circ g)\in R_{e'}$;

\item The category $\cal C$ has equalizers, for any object $c$ of $\cal C$ there exist arbitrary intersections of subobjects of $c$ in $\cal C$, and for any pair of arrows $h, k:c\to e$ in $\cal C$ with equalizer $m:d\mono c$ we have that for any pair of arrows $l,n:c\to e'$, $l\circ m=n\circ m$ if and only if there exist an arrow $s:e'\to e''$ such that $(s\circ l, s\circ n)$ belongs to the equivalence relation on $Hom_{\cal C}(c, e'')$ generated by the relation consisting of the pairs of the form $(t\circ h, t\circ k)$ for an arrow $t:e \to e''$. 
\end{enumerate}  
\end{theorem}

\begin{proofs}
The equivalence between $(i)$ and $(ii)$ was proved above as Proposition \ref{dualatom}; we shall prove that $(i)$ implies $(iii)$ and that $(iii)$ implies $(ii)$.

$(i)\imp (iii)$ Let $\cal D$ be the full subcategory of the topos $\Sh({\cal C}^{\textrm{op}}, J_{at})$ on its atoms; if $\cal C$ is atomically complete then it is dually equivalent to $\cal D$. The condition that $\cal C$ has equalizers is thus equivalent to the condition that $\cal D$ has coequalizers, and this is obviously satisfied since $\cal D$ is closed in $\Sh({\cal C}^{\textrm{op}}, J_{at})$ under coequalizers, as epimorphic images of atoms are atoms. Concerning the property that arbitrary intersection of subobjects exist in $\cal C$, this is clearly equivalent to the condition that for any object $d$ of $\cal D$ and any family of arrows $\{s_{i}:d\to d_{i} \textrm{ | } i\in I\}$ in $\cal D$ there exists an arrow $s:d\to e$ which factors (necessarily uniquely, as all arrows in $\cal D$ are epimorphisms) through each of the $s_{i}$ and such that any arrow $s':d\to e'$ with the same property factors through $s$. To prove this property of $\cal D$, we use the fact that for any $i\in I$, the arrow $s_{i}:d\to d_{i}$ is isomorphic to the quotient arrow $d\to d\slash R_{i}$, where $R_{i}\mono d\times d$ is the kernel pair of $s_{i}$; if we denote by $R\mono d\times d$ the union of the subobjects $R_{i}\mono d\times d$ in the topos $\Sh({\cal C}^{\textrm{op}}, J_{at})$ then the coequalizer $s:d\to d\slash R$ is an arrow in $\cal D$ (since $d\slash R$ is a quotient of an atom, namely $d$, and hence it is itself an atom) which is easily seen to satisfy the required property. Let $h,k:c\to e$ two arrows in $\cal C$; since $\cal C$ is atomically complete then condition $(ii)$ holds for the equivalence relation $R_{h,k}$ on $yc$ generated by the image of the arrow $<yh, yk>:ye\to yc\times yc$, yielding an arrow $m:d\to c$ with the property that for any arrows $f,g:c\to e$ in $\cal C$, $f\circ m=g\circ m$ if and only if there exists $h:e\to e'$ such that $(h\circ f, h\circ g)\in R_{e'}$; by using the explicit description of the relation $R_{h,k}$ and the fact that all the arrows in $\cal C$ are monomorphisms, one immediately realizes that one can suppose without loss of generality that $m$ is the equalizer of $h$ and $k$ in $\cal C$; condition $(iii)$ is thus satisfied. 

$(iii)\imp (ii)$ Given an equivalence relation $R$ on $yc$ in $\Sh({\cal C}^{\textrm{op}}, J_{at})$, since the topos $\Sh({\cal C}^{\textrm{op}}, J_{at})$ is atomic, $R$ can be expressed as a (disjoint) union of relations $R_{i}\mono yc\times yc$ which are atoms of the topos $\Sh({\cal C}^{\textrm{op}}, J_{at})$; every such relation $R_{i}$ is of the form $R_{h,k}$ for some arrows $h,k:c\to e$ in $\cal C$, as there exists an arrow $ye\to R_{i}$ in $\Sh({\cal C}^{\textrm{op}}, J_{at})$, which is necessarily an epimorphism. Now, condition $(iii)$ ensures that every relation $R_{i}$ admits a coequalizer of the form $ym_{i}:yc \to yd_{i}$ for an arrow (in fact, a monomorphism) $m_{i}:d_{i}\to c$; but from the above discussion we know that the arrow $ym:tc\to ydom(m)$, where $m:dom(m)\mono c$ is the intersection of all the subobjects $m_{i}:d_{i}\mono c$ in $\cal C$, is a coequalizer of the relation $R$ in $\Sh({\cal C}^{\textrm{op}}, J_{at})$. Condition $(ii)$ is thus satisfied.        
\end{proofs}

\begin{remark}
\begin{enumerate}[(a)]
\item It is clear from the proof of the theorem that in condition $(ii)$ one can suppose without loss of generality $m$ to be the intersection of the equalizers of the pairs $(f,g)\in R_{e}$ for some $e\in {\cal C}$.

\item If $\cal C$ has equalizers and arbitrary intersections of subobjects then the embedding ${\cal C}\hookrightarrow {\cal C}_{at}$ is an equivalence (i.e., $\cal C$ is atomically complete) if and only if it preserves equalizers. Indeed, by condition $(ii)$ in the statement of the theorem, $\cal C$ is atomically complete if and only if for any arrows $f,g:a\to c$ in $\cal C$ with equalizer $m:d\to a$, the arrow $ym$ is a coequalizer of the arrows $yf$ and $yg$ in $\Sh({\cal C}^{\textrm{op}}, J_{at})$, where $y:{\cal C}^{\textrm{op}}\to [{\cal C}, \Set]$ is the Yoneda embedding. Now, as ${\cal C}_{at}$ is atomically complete, for any arrows $f,g:a\to c$ in $\cal C$ with equalizer $m:d\to a$, the arrow $y'm$ is a coequalizer of the arrows $y'f$ and $y'g$ in $\Sh({{\cal C}_{at}}^{\textrm{op}}, J_{at})$, where $y':{\cal C}_{at}\to [{{\cal C}_{at}}^{\textrm{op}}, \Set]$ is the Yoneda embedding; but $\Sh({{\cal C}_{at}}^{\textrm{op}}, J_{at})\simeq \Sh({\cal C}^{\textrm{op}}, J_{at})$, and under this equivalence every representable of the form $y'e$ (for $e\in {\cal C}$) corresponds to the representable $ye$, and every arrow of the form $y'k$ for $k$ in $\cal C$ is sent to the arrow $yk$. From which it follows that if the embedding ${\cal C}\hookrightarrow {\cal C}_{at}$ preserves equalizers then $ym$ is the equalizer of $yf$ and $yg$ in $\Sh({\cal C}^{\textrm{op}}, J_{at})$. Since this argument holds for arbitrary $f,g$ we can conclude that the fact that the embedding ${\cal C}\hookrightarrow {\cal C}_{at}$ preserves equalizers implies the atomic completeness of $\cal C$.         
\end{enumerate}
\end{remark}

The following result provides some properties of atomically complete categories. Recall that a multilimit for a diagram $D:{\cal I}\to {\cal C}$ in a category $\cal C$ is a set $T$ of cones on $D$ such that for every cone $x$ over $D$ there exists a cone $t$ in $T$ such that $x$ factors through $t$; the notion of multicolimit is the dual one.

\begin{theorem}
Let $\cal C$ be an atomically complete category. Then
\begin{enumerate}[(i)]
\item Let $D:{\cal I}\to {\cal C}$ be a diagram defined on a small non-empty connected category $\cal I$. Then $D$ has a limit in $\cal C$;

\item Let $D:{\cal I}\to {\cal C}$ be a diagram defined on a small category $\cal I$. Then $D$ has a multicolimit in $\cal C$.
\end{enumerate}
\end{theorem}

\begin{proofs}
Let $\cal D$ be the full subcategory of the topos $\Sh({\cal C}^{\textrm{op}}, J_{at})$ on its atoms. Then $\cal C$ is dually equivalent to $\cal D$. 

$(i)$ To prove condition $(i)$ we show that every diagram $D:{\cal I}\to {\cal D}$ defined on a small non-empty connected category $\cal I$ has a colimit in $\cal D$. Let us consider the colimit $colim(D)$ of $D$ in the topos $\Sh({\cal C}^{\textrm{op}}, J_{at})$. Since the topos $\Sh({\cal C}^{\textrm{op}}, J_{at})$ is atomic, the object $colim(D)$ can be expressed as a coproduct of atoms $A_{k}$ (for $k\in K$); let us denote by $j_{k}:A_{k}\to colim(D)$ the canonical coproduct arrows. For each $i\in {\cal I}$ we have a colimit arrow $s_{i}:D(i)\to colim(D)$ in $\Sh({\cal C}^{\textrm{op}}, J_{at})$; the fact that $D(i)$ is an atom and that coproducts are stable under pullback in a topos implies that there exists a unique $k_{i}\in K$ such that $s_{i}$ factors through $j_{k_{i}}$. Notice that if there is an arrow $f:i\to i'$ in $\cal I$ then we have an arrow $D(f):D(i)\to D(i')$ in $\Sh({\cal C}^{\textrm{op}}, J_{at})$ such that $s_{i'}\circ D(f)=s_{i}$ and hence $k_{i}=k_{i'}$. Therefore, since $\cal I$ is connected, there exists a element $k\in K$ such that all the arrows $s_{i}$ factor through $j_{k}$. Now, by the universal property of the colimit, the arrows $s_{i}$ are jointly epimorphic, whence the coproduct arrow $j_{k}:A_{k}\to colim(D)$ is an epimorphism; but this implies that $colim(D)\cong A_{k}$, as $j_{k}$ is also a monomorphism; in other words, the cone on $D$ in $\cal D$ with vertex $A_{k}$ given by the factorizations of the arrows $s_{i}$ through $j_{k}$ is a colimiting cone for $D$ in $\cal D$.          

$(ii)$ To prove condition $(ii)$ we show that every diagram $D:{\cal I}\to {\cal D}$ be a diagram defined on a small category $\cal I$ has a multilimit in $\cal D$. Let us consider the limit $lim(D)$ of the diagram $D$ in the topos $\Sh({\cal C}^{\textrm{op}}, J_{at})$, and decompose the object $lim(D)$ as a coproduct of atoms $A_{k}$ (for $k\in K$). Then the composition of the coproduct arrows $A_{k}\to lim(D)$ with the limit arrows $lim(D)\to D(i)$ defines a set of cones in $\cal D$ on the diagram $D$ with vertexes $A_{k}$; to verify that this family of cones defines a multilimit for $D$ in $\cal D$ it suffices to observe that for any cone $(B, \{t_{i}:B\to D(i) \textrm{ | } i\in I\}$ over $D$ in $\cal D$ the arrow $B\to lim(D)$ induced by the universal property of the limit $lim(D)$ factors through exactly one of the coproduct arrows $A_{k}\mono A$, as $B$ is an atom and coproducts are stable under pullback in a topos.  
\end{proofs}

In the following theorem, we denote by $Subgr(Aut_{\cal D}(u))$ the preorder category consisting of the subgroups of $Aut_{\cal D}(u)$.

\begin{theorem}\label{galoiscomplete}
Let $\cal C$ be an atomically complete category and ${\cal C}\hookrightarrow {\cal D}$ be an embedding of $\cal C$ into a category $\cal D$ containing a $\cal C$-universal and $\cal C$-ultrahomogeneous $u$ such that there is an equivalence
\[
\Sh({\cal C}^{\textrm{op}}, J_{at})\simeq \Cont(Aut_{\cal D}(u))
\]
induced by a functor $F:{\cal C}^{\textrm{op}}\slash u\to Subgr(Aut_{\cal C}(u))$ sending to any object $\chi:d\to u$ of ${\cal C}\slash u$ the subgroup ${\cal I}_{\chi}:=\{f:u\cong u \textrm{ | } f\circ \chi=\chi\}$ of $Aut_{\cal C}(u)$ (as provided for instance by Theorem \ref{maincategorical} or by Theorem \ref{catd}). Then $F$ induces a bijection between the isomorphism classes of objects of ${\cal C}\slash u$ and the open subgroups of $Aut_{\cal D}(u)$; for any open subgroup $U$ of $Aut_{\cal D}(u)$ the arrow $\chi_{U}:c_{U}\to u$ in $\cal D$ corresponding to it in this bijection satisfies the following universal property: it is fixed by all the automorphisms in $U$ and any other arrow $a\to u$ in $\cal D$ which is fixed by all the automorphisms in $U$ factors uniquely through it.  
\end{theorem}\qed

\begin{proposition}
In Theorem \ref{galoiscomplete}, if the subgroup $U$ corresponds to a pair $(\xi:c\to u, R)$, where $R$ is an equivalence relation on $yc$ in $[{\cal C}, \Set]$ as in Theorem \ref{atoms} (i.e., $U=\{z \in Aut_{\cal C}(u) \textrm{ | } \xi=\chi\circ f \textrm{ and } z\circ \xi =\chi\circ g \textrm{ for some } (f,g)\in R \textrm{ and  } \chi:cod(f)\to u\}$) then the arrow $\chi_{U}:c_{U}\to u$ can be identified with the composite $\xi\circ m$ of $\xi$ with the monomorphism $m:d\mono c$ given by the intersection of the equalizers of the pairs $(f,g)\in R_{e}$ for some $e\in {\cal C}$.     
\end{proposition}

\begin{proofs}
First, we notice that the existence of the equalizers and intersection of subobjects in $\cal C$ is guaranteed by Theorem \ref{atomcompl}. 

To prove that $\xi \circ m$ is isomorphic to $\chi_{U}$ in ${\cal C}\slash u$, we exploit the universal property of $\chi_{U}$ given by Theorem \ref{galoiscomplete}; that is, we verify that
$(1)$ $\xi\circ m$ is fixed by all the automorphisms in $U$;
$(2)$ for any arrow $\beta:b\to u$ which is fixed by all the automorphisms in $U$, $\beta$ factors uniquely through $\xi\circ m$. 

To show $(1)$, we have to verify that for every $z\in U$, $z\circ \xi \circ m= \xi\circ m$. By definition of $U$, there exists a pair $(f,g)\in R_{e}$ and an arrow $\chi:e\to u$ such that $(\xi, z\circ \xi)=(\chi\circ f, \chi\circ g)$. Therefore, as $f\circ m=g\circ m$ by definition of $m$, we have that $z\circ \xi \circ m= \xi\circ m$, as required.

To show $(2)$, we remark that the uniqueness of the factorization of $\beta$ through $\xi$ is guaranteed by condition $(b)$ in Theorem \ref{catd} (cf. Remark \ref{necessarycondition}(a)). To prove its existence, we argue as follows. First, we note that $\beta$ factorizes (uniquely) through $\xi$, say $\beta=\xi\circ w$; indeed, ${\cal I}_{\xi}\subseteq U \subseteq {\cal I}_{\beta}$. Next, we observe that $w$ factors (uniquely) through $m$. To prove this, we verify that for any $(f,g)\in R_{e}$, $w\circ f=w\circ g$. Now, the fact that $u$ is $\cal C$-universal, $\cal C$-homogeneous and $\cal C$-ultrahomogeneous implies that there exists an arrow $\gamma:e\to u$ and an automorphism $z$ of $u$ such that $\gamma \circ f=\xi$ and $z\circ \gamma \circ g=\xi$. From these identities it follows at once that $z^{-1}\in U$, which in turn implies that $z^{-1}\circ \beta=\beta$; but $\beta=\xi \circ w$ and $\xi=\gamma \circ f=z\circ \gamma \circ g$, whence $z^{-1}\circ (z\circ \gamma \circ g) \circ w=\gamma \circ f\circ w$. As $\gamma$ satisfies condition $(b)$ in Theorem \ref{catd}, we conclude that $w\circ f=w\circ g$, as required.    
\end{proofs}

For any atomically complete category $\cal C$, the topos $\Sh({\cal C}^{\textrm{op}}, J_{at})$ admits a simple characterization. 

\begin{proposition}\label{freecoproduct}
Let $\cal C$ be an atomically complete category. Then the canonical embedding ${\cal C}^{\textrm{op}}\hookrightarrow \Sh({\cal C}^{\textrm{op}}, J_{at})$ realizes $\Sh({\cal C}^{\textrm{op}}, J_{at})$ as the free small-coproduct completion of ${\cal C}^{\textrm{op}}$.
\end{proposition}

\begin{proofs}
It suffices to verify that the embedding ${\cal C}^{\textrm{op}}\hookrightarrow \Sh({\cal C}^{\textrm{op}}, J_{at})$ satisfies the universal property of the free small-coproduct completion of the category ${\cal C}^{\textrm{op}}$. We note that ${\cal C}^{\textrm{op}}$ is equivalent to the full subcategory of the topos $\Sh({\cal C}^{\textrm{op}}, J_{at})$ on its atoms. The thesis follows from the fact that, as $\Sh({\cal C}^{\textrm{op}}, J_{at})$ is an atomic topos, every object $A$ of it can be written, in a unique way, as a disjoint union (=coproduct) of atoms $A_{k}$; also, since coproducts in a topos are stable under pullback, for any arrow $f:A\to B$ in $\Sh({\cal C}^{\textrm{op}}, J_{at})$, where $A=\coprod_{k\in K}A_{k}$ and $B=\coprod_{j\in J}B_{j}$ with the $A_{k}$ and the $B_{j}$ atoms, there exists a unique function $s_{f}:K\to J$ and for any $k\in K$ a unique arrow $f_{k}:A_{k}\to B_{s_{f}(k)}$ such that $f\circ \alpha_{k}=\beta_{s_{f}(k)}\circ f_{k}$ for any $k$, where $\alpha_{k}:A_{k}\to A$ and $\beta_{j}:B_{j}\to B$ are the canonical coproduct arrows. Any functor $F:{\cal D}\to {\cal B}$ from $\cal D$ to a category with small coproducts $\cal B$ can thus be extended, in a unique way, to a small coproduct-preserving functor $\tilde{F}:\Sh({\cal C}^{\textrm{op}}, J_{at})\to {\cal B}$. Indeed, for any object $A$ of $\Sh({\cal C}^{\textrm{op}}, J_{at})$, we set $\tilde{F}(A)$ equal to the coproduct $\coprod_{k\in K}F(A_{k})$ and for any arrow $f:A\to B$ of $\Sh({\cal C}^{\textrm{op}}, J_{at})$ we set $\tilde{F}(f)$ equal to the unique arrow $h:\coprod_{k\in K}F(A_{k}) \to \coprod_{j\in J}F(B_{j})$ in $\cal B$ such that for any $k\in K$, $h\circ \xi_{k}=\chi_{s_{f}(k)} \circ F(f_{k})$, where $\xi_{k}:F(A_{k})\to \coprod_{k\in K}F(A_{k})$ and $\xi_{j}:F(B_{j})\to \coprod_{j\in J}F(B_{j})$ are the canonical coproduct arrows. 
\end{proofs}

\section{Other insights from the `bridge' technique}

\subsection{Irreducibility and discreteness}\label{discreteness}

Let us analyze the notion of irreducible object in the context of a Galois-type equivalence
\[
\Sh({\cal C}^{\textrm{op}}, J_{at})\simeq \Cont(Aut_{\cal C}(u))
\]
provided by Theorem \ref{maincategorical} or Theorem \ref{catd}.

In a general Grothendieck topos, an object is said to be irreducible if any epimorphic sieve on it contains the identity; in a Boolean topos, an object is irreducible if and only if it is an atom and every epimorphism to it is an isomorphism. Therefore, for any topological group $G$, the irreducible objects of the topos $\Cont(G)$ are precisely the $G$-sets of the form $G\slash U$ where $U$ is an open subgroup which does not contain any proper open subgroup. The irreducible objects of the topos $\Cont(Aut_{\cal C}(u))$ are thus the open sets of the form ${\cal I}_{\chi}$ such that for any $\xi:c\to u$, ${\cal I}_{\xi}\subseteq {\cal I}_{\chi}$ implies ${\cal I}_{\xi}={\cal I}_{\chi}$ (since any open subgroup of $Aut_{\cal C}(u)$ contains a subgroup of the form ${\cal I}_{\chi}$). This implies, in view of our Morita-equivalence, that the irreducible objects of the topos $\Sh({\cal C}^{\textrm{op}}, J_{at})$ are all of the form $l(c)$ (for $c\in {\cal C}$). If $J_{at}$ is subcanonical then we can further characterize them as follows: these are the objects of the form $l(c)$ with the property that any arrow in ${\cal C}$ with domain $c$ is an isomorphism. 

Next, let us consider a strengthening of the property of being irreducible: the property of being irreducible and a generator for the topos. For any topological group $G$, the irreducible objects which generate the topos $\Cont(G)$ are precisely the open subgroups which do not contain any proper open subgroups and which are contained in some conjugate of any open subgroup of $G$. If $J_{at}$ is subcanonical, the irreducible generators for the topos $\Sh({\cal C}^{\textrm{op}}, J_{at})$ are precisely those of the form ${\cal I}_{\chi}$ where $dom(\chi)$ is an object such that any arrow in $\cal C$ with domain $dom(\chi)$ is an isomorphism and any object of $\cal C$ admits an arrow to $dom(\chi)$.  

Summarizing, we have the following characterization result.

\begin{proposition}\label{propirr}
Let $({\cal C}, u)$ be a pair satisying the hypotheses of Theorem \ref{maincategorical} such that every arrow of $\cal C$ is a strict monomorphism. Then
\begin{enumerate}[(i)]
\item The open subgroups of $Aut_{\cal C}(u)$ which do not contain any proper open subgroup are precisely the ones of the form ${\cal I}_{\chi}$ where $dom(\chi)$ is an object such that any arrow in $\cal C$ with domain $dom(\chi)$ is an isomorphism; 
\item The open subgroups which moreover are contained, up to conjugation, in any open subgroup of $Aut_{\cal C}(u)$ are those of the form ${\cal I}_{\chi}$ where $dom(\chi)$ is an object such that any arrow in $\cal C$ with domain $dom(\chi)$ is an isomorphism and any object of $\cal C$ admits an arrow to $dom(\chi)$.
\end{enumerate}    
\end{proposition}\qed

The following proposition expresses a link between the concept of irreducible generator in a topos and the property of discreteness of a topological group. 

\begin{proposition}
Under the hypotheses of Proposition \ref{propirr}, the following conditions are equivalent.
\begin{enumerate}[(i)]
\item There exists an open subgroup of $Aut_{\cal C}(u)$ which does not contain any proper open subgroup and which is contained, up to conjugation, in any open subgroup of $Aut_{\cal C}(u)$; 

\item There exists an object $c$ of $\cal C$ with the property that any arrow in $\cal C$ with domain $dom(\chi)$ is an isomorphism and any object of $\cal C$ admits an arrow to $c$;

\item The topological group $Aut_{\cal C}(u)$ is Morita-equivalent to a discrete group.
\end{enumerate}

\end{proposition}

\begin{proofs}
The equivalence of $(i)$ and $(ii)$ has been established above. The fact that $(iii)$ implies $(i)$ is obvious. To prove that $(ii)$ implies $(iii)$ we observe that, if $c$ is an object of $\cal C$ with the property that any arrow in $\cal C$ with domain $dom(\chi)$ is an isomorphism and any object of $\cal C$ admits an arrow to $c$ then the full subcategory $Aut{\cal C}(c)$ of ${\cal C}^{\textrm{op}}$ on $c$ is $J_{at}$-dense whence the Comparison Lemma yields an equivalence $\Sh({\cal C}^{\textrm{op}}, J_{at})\simeq \Cont(Aut_{\cal C}(c))$, where $Aut_{\cal C}(c)$ is endowed with the discrete topology. 
\end{proofs}

Notice that if the group $Aut_{\cal C}(u)$ is discrete then $u\cong c$; indeed, since groups are (trivially) Cauchy-complete categories, we have an isomorphism $Aut_{\cal C}(u)\simeq Aut_{\cal C}(c)$, from which it follows (cf. section \ref{univ} below) that $u\cong c$. 

In other words, if $\cal C$ contains an object $c$ with the above-mentioned property then any $\cal C$-ultrahomogeneous and $\cal C$-universal object in $\Ind{\cal C}$ is isomorphic to $c$. Notice that $c$ is $\cal C$-ultrahomogeneous (this follows from the equivalence $\Sh({\cal C}^{\textrm{op}}, J_{at})\simeq \Cont(Aut_{\cal C}(c))$ or can also be proved directly by using the fact that $\cal C$ satisfies the amalgamation property).   

The following result represents a discrete `Galois-type' theorem.

\begin{theorem}\label{discreteGalois}
Let $\cal C$ be a category satisfying AP which contains an object $c$ with the property that all arrows $c\to d$ in $\cal C$ are isomorphisms and for any $e\in {\cal C}$ there exists an arrow $e\to c$ in $\cal C$. Then $c$ is $\cal C$-ultrahomogeneous and we have an equivalence $\Sh({\cal C}^{\textrm{op}}, J_{at})\simeq \Cont(Aut_{\cal C}(c))$, where $Aut_{\cal C}(c)$ is the discrete group of automorphisms of $c$ in $\cal C$. Moreover, if all the arrows of $\cal C$ are strict monomorphisms then the functor $F:{\cal C}^{\textrm{op}}\slash c\to Subgr(Aut_{\cal C}(c))$ sending to any object $\chi:d\to c$ of $\cal C$ the subgroup $\{f:c\cong c \textrm{ | } f\circ \chi=\chi\}$ of $Aut_{\cal C}(c)$ yields a bijection between the isomorphism classes of objects of ${\cal C}\slash c$ and the subgroups of $Aut_{\cal C}(c)$ in the image of the functor $F$.
\end{theorem}\qed

\begin{remark}
The second part of the theorem (that is, the particular case of it for subcanonical sites) is essentially equivalent to the discrete Galois theory of \cite{dubuc3}.
\end{remark}

\subsection{Homogeneous structures, completions and torsors}

By considering the invariant notion of point of a topos in the context of the Morita-equivalence of Theorem \ref{main}, and recalling that any continuous (resp. cocontinuous) functors between Grothendieck toposes has a left adjoint (resp. a right adjoint), we immediately obtain the following result.

\begin{theorem}
Any $\cal C$-homogeneous object of $\Ind{\cal C}$, regarded as a functor ${\cal C}^{\textrm{op}}\to \Set$, can be extended via $F:{\cal C}^{\textrm{op}}\to \Cont(Aut_{\cal C}(u))$, uniquely up to isomorphism, to a cartesian cocontinuous functor $\Cont(Aut_{\cal C}(u))\to \Set$; conversely, any such functor restricts, via $F$, to a $\cal C$-homogeneous object of $\Ind{\cal C}$.  
\end{theorem}\qed

In the discrete case, that is given a Morita-equivalence of the form 
\[
\Sh({\cal C}^{\textrm{op}}, J_{at})\simeq \Cont(Aut_{\cal C}(c)),
\]
the situation simplifies and we obtain the following characterization of homogeneous structures in terms of torsors: the $\cal C$-homogeneous structures can be identified precisely with the $Aut_{\cal C}(c)$-torsors.

Notice that in general, by Diaconescu's theorem, for any atomic site $({\cal D}, J_{at})$ the canonical functor ${\cal D}\to \Sh({\cal D}, J_{at})$ is characterized by the following universal property: for any Grothendieck topos $\cal E$ and any filtering (in particular, cartesian, if $\cal D$ is cartesian) functor ${\cal D}\to {\cal E}$ which sends arrows in $\cal D$ to epimorphisms can be extended in a unique way (up to isomorphism) to a cartesian cocontinuous functor $\Sh({\cal D}, J_{at}) \to {\cal E}$.

\subsection{Ultrahomogeneous structures as universal models}\label{universal}

Let $\mathbb T$ be an atomic complete theory with a special model $M$. Then, by Theorem \ref{main}, we have a representation
\[
\Sh({\cal C}_{\mathbb T}, J_{\mathbb T})\simeq \Sh({\cal C}^{at}_{\mathbb T}, J_{at})\simeq \Cont(Aut(M))
\]
of its classifying topos.

\begin{theorem}\label{univ}
Let $\mathbb T$ be an atomic complete theory with a special model $M$. Then the model $M$, endowed with the (continuous) canonical action of $Aut(M)$, is a universal model of $\mathbb T$ in the topos $\Cont(Aut(M))$.
\end{theorem}

\begin{proofs}
Let $\Sigma$ be the signature of $\mathbb T$. Consider the $\Sigma$-structure $\tilde{M}$ of $\mathbb T$ in the topos $\Cont(Aut(M))$ given by the canonical (continuous) action of $Aut(M)$ on $M$. Then $\tilde{M}$ is a $\Sigma$-structure is a model of $\mathbb T$ in $\Cont(Aut(M))$; indeed, the canonical point $\Set \to \Cont(Aut(M))$ of $\Cont(Aut(M))$ is a surjection, and the image of $\tilde{M}$ under the inverse image of this point is isomorphic to $M$, which, by our hypothesis, is a model of $\mathbb T$ in $\Set$.  

Let ${\cal C}^{c}_{\mathbb T}$ be the full subcategory of the geometric syntactic category ${\cal C}_{\mathbb T}$ of $\mathbb T$ on the $\mathbb T$-complete formulae.

For any Grothendieck topos $\cal E$, since $\Sh({\cal C}^{c}_{\mathbb T}, J_{at})$ can be identified, via the equivalence $\Sh({\cal C}_{\mathbb T}, J_{\mathbb T})\simeq \Sh({\cal C}^{c}_{\mathbb T}, J_{at})$, with the classifying topos for $\mathbb T$, we have a correspondence between the $\mathbb T$-models in $\cal E$ and the geometric morphisms ${\cal E}\to \Sh({\cal C}^{c}_{\mathbb T}, J_{at})$, which in turn can be identified with the flat $J_{at}$-continuous functors ${\cal C}^{c}_{\mathbb T} \to {\cal E}$; the flat functor corresponding to a model $N$ of $\mathbb T$ in $\cal E$ is given by the functor assigning to any $\mathbb T$-complete formula $\phi(\vec{x})$ its interpretation $[[\vec{x}. \phi]]_{N}$ (and acting on the arrows accordingly). Now, the universal model of $\mathbb T$ in $\Sh({\cal C}^{c}_{\mathbb T}, J_{at})$ corresponds to the flat functor $l:{\cal C}^{c}_{\mathbb T} \to \Sh({\cal C}^{c}_{\mathbb T}, J_{at})$ given by the Yoneda embedding, while the model $M$ corresponds to the functor ${\cal C}^{c}_{\mathbb T}\to \Cont(Aut(M))$ sending any $\mathbb T$-complete formula $\phi(\vec{x})$ to its interpretation $[[\phi(\vec{x})]]_{M}$ in $M$. Hence the two flat functors correspond to each other under the equivalence defined in the proof of Theorem \ref{main}; the $\Sigma$-structure $\tilde{M}$ thus corresponds to the universal model of $\mathbb T$ in $\Sh({\cal C}_{\mathbb T}, J_{\mathbb T})\simeq \Sh({\cal C}^{c}_{\mathbb T}, J_{at})$ under this equivalence, and hence it is itself a universal model of $\mathbb T$ in the topos $\Cont(Aut(M))$, as required. 
\end{proofs}   

The fact that $\tilde{M}$ is a universal model of $\mathbb T$ has several remarkable consequences, notably including the following two ones.

\begin{theorem}\label{specialuniversal}
Let $\mathbb T$ be an atomic complete theory with a special model $M$. Then 

\begin{enumerate}[(i)]
\item For any subset $S\subseteq MA_{1}\times \cdots \times MA_{n}$ which is closed under the action of $Aut(M)$, there exists a (unique up to $\mathbb T$-provable equivalence) geometric formula $\phi(\vec{x})$ over the signature of $\mathbb T$ (where $\vec{x}=(x_{A_{1}}, \ldots, x_{A_{n}})$) such that $S=[[\phi(\vec{x})]]_{M}$; 

\item For any $Aut(M)$-equivariant map $f:S\to T$ between invariant subsets $S$ and $T$ as in $(i)$ there exists a (unique up to $\mathbb T$-provable equivalence) $\mathbb T$-provably functional geometric formula $\theta(\vec{x}, \vec{y})$ from $\phi(\vec{x})$ to $\psi(\vec{y})$, where $S=[[\phi(\vec{x})]]_{M}$ and $T=[[\psi(\vec{y})]]_{M}$, such that the graph of $f$ coincides with $[[\theta(\vec{x}, \vec{y})]]_{M}$.
\end{enumerate}

\end{theorem}

\begin{proofs}
This immediately follows from Theorem \ref{main} and Theorem 2.2 \cite{OCU}.
\end{proofs}

\begin{remark}\label{orbits}
It easily follows from the theorem that for any finite string $A_{1}, \ldots, A_{n}$ of sorts of the signature of the theory $\mathbb T$, the orbits of the action of $Aut(M)$ on $MA_{1}\times \cdots \times MA_{n}$ coincide precisely with the interpretations $[[\phi(\vec{x})]]_{M}$ of $\mathbb T$-complete formulae $\phi(\vec{x})$, where $\vec{x}=(x^{A_{1}}, \ldots, x^{A_{n}})$, that is they correspond exactly to the $\mathbb T$-provable equivalence classes of $\mathbb T$-complete formulae in the context $\vec{x}$.
\end{remark}

It is interesting to investigate to which extent a structure is determined by its automorphism group. We can immediately deduce, from Remark \ref{orbits}, is that if $M$ and $N$ are two special models of an atomic complete theory then for any finite string $A_{1}, \ldots, A_{n}$ of sorts of the signature of the theory $\mathbb T$, the orbits of the action of $Aut(M)$ on $MA_{1}\times \cdots \times MA_{n}$ are in bijective correspondence with the orbits of the action of $Aut(N)$ on $NA_{1}\times \cdots \times NA_{n}$. The following result shows that, if the topological group $Aut(M)$ is discrete then $M$ is uniquely determined by it.

\begin{proposition}
Let $\mathbb T$ be an atomic complete theory with two special models $M$ and $N$. If $Aut(M)$ is discrete then $M\cong N$.
\end{proposition} 

\begin{proofs}
By Theorem \ref{main}, the classifying topos of $\mathbb T$ can be represented both as $\Cont(Aut(M))$ and as $\Cont(Aut(N))$, and the canonical points $p_{M}:\Set \to \Cont(Aut(M))$ and $\Cont(Aut(N))$ are essential. But there is only one essential point, up to isomorphism, of any topos $[G^{\textrm{op}}, \Set]$ of actions of a discrete group $G$, whence $M\cong N$, as required. 
\end{proofs}

Thanks to Remark \ref{orbits}, we can now characterize the continuous homomorphisms $h:Aut(M')\to Aut(M)$ induced by an interpretation ${\cal C}_{\mathbb T}\to {\cal C}_{{\mathbb T}'}$ (cf. Corollary \ref{int}): they are exactly the homomorphisms $h$ such that for any finite string $A_{1}, \ldots, A_{n}$ of sorts over the signature of $\mathbb T$, the action of $Aut(M')$ via $h$ on a $Aut(M)$-invariant subset of $MA_{1}\times \cdots \times MA_{n}$ is isomorphic (in $\Cont(Aut(M'))$) to the action of $Aut(M')$ on a $Aut(M')$-invariant subset of $M'A_{1}\times \cdots \times M'A_{n}$. A variant of this characterization for coherent theories is obtained by replacing the invariant subsets with invariant subsets with a finite number of orbits.

Next, let us describe the universal model of the theory $\mathbb T$ in terms of the different representations of its classifying topos. If $\mathbb T$ is the theory of homogeneous $\mathbb S$-models for a theory of presheaf type $\mathbb S$ (in the sense of \cite{OC2}), then the classifying topos of $\mathbb T$ can be represented as $\Sh({\cal C}^{\textrm{op}}, J_{at})$, where $\cal C$ is a dense subcategory of finitely presentable $\mathbb S$-models, and the image of any universal model of $\mathbb S$ under the associated sheaf functor $a_{J_{at}}$ yields a universal model of $\mathbb T$; in particular, if the forgetful functor $UA_{1}\times \cdots \times UA_{n}:\textrm{f.p.}{\mathbb S}\textrm{-mod}(\Set)\to \Set$ is a $J_{at}$-sheaf then $UA_{1}\times \cdots \times UA_{n}$ is a universal model of $\mathbb T$ (cf. \cite{OCG}). Notice that $U$ is a $J_{at}$-sheaf if and only if for any arrow $f:c\to d$ in the category $\textrm{f.p.}{\mathbb S}\textrm{-mod}(\Set)$ and any element $x\in dA_{1}\times \cdots \times dA_{n}$ such that for any arrows $g,h:d\to e$ such that $g\circ f=h\circ f$, $(hA_{1}\times \cdots \times hA_{n})(x)=(gA_{1}\times \cdots \times gA_{n})(x)$, there exists a unique $y\in c$ such that $(fA_{1}\times \cdots \times fA_{n})(y)=x$. In fact, this condition is often satisfied; indeed, if for any sort $A$ in the signature of $\mathbb T$ the formula $\{x^{A}. \top\}$ presents a $\mathbb S$-model $M_{A}$ and any morphism in $\textrm{f.p.}{\mathbb S}\textrm{-mod}(\Set)$ with domain $M_{A}$ is a strict monomorphism then the underlying objects of the universal model of $\mathbb S$ in $[\textrm{f.p.}{\mathbb S}\textrm{-mod}(\Set), \Set]$ are $J_{at}$-sheaves. An example of this kind of situations will be given in section \ref{examples}.     

If $UA_{1}\times \cdots \times UA_{n}$ is a $J_{at}$-sheaf then the subobjects of $UA_{1}\times \cdots \times UA_{n}$ in $\Sh({\cal C}^{\textrm{op}}, J_{at})$ are precisely the $J_{at}$-closed subobjects of $UA_{1}\times \cdots \times UA_{n}$ in $[\textrm{f.p.}{\mathbb S}\textrm{-mod}(\Set), \Set]$, that is the subfunctors $H$ of $UA_{1}\times \cdots \times UA_{n}$ with the property that for any arrow $f:c \to d$ in $\cal C$ and $x\in cA_{1}\times \cdots \times cA_{n}$, $(fA_{1}\times \cdots \times fA_{n})(x)\in H(d)$ if and only if $x\in H(c)$. From Theorem 2.2 \cite{OCU} we thus obtain that, under this hypothesis, any such subfunctor is of the form $d\to \overline{[[\phi(\vec{x})]]_{d}}$ for some geometric formula $\phi(x^{A_{1}}, \ldots, x^{A_{n}})$ over the signature of $\mathbb T$, where $\overline{[[\phi(\vec{x})]]_{d}}$ is the set of elements $z$ of $dA_{1}\times \cdots \times dA_{n}$ such that there exists an arrow $c\to d$ and an element $x\in [[\phi(\vec{x})]]_{c}$ such that $z=(fA_{1}\times \cdots \times fA_{n})(x)$. By Theorem \ref{specialuniversal}, we thus have a bijective correspondence between the $J_{at}$-closed subfunctors of $UA_{1}\times \cdots \times UA_{n}$ and the $Aut(M)$-invariant subsets of $MA_{1}\times \cdots \times MA_{n}$; one half of this correspondence can be described as the map sending a subfunctor $H$ to the subset of $MA_{1}\times \cdots \times MA_{n}$ consisting of the elements $z$ of the form $(gA_{1}\times \cdots \times gA_{n})(x)$, where $g:N\to M$ is a $\mathbb S$-model homomorphism and $x\in NA_{1}\times \cdots \times NA_{n}$.    

A similar result for the arrows between subobjects of the universal model holds. 

An important problem in Model Theory is to understand to which extent models $M$ are uniquely determined by their topological automorphism $Aut(M)$. For instance, a natural question is: given a geometric theory $\mathbb T$ and two $\mathbb T$-models $M$ and $N$ in $\Set$, is it true that $Aut(M)\cong Aut(N)$ as topological groups implies that $M\cong N$ as $\mathbb T$-models? The answer to this questions is in general negative (many countexamples are known), but there are some classes of theories and models of them for which such a statement is true. Below, we shall see that this problem admits, in the case $\mathbb T$ is an atomic complete theory and $M$ and $N$ are special models for $\mathbb T$, a natural geometric interpretation in terms of automorphisms of toposes. 

First, we notice that if $h:Aut(M)\to Aut(N)$ is a continuous group homomorphism then the set $N$, with the induced $Aut(M)$-action, gives rise to a model $\overline{N}$ of $\mathbb T$ in the topos $\Cont(Aut(M))$, which is the image of the universal model $\tilde{N}$ of $\mathbb T$ in $\Cont(Aut(N))$ under the inverse image of the geometric morphism $\Cont(h):\Cont(Aut(M))\to \Cont(Aut(N))$ induced by $h$. By the universal property of classifying toposes, we thus have a geometric morphism $g:\Cont(Aut(M)) \to \Cont(Aut(M))$ such that $g^{\ast}(\tilde{M})\cong \overline{N}$, where $\tilde{M}$ is the universal model of $\mathbb T$ in $\Cont(Aut(M))$. Notice that if $h$ is a homeomorphism then the induced morphism $\Cont(h)$ is an equivalence and hence $\overline{N}$ is a universal model of $\mathbb T$ in $\Cont(Aut(M))$ and $g$ is an equivalence. The condition that $N$ should be isomorphic to $M$ as $\mathbb T$-models can thus be reformulated as the requirement that the morphism $g$ should be isomorphic to the identity on $\Cont(Aut(M))$. In general, the relationship between the geometric endomorphisms on the classifying topos of the theory and the models of the theory is expressed by the following result.

\begin{theorem}\label{endo}
Let $\mathbb T$ be an atomic complete theory with a special model $M$. Then there is a bijective correspondence between the (isomorphism classes of) geometric endomorphisms $f:\Set[{\mathbb T}]\to \Set[{\mathbb T}]$ on the classifying topos of $\mathbb T$ and the pairs $(N, h)$, where $h:Aut(M)\to Aut(N)$ is a continuous group homomorphism. Under this bijection, the automorphisms on $\Set[{\mathbb T}]$ correspond to the pairs $(N, h)$ where $N$ is a special model of $\mathbb T$ and $h$ is a continuous group homomorphism such that the geometric morphism $\Cont(h):\Cont(Aut(M))\to \Cont(Aut(N))$ is an equivalence. 
\end{theorem}
 
\begin{proofs}
Let us define the correspondence as follows. Given a pair $(N, h)$, where $N$ is a model of $\mathbb T$ in $\Set$ and $h$ is a continuous homomorphism $h:Aut(M)\to Aut(N)$, we have a geometric morphism
\[
\Cont(h):\Cont(Aut(M))\to \Cont(Aut(N)),
\]
and hence a $\mathbb T$-model $\Cont(h)^{\ast}(\tilde{N})$ in $\Cont(Aut(M))$. Notice that the image $\Cont(h)^{\ast}(\tilde{N})$ of $\tilde{N}$ under $\Cont(h)^{\ast}$ can be identified with the action $\alpha:Aut(M)\times N \to N$ of $Aut(M)$ on $N$ induced by the canonical action of $Aut(N)$ on $N$ via the homomorphism $h$. By the universal property of $\Cont(Aut(M))$ as a classifying topos of $\mathbb T$, we thus have a unique geometric morphism $f:\Cont(Aut(M))\to \Cont(Aut(M))$ (up to isomorphism) such that $f^{\ast}(\tilde{M})\cong \Cont(h)^{\ast}(\tilde{N})$. We associate $f$ to $(N, h)$ in our correspondence. Conversely, given a geometric endomorphism $f:\Cont(Aut(M))\to \Cont(Aut(M))$, the image $f^{\ast}(\tilde{M})$ is a $\mathbb T$-model in the topos $\Cont(Aut(M))$; notice that a $\mathbb T$-model in the topos $\Cont(Aut(M))$ can be regarded as a $\mathbb T$-model $N$ in $\Set$ endowed with a continuous action $\alpha:Aut(M)\times N \to N$ such that for any relation symbol $R$ over the signature of $\mathbb T$, $R_{N}$ is $Aut(M)$-invariant and for any function symbol $k$ of arity $n$ over the signature of $\mathbb T$ the following diagram commutes:
\[  
\xymatrix {
Aut(M)\times N^{n} \ar[d]^{\alpha_{n}} \ar[r]^{1\times k_{N}} & Aut(M)\times N \ar[d]^{\alpha}\\
N^{n} \ar[r]_{k_{N}} & N, } 
\] 
where $\alpha_{n}$ is the map $Aut(M)\times N^{n} \to N^{n}$ sending any string $(r, x_{1}, \ldots, x_{n})$ to $(\alpha(r, x_{n}), \ldots, \alpha(r, x_{n}))$.

From this explicit description we see that the function defined on $Aut(M)$ obtained by `currying' $\alpha:Aut(M)\times N \to N$ takes values in the set of $\mathbb T$-model homomorphisms $N \to N$; moreover, since $\alpha$ is an action of a group, all such homomorphisms are bijective, that is we have a function $h:Aut(M)\to Aut(N)$. This function is continuous as the action $\alpha:M\times N \to N$ is continuous, and the action $\alpha$ coincide with the action of $Aut(M)$ on $N$ induced by the canonical action of $Aut(N)$ on $N$ via the homomorphism $h$. We assign the pair $(N, f)$ to $h$.

It is clear that the assignments that we have just defined are inverse to each other and hence define a bijection between the (isomorphism classes of) geometric endomorphisms $f:\Set[{\mathbb T}]\to \Set[{\mathbb T}]$ on the classifying topos of $\mathbb T$ and the pairs $(N, h)$, where $h:Aut(M)\to Aut(N)$ is a continuous group homomorphism. This proves the first part of the theorem; it remains to show that such correspondence restrict to a bijection between the (isomorphism classes of) equivalences $f$ on $\Set[{\mathbb T}]$ and the pairs $(N, h)$ such that $N$ is special and $\Cont(h)$ is an equivalence.  

If $N$ is special then $\tilde{N}$ is a universal model of $\mathbb T$ in the topos $\Cont(Aut(N))$ (by Theorem \ref{univ}) and hence if $\Cont(h)$ is an equivalence $\Cont(h)^{\ast}(\tilde{N})$ is a universal model of $\mathbb T$ in $\Cont(Aut(M))$; therefore the corresponding geometric morphism $f$ is an equivalence (by the universal property of classifying toposes). Conversely, suppose that $f:\Cont(Aut(M))\to \Cont(Aut(M))$ is an equivalence. Then it sends atoms to atoms and hence for any $\mathbb T$-complete formula $\phi(\vec{x})$, the set $[[\phi(\vec{x})]]_{N}$, equipped with the action $\alpha:Aut(M)\times [[\phi(\vec{x})]]_{N} \to [[\phi(\vec{x})]]_{N}$, is an atom of the topos $\Cont(Aut(M))$, i.e. it is a non-empty transitive action; but this action can be regarded as the composite of the canonical action of $Aut(N)$ on $N$ with the homomorphism $h:Aut(M)\to Aut(N)$, from which it follows that the action of $Aut(N)$ on $[[\phi(\vec{x})]]_{N}$ is transitive (and non-empty); in other words, $N$ is special. The fact that $\Cont(h):\Cont(Aut(M))\to \Cont(Aut(N))$ is an equivalence follows from the fact that, $N$ being special, $\tilde{N}$ is a universal model of $\mathbb T$ in $\Cont(Aut(N))$ and $\Cont(h)^{\ast}(\tilde{N})$ is a universal model of $\mathbb T$ in $\Cont(Aut(M))$, it being the image of the universal model $\tilde{M}$ in $\Cont(Aut(M))$ under the equivalence $f$.       
\end{proofs} 

Theorem \ref{endo} provides us with a useful geometric perspective on the relationships between the automorphism groups of models of atomic complete theories, which makes it possible to investigate such relationships by analyzing the endomorphisms of the classifying topos of the theory, a task of entirely categorical/geometric nature. As an example, consider a topos of the form $\Cont(G)$, where $G$ is a \emph{discrete} group. Then there is exactly one isomorphism class of geometric equivalences on $\Cont(G)$. Indeed, $\Cont(G)$ coincides with the presheaf topos $[G^{\textrm{op}}, \Set]$ and any geometric equivalence on $[G^{\textrm{op}}, \Set]$ must be essential (as it is an equivalence). The left adjoint to the inverse image must therefore preserve indecomposable projective objects (as it is an equivalence) and hence must be induced by a functor $G \to G$; but all such functors are isomorphic to the identity, as $G$ is a group. Notice that the fact that there are no non-trivial autormorphisms on such toposes corresponds to the fact that in the context of discrete Galois theory there is exactly one (up to isomorphism) special model of the theory.

\subsection{Coherence}\label{coherence}

Recall from \cite{blasce} that a topological group $G$ is said to be \emph{coherent} if for any open subgroup $H$ of $G$ the number of subsets of the form $HgH$ for $g\in G$ is finite; it is proved in \cite{blasce} that a topological group $G$ is coherent if and only if the topos $\Cont(G)$ is coherent. A prodiscrete topological group is coherent if and only if it is profinite (cf. section D3.4 \cite{El}).   

Let $\mathbb T$ be an atomic complete geometric theory with a special model $M$. Then $\mathbb T$ is syntactically equivalent to a coherent theory (over its signature) if and only if for any finite string of sorts $A_{1}, \ldots, A_{n}$ over the signature of $\mathbb T$, the action of $Aut(M)$ on $MA_{1}\times \cdots \times MA_{n}$ has only a finite number of orbits (see section \ref{univ} and the criterion for a geometric theory to be coherent established in \cite{OCS}).

From Theorem \ref{main} we thus deduce that for any atomic and complete coherent theory $\mathbb T$ with a special model $M$, the topological group $Aut(M)$ is coherent; in particular, for any string of elements $a_{1}, \ldots, a_{n}$ of $M$ there exists a finite number of automorphisms $f_{1}, \ldots, f_{m}$ of $M$ such that every automorphism $f$ of $M$ can be written as $gf_{j}h$ for some $j\in \{1, \ldots m\}$, where $g$ and $h$ are automorphisms which fix all the $a_{i}$. Anyway, this property holds more generally for any atomic complete theory which is \emph{Morita-equivalent} to a coherent (atomic and complete) theory with a special model $M$; for example, $\mathbb T$ might be the theory of homogeneous $\mathbb S$-models where $\mathbb S$ is a theory of presheaf type such that its category of finitely presentable models satisfies AP and JEP and has all fc finite colimits (notice that the latter condition is automatically satisfied if $\mathbb S$ is coherent, cf. \cite{flatcoh}) and such that there exists a ultrahomogeneous $\mathbb S$-model.  

For any topological group $G$, the coherent objects of the topos $\Cont(G)$ are exactly the compact objects, that is the actions with a finite number of orbits (cf. section D3.4 \cite{El}). Then, recalling that for any coherent theory its syntactic pretopos can be characterized, up to equivalence, as the full subcategory of its classifying topos on its coherent objects, we deduce the following characterization of the continuous homomorphisms $Aut(M')\to Aut(M)$ between the automorphism groups $Aut(M)$ and $Aut(M')$ of special models $M$ and $M'$ respectively of atomic complete theories $\mathbb T$ and ${\mathbb T}'$ which are induced by an interpretation of $\mathbb T$ into ${\mathbb T}'$ (in the sense of section \ref{functorialization}) given by a coherent functor ${\cal P}_{\mathbb T}\to {\cal P}_{{\mathbb T}'}$ between their syntactic pretoposes: these are precisely the continuous homomorphisms such that for any continuous $Aut(M)$-action on a set $X$ with a finite number of orbits, the induced action of $Aut(M')$ on $X$ also has a finite number of orbits.

\section{Examples}\label{examples}

In this section we present some examples applications of the theory developed above in different mathematical fields.

\subsection{Preliminaries on theories of presheaf type}\label{prelim}

Let $\mathbb T$ be a theory of presheaf type. It is natural to wonder if there is a theory of presheaf type ${\mathbb T}'$ whose models in $\Set$ can be identified with those of $\mathbb T$ in $\Set$ and whose model homomorphisms between them can be identified with the injective $\mathbb T$-model homomorphism. This problem is relevant for our purposes because the theories of homogeneous $\mathbb S$-models are interesting especially when the arrows between the finitely presentable $\mathbb S$-models are monic.

The most natural candidate for a theory ${\mathbb T}'$ satisfying the above-mentioned property is the theory obtained from $\mathbb T$ by adding to its signature a binary predicate $D$ and the coherent sequents expressing the fact that $D$ is a provable complement to the equality relation; we shall call this theory the \emph{injectivization} of $\mathbb T$.

In general, there are various criteria that one can use for investigating whether a theory is of presheaf type.
   
\begin{theorem}\label{crit}
Let $\mathbb T$ be a theory of presheaf type and ${\mathbb T}'$ be a quotient of $\mathbb T$. Suppose that there exists a set $\cal A$ of finitely presentable ${\mathbb T}'$-models which are finitely presentable as $\mathbb T$-models and which are jointly conservative for ${\mathbb T}'$. Then ${\mathbb T}'$ is of presheaf type, and every finitely presentable ${\mathbb T}'$-model is a retract of a model in $\cal A$. 
\end{theorem}

\begin{proofs}
Since every model in $\cal A$ is finitely presentable as a $\mathbb T$-model, we have a geometric inclusion $i:[{\cal A}, \Set] \hookrightarrow \Set[{\mathbb T}]\simeq [\textrm{f.p.} {{\mathbb T}}\textrm{-mod}(\Set), \Set]$ induced by the canonical inclusion ${\cal A} \hookrightarrow \textrm{f.p.} {{\mathbb T}}\textrm{-mod}(\Set)$. This subtopos corresponds, by the duality theorem of \cite{OCthesis}, to a quotient ${\mathbb T}'$ of $\mathbb T$ classified by the topos $[{\cal A}, \Set]$, which can be characterized as the collection of all geometric sequents which hold in every model in $\cal A$, that is, as ${\mathbb T}'$ itself (as these models are jointly conservative for ${\mathbb T}'$ by our hypothesis). Therefore ${\mathbb T}'$ is of presheaf type classified by the topos $[{\cal A}, \Set]$; but the finitely presentable ${\mathbb T}'$-models are the finitely presentable objects of $\Ind{\cal A}$, that is as the retracts of objects of $\cal A$ in $\Ind{\cal A}\simeq {{\mathbb T}'}\textrm{-mod}({\Set})$. This completes the proof of the theorem. 
\end{proofs}
  
To proceed with our analysis, we need to recall some basic facts about filtered colimits in categories of structures. 

Let $\Sigma$ be a first-order signature. It is well-known that the category $\Sigma\textrm{-str}(\Set)$ of $\Sigma$-structures in $\Set$ has all filtered colimits (cf. for example \cite{El}), which can be constructed as follows. For any diagram $D:{\cal I}\to \Sigma\textrm{-str}(\Set)$ defined on a filtered category ${\cal I}$, its colimit $colim(D)$ is defined by taking, for any sort $A$ of $\Sigma$, $colim(D)A$ equal to the colimit of the functor $DA:{\cal I}\to \Set$, that is (by Proposition 2.13.3 \cite{borceux}) the disjoint union $\coprod D(i)A$ for $i\in {\cal I}$ modulo the equivalence relation $\approx_{A}$ defined as follows: given $x\in D(i)A$ and $x'\in D(i')A$, $x\approx_{A} x'$ if and only if there exist arrows $f:i\to i''$ and $g:i'\to i''$ in $\cal I$ such that $(DA f)(x)=(DA g)(x')$. Denoted by $J_{i}:D(i)\to colim(D)$ the canonical colimit maps (for $i\in {\cal I}$), for any function symbol $f:A_{1}, \ldots, A_{n}\to B$ in $\Sigma$, we define $colim(D)f:colim(D)A_{1}\times \cdots \times colim(D)A_{n}\to colim(D)B$ as follows. For any $(x_{1}, \ldots, x_{n})\in colim(D)A_{1}\times \cdots \times colim(D)A_{n}$ we can suppose, thanks to induction and the fact that $\cal I$ satisfies the joint embedding property, that there exists $i\in I$ such that for any $k\in \{1, \ldots, n\}$, $x_{k}\in DA_{k}(i)$; we thus set $colim(D)f((x_{1}, \ldots, x_{n}))$ equal to $J_{i}B(f_{i}((x_{1}, \ldots, x_{n})))$, where $f_{i}$ is the interpretation of the function symbol $f$ in the structure $D(i)$. It is easy to verify, again by induction using the fact that $\cal I$ satisfies the joint embedding property, that this map is well-defined. For any relation symbol $R\mono A_{1}, \ldots, A_{n}$ in $\Sigma$, we define $colim(D)R$ as the subset of $colim(D)A_{1}\times \cdots \times colim(D)A_{n}$ consisting of the tuples of the form $(x_{1}, \ldots, x_{n})\in D(i)R$ for some $i\in {\cal I}$ (where $D(i)R$ is the interpretation of the relation symbol $R$ in the $\Sigma$-structure $D(i)$). Notice that, since $\cal I$ satisfies the weak coequalizer property, for any $i\in {\cal I}$ and any $x, x'\in D(i)A$, $x\approx_{A} x'$ if and only if there exists an arrow $s:i\to j$ in $\cal I$ such that $D(s)A(x)=D(s)A(x')$; note also that if $x\approx_{A} x'$ then for any arrow $t:i\to k$ in $\cal K$, $D(t)A(x)\approx_{A} D(t)A(x')$.    

We recall from \cite{El} (Lemma D2.4.9) that for any geometric theory $\mathbb T$ over a signature $\Sigma$, ${\mathbb T}\textrm{-mod}(\Set)$ of $\mathbb T$-models in $\Set$ is closed in ${\Sigma}\textrm{-str}(\Set)$ under filtered colimits.

\begin{lemma}\label{fp}
Let $\cal D$ be a category with filtered colimits and $M$ be an object of $\cal D$. Then $M$ is finitely presentable in $\cal D$ if and only if for any filtered diagram $D:{\cal I}\to {\cal D}$, any arrow $M\to colim(D)$ factors through one of the canonical colimit arrows $J_{i}:D(i)\to colim(D)$ and for any arrows $f:M\to D(i)$ and $g:M\to D(j)$ in $\cal D$ such that $J_{i}\circ f=J_{j}\circ g$ then there exists $k\in {\cal I}$ and two arrows $s:i\to k$ and $t:j\to k$ such that $D(s)\circ f=D(t)\circ g$. 
\end{lemma}

\begin{proofs}
By definition, $M$ is finitely presentable if and only if the hom functor $Hom_{\cal D}(M, -):{\cal D}\to \Set$ preserves filtered colimits, that is for any diagram $D:{\cal I}\to {\cal D}$ defined on a filtered category $\cal I$, the canonical map $\xi$ from the colimit in $\Set$ of the composite diagram $Hom_{\cal D}(M, -)\circ D$ to $Hom_{\cal D}(M, colim(D))$ is a bijection. The former can be described, by Proposition 2.13.3 \cite{borceux}, as the disjoint union of the $Hom_{\cal D}(M, D(i))$ (for $i\in {\cal I}$) modulo the equivalence relation $\approx$ defined by: for any $f:M\to D(i)$ and $g:M\to D(j)$, $f\approx g$ if and only if there exists $k\in {\cal I}$ and two arrows $s:i\to k$ and $t:j\to k$ such that $D(s)\circ f=D(t)\circ g$. The map $\xi$ sends an equivalence class $[f]$ of an arrow $f:M\to D(i)$ to the composite arrow $J_{i}\circ f$, where $J_{i}$ is the canonical colimit arrow $D(i)\to colim(D)$; hence the first and second condition in the statement of the lemma are immediately seen to correspond respectively to the requirements that $\xi$ be surjective and injective.
\end{proofs}

Another useful result which follows straightforwardly from the explicit construction of filtered colimits of model categories of geometric theories is the following.

\begin{proposition}\label{injective}
Let $\mathbb T$ be a geometric theory. Then the category of $\mathbb T$-models and injective homomorphisms between them is closed under filtered colimits in ${\mathbb T}\textrm{-mod}(\Set)$. Explicitly, if $M$ is a model of $\mathbb T$ which is a filtered colimit of a diagram $D:{\cal I}\to {\mathbb T}\textrm{-mod}(\Set)$ and for any sort $A$ over the signature of $\mathbb T$ and any arrow $f:i\to j$ in $\cal I$, $D(f)A:D(i)A\to D(j)A$ is injective then the universal colimit arrows $D(i)A\to MA$ (for $i\in {\cal I}$) are also injective; and for any cocone on the diagram $D$ an object $N$ of ${\mathbb T}\textrm{-mod}(\Set)$ whose legs are injective, the universal arrow $M\to N$ is also injective. 
\end{proposition}   

\begin{proofs}
The theory ${\mathbb T}'$ obtained from $\mathbb T$ by adding a binary predicate and the coherent sequents asserting that it is complemented to the equality relation is clearly geometric and hence its category of models, which coincides with the category of $\mathbb T$-models in $\Set$ and injective homomorphisms between them, has filtered colimits (Lemma D2.4.9 \cite{El}), which are computed as in $\Set$; therefore the category of ${\mathbb T}'$-models in $\Set$, regarded as a subcategory of the category of $\mathbb T$-models in $\Set$, is closed under filtered colimits, as required.  
\end{proofs}

The next result shows that (set-theoretic) finiteness implies finite presentability in categories of set-based models of geometric theories. Recall that a model $M$ of a geometric theory $\mathbb T$ in the category $\Set$ is said to be finitely presentable if the hom functor $Hom_{{\mathbb T}\textrm{-mod}(\Set)}(M, -):{\mathbb T}\textrm{-mod}(\Set)\to \Set$ preserves filtered colimits. We shall say that a model of a geometric theory over a signature $\Sigma$ is \emph{finite} if the disjoint union of the $MA$ for $A$ sort over $\Sigma$ is a finite set.

\begin{theorem}\label{finite}
Let $\mathbb T$ be a geometric theory over a signature consisting of a finite number of sorts, function symbols and relation symbols. Then any finite model $M$ of $\mathbb T$ is finitely presentable.
\end{theorem}

\begin{proofs}
Let us use the criterion for finite presentability provided by Lemma \ref{fp} and the description of filtered colimits in categories of set-based models of geometric theories given above. Let $D$ be a diagram defined on a filtered category $\cal I$ with values in the category ${\mathbb T}\textrm{-mod}(\Set)$ and let $M$ be a finite model of $\mathbb T$. Let $\chi:M\to colim(D)$ be a $\mathbb T$-model homomorphism; first, let us prove that there exists $i\in {\cal I}$ and a $\mathbb T$-model homomorphism $\chi':M\to D(i)$ such that $J_{i}\circ \chi'=\chi$.

Let us suppose that we have for any sort $A$ of $\Sigma$ a \emph{function} $\xi_{A}:MA\to D(i)A$ (for some fixed $i\in {\cal I}$) such that $\chi A=J_{i}\circ \xi_{A}$ for all $A$, and show that, through an inductive process, we can modify the functions $\xi_{A}$ in such a way so to make them the components of a homomorphism of $\mathbb T$-models $\xi:M\to D(i)$ such that $J_{i}\circ \xi=\chi$. Recall that, for any $\mathbb T$-models $M$ and $N$ in $\Set$, a bunch of functions $\xi_{A}:MA\to NA$ is the set of components of a $\mathbb T$-model homomorphism $M\to N$ if and only if the following conditions are satisfied: $(1)$ for any function symbol $f:A_{1}, \ldots, A_{n}\to B$ over $\Sigma$ and any $(x_{1}, \ldots, x_{n})\in MA_{1}\times \cdots \times MA_{n}$, if $y=fM(x_{1}, \ldots, x_{n})$ then $\xi_{B}(y)=fN(\xi_{A_{1}}(x_{1}), \ldots, \xi_{A_{n}}(x_{n}))$ and $(2)$ for any relation symbol $R\mono A_{1}, \ldots, A_{n}$ over $\Sigma$ and any $(x_{1}, \ldots, x_{n})\in R_{M}$, $(\xi_{A_{1}}(x_{1}), \ldots, \xi_{A_{n}}(x_{n}))\in R_{N}$. Notice that, in our case, the condition that the $\xi_{A}$ should be the components of a $\mathbb T$-model homomorphism $M\to D(i)$ can be formulated as a \emph{finite} set of conditions involving elements of the $MA$ and the $D(i)A$. Let us now show that, starting from an arbitrary bunch of functions $\xi_{A}:MA\to D(i)A$, any such condition can be achieved at the cost of replacing each $\xi_{A}:MA\to D(i)A$ with $D(s)A\circ \xi_{A}:MA\to D(j)A$ for some arrow $s:i\to j$ in $\cal I$. Take a condition $C$ of type $(1)$ associated to a tuple $(x_{1}, \ldots, x_{n})\in MA_{1}\times \cdots \times MA_{n}$; since $\chi:M\to colim(D)$ is a $\mathbb T$-model homomorphism, we have that $f_{D(i)}(\xi_{A_{1}}(x_{1}), \ldots, \xi_{A_{n}}(x_{n})) \approx_{B} \xi_{B}(f_{M}(x_{1}, \ldots, x_{n}))$; indeed, by considering the images of both elements under $J_{i}B$ we obtain that
\[ 
J_{i}B(f_{D(i)}(\xi_{A_{1}(x_{1}), \ldots, \xi_{A_{n}}(x_{n})}))=f_{colimD}(J_{i}A_{1}\times \cdots \times J_{i}A_{n}((\xi_{A_{1}}(x_{1}), \ldots, \xi_{A_{n}}(x_{n}))))
\]
\[
=f_{colimD}((\chi A_{1}\times \cdots \times \chi A_{n})(x_{1}, \ldots, x_{n}))=\chi B(f_{M}(x_{1}, \ldots, x_{n}))
\]
\[
=J_{i}B(\xi_{B}(f_{M}(x_{1}, \ldots, x_{n}))),
\]
where the central equalities follow from the fact that $\chi$ and $J_{i}$ are $\mathbb T$-model homomorphisms and from the definition of $f_{colim(D)}$. Therefore there exists an arrow $s:i\to j$ in $\cal I$ such that $D(s) B(f_{D(i)}(\xi_{A_{1}}(x_{1}), \ldots, \xi_{A_{n}}(x_{n})))=D(s)B (\xi_{B}(f_{M}(x_{1}, \ldots, x_{n})))$; but, $D(s)$ being an homomorphism of $\mathbb T$-models $D(i)\to D(j)$, we have $D(s) B(f_{D(i)}(\xi_{A_{1}}(x_{1}), \ldots, \xi_{A_{n}}(x_{n})))=f_{D(j)}((D(s)A_{1}\circ \xi_{A_{1}})(x_{1}), \ldots, (D(s)A_{n}\circ \xi_{A_{n}})(x_{n}))$, whence the bunch of functions $D(s)A\circ \xi_{A}:MA\to D(j)A$ satisfies our condition $C$, as required. Now, let us consider a condition $B$ of type $(2)$. Let $R$ be any relation symbol $R\mono A_{1}, \ldots, A_{n}$ over $\Sigma$ and let $(x_{1}, \ldots, x_{n})$ be in $R_{M}$. Since $\chi:M\to colim(D)$ is a $\mathbb T$-model homomorphism then $(\chi A_{1}(x_{1}), \ldots, \chi A_{n}(x_{n}))\in R_{colim(D)}$; this means that there exists $j\in {\cal I}$ and a tuple $(y_{1}, \ldots, y_{n})\in D(j)A_{1}\times \cdots \times D(j)A_{n}$ such that $(\chi A_{1}(x_{1}), \ldots, \chi A_{n}(x_{n}))=(J(j)A_{1}(y_{1}), \ldots, J(j)A_{n}(y_{n}))$. Clearly, since $\cal I$ satisfies the joint embedding property and for any arrow $s$ in $\cal I$, $D(s):D(i)\to D(j)$ is a $\mathbb T$-model homomorphism, we can suppose that there is an arrow $s:i\to j$ in $\cal I$. So have that $\xi_{A_{1}}(x_{1})\approx_{A_{1}} y_{1}, \ldots, \xi_{A_{n}}(x_{1})\approx_{A_{n}} y_{n}$; but since $s$ is an arrow in $\cal I$ we have that $D(s)A_{1}(\xi_{A_{1}}(x_{1}))\approx_{A_{1}} y_{1}, \ldots, D(s)A_{n}(\xi_{A_{n}}(x_{n}))\approx_{A_{n}} y_{n}$. Therefore there exists $s_{1}:j\to k_{1}$ such that $D(s_{1})A_{1}(D(s)A_{1}(\xi_{A_{1}}(x_{1})))=D(s_{1})A_{1}(y_{1})$. Now, the fact that $D(s)A_{2}(\xi_{A_{2}}(x_{2}))\approx_{A_{2}} y_{2}$ implies that we have $D(s_{1})A_{2}(D(s)A_{2}(\xi_{A_{2}}(x_{2}))) \approx_{A_{2}} D(s_{1})A_{2}(y_{2})$; if we go on inductively for $n$ steps, at the end we will obtain an arrow $t$, namely the composition $s_{n}\circ s_{n-1}\circ \cdots \circ s_{1}\circ s$, such that 
\[
D(t)A_{1}(\xi_{A_{1}}(x_{1}))=D(t')A_{1}(y_{1}), \ldots, D(t)A_{n}(\xi_{A_{n}}(x_{n}))=D(t')A_{n}(y_{n}),
\]
where $t'=s_{n}\circ s_{n-1}\circ \cdots \circ s_{1}$. Since, $D_{t'}$ being a $\mathbb T$-model homomorphism, we have that $(D(t')A_{1}(y_{1}), \ldots, D(t')A_{n}(y_{n}))\in R_{D(cod(t))}$, the bunch of maps of the form $D(t)A \circ \xi_{A}$ satisfies condition $B$.    

Next, we note for any of the conditions of type $(1)$ or $(2)$ (which are indexed by strings of elements of the $MA$), if a bunch of functions $\xi_{A}:MA\to D(i)A$ satisfies it then for any arrow $s:i\to j$ in $\cal I$, the bunch of functions $D(s)A\circ \xi_{A}:MA\to D(j)A$ also satisfies it; indeed, since $D$ is a diagram defined on $\cal I$ with values in the category ${\mathbb T}\textrm{-mod}(\Set)$, $D(s)$ is a $\mathbb T$-model homomorphism $D(i)\to D(j)$. Also, if the bunch $\xi_{A}:MA\to D(i)A$ represents a factorization of the function $\chi:M\to colim(D)$ across a canonical map $J(i):D(i)\to colim(D)$ then for any $s:i\to j$ in $\cal I$ the bunch $D(s)A \circ \xi_{A}$ represents a factorization of $\chi$ across the canonical map $J_{j}:D(j)\to colim(D)$. This ensures that we can, through an inductive process, modify our initial bunch of functions $\xi_{A}$ by replacing it with a bunch of arrows of the form $D(s)A \circ \xi_{A}$ (for some arrow $s$ in $\cal I$) until we arrive, at the last step, to a bunch of the form $D(s)A\circ \xi_{A}$ which satisfies all the conditions of either type $(1)$ or type $(2)$ and which therefore yields a factorization of the homomorphism $\chi:M\to colim(D)$ through the canonical map $J(cod(s)): D(cod(s))\to colim(D)$ which is a $\mathbb T$-model homomorphism $\xi:M\to D(cod(s))$. 

To complete the proof that any $\mathbb T$-model homomorphism $\chi: M\to colim(D)$ factors in ${\mathbb T}\textrm{-mod}(\Set)$ through a canonical homomorphism $J_{i}:D(i)\to colim(D)$, it just suffices to prove that there exists, for some $i\in {\cal I}$, a bunch of functions $\xi_{A}:MA\to D(i)A$ such that $J_{i}A\circ \xi_{A}=\chi A$. Since there is only a finite number of sorts and $\cal I$ satisfies the joint embedding property, it suffices to prove that for any sort $A$ over $\Sigma$ there exists a function $\xi_{A}:MA\to D(i)A$, for some $i\in {\cal I}$, such that $\chi A=J_{i}A \circ \xi_{A}$. Since $MA$ is finite, we can number its elements $x_{1}, \ldots, x_{n}$. Since $\chi A$ takes values in $colim(D)A$ there exists $k_{1}\in {\cal K}$ and $z_{1}\in D(k_{1})A$ such that $\chi A(x_{1})=J_{k_{1}}A(z_{1})$, and there exists $k_{2}\in {\cal K}$ and $z_{2}\in D(k_{2})A$ such that $\chi A(x_{2})=J_{k_{2}}A(z_{2})$. Since $\cal I$ satisfies the joint embedding property, there exists $a_{1}\in {\cal I}$ and arrows $h_{1}:k_{1}\to a_{1}$ and $h_{1}':k_{2}\to a_{1}$. In the same line of thought, there exists $k_{3}\in {\cal K}$ and $z_{3}\in D(k_{3})A$ such that $\chi A(x_{3})=J_{k_{3}}A(z_{3})$ and hence, again by the fact that $\cal I$ satisfies the joint embedding property, we have an object $a_{2}$ of $\cal I$ and arrows $h_{2}:a_{1}\to a_{2}$ and $h_{2}':k_{3}\to a_{2}$. Let us repeat this process for exactly $n-1$ steps; at the end we will have a sequence of arrows $h_{j}:a_{j-1}\to a_{j}$ for $j=1, \ldots, n-1$ (where $a_{0}=z_{1}$), which allows us to define $\xi_{A}:MA\to D(a_{n-1})A$ as follows: $\xi_{A}(x_{1})=D(h_{n-1}\circ h_{n-2} \circ \cdots \circ h_{1})A(x_{1})$, $\xi_{A}(x_{2})=D(h_{n-1}\circ h_{n-2} \circ \cdots \circ h_{2} \circ h_{1}')A(x_{2})$, $\xi_{A}(x_{3})=D(h_{n-1}\circ h_{n-2} \circ \cdots \circ h_{3} \circ h_{2}')A(x_{3}), \ldots, \xi_{A}(x_{n})=D(h_{n-1}')A(x_{n})$.             

To conclude the proof of the theorem, it remains to verify that for any $i, j\in {\cal I}$ and $\mathbb T$-model homomorphisms $f:M\to D(i)$ and $g:N\to D(j)$ such that $J_{i}\circ f=J_{j}\circ g$, there exists $k\in {\cal I}$ and arrows $s:i\to k$ and $t:j\to k$ such that $D(s)\circ f=D(t)\circ g$. First, we notice that it suffices to prove that for any sort $A$ over $\Sigma$ there exists $k\in {\cal I}$ and arrows $s:i\to k$ and $t:j\to k$ such that $D(s)A\circ fA=D(t)A\circ gA$. Indeed, there is only a finite number of sorts and for any two sorts $A$ and $B$, arrows $s:i\to k$ and $t:j\to k$ such that $D(s)A\circ fA=D(t)A\circ gA$ and arrows $s':i\to k'$ and $t':j\to k'$ such that $D(s')B\circ fB=D(t)B\circ gB$, since $\cal C$ satisfies the joint embedding property and the weak coequalizer property, there exist $k''\in {\cal I}$ and arrows $a:k\to k''$ and $b:k'\to k''$ such that $a\circ s=b\circ s'$ and $a\circ t=b\circ t'$. Let us therefore fix a sort $A$ of $\Sigma$ and consider the functions $fA:MA\to D(i)A$ and $gA:MA\to D(j)A$; since $MA$ is finite, we can number its elements as $x_{1}, \ldots, x_{n}$. Consider $fA(x_{1})$ and $gA(x_{1})$; since $J_{i}\circ f=J_{j}\circ g$, we have that $fA(x_{1}) \approx_{A} gA(x_{1})$, from which it follows that there exist arrows $s:i\to k$ and $t:j\to k$ such that $(D(s)A\circ fA)(x_{1})=(D(t)A\circ gA)(x_{1})$; now, since $fA(x_{2}) \approx_{A} gA(x_{2})$ then $(D(s)A\circ fA)(x_{2})\approx_{A} (D(t)A\circ gA)(x_{2})$ and hence there exists an arrow $r:k\to k'$ such that $D(r)A(D(s)A\circ fA)(x_{2})=D(r)A(D(t)A\circ gA)(x_{2})$; by iterating the argument for other $n-2$ steps we obtain a sequence of composable arrows in $\cal I$ such that their composite $z$ is such that $(D(z\circ s)A\circ fA)(x_{i})=(D(z\circ t)A\circ gA)(x_{i})$ for all $i$, that is such that $D(z\circ s)A\circ fA= D(z\circ t)A\circ gA$, as required.        
\end{proofs}

\begin{remark}
The theorem specializes to Example 1.2 (4) \cite{ar} in the case where there are no function symbols.
\end{remark}

As a corollary, we immediately obtain the following result.

\begin{corollary}
Let $\Sigma$ be a signature with only a finite number of function and relation symbols. Then for any category $\cal A$ of finite $\Sigma$-structures and $\Sigma$-structure homomorphisms between them, there exists exactly one theory over the signature $\Sigma$ classified by the topos $[{\cal A},\Set]$ for which the structures in $\cal A$ are jointly conservative, and any finitely presentable model of it is a retract of a model in $\cal A$; in particular, any geometric theory admitting a set $\cal A$ of jointly conservative finite models is of presheaf type classified by the topos $[{\cal A}, \Set]$ and its finitely presentable models are precisely the finite models of the theory. 
\end{corollary}

\begin{proofs}
This follows immediately from Theorem \ref{crit} and Proposition \ref{finite}.
\end{proofs}

\begin{remark}\label{injectivization}
The corollary can be notably applied to theories $\mathbb T$ with enough set-based models such that every model of $\mathbb T$ can be expressed as a filtered colimit of finite models of $\mathbb T$. Notice that, assuming the axiom of choice, every coherent theory has enough models; in particular, it follows from the corollary that, assuming the axiom of choice, the injectivization of any finitary theory of presheaf type whose finitely presentable models are all finite is again of presheaf type.  
\end{remark}

We mention that another general, fully constructive, method for proving that a theory $\mathbb T$ is of presheaf type is to establish a geometrically constructive correspondence between the models of $\cal T$ and the flat functors on $\textrm{f.p.}{\mathbb T}\textrm{-mod}(\Set)$; in fact, this method is applied in \cite{vickers} yielding a criterion for a theory to be of presheaf type which is applicable to a wide class of theories whose finitely presentable models are finite (notably including the theory of decidably linear orders and the theory of decidable sets). A concrete description of an equivalence established by using this method is provided by Exercise VIII (8) \cite{MM}.      

Let us consider a theory of presheaf type $\mathbb T$ and a quotient ${\mathbb T}'$ of ${\mathbb T}$ corresponding, via the duality theorem of \cite{OC6}, to a Grothendieck topology $J$ on $\textrm{f.p.}{\mathbb T}\textrm{-mod}(\Set)^{\textrm{op}}$. The following result gives a characterization of the models of ${\mathbb T}'$ which are finitely presentable as models of $\mathbb T$ in terms of the topology $J$.

The following result is useful for recognizing the quotients of a theory of presheaf type which are again of presheaf type.

\begin{theorem}
Let ${\mathbb T}'$ be a quotient of a theory of presheaf type ${\mathbb T}$ corresponding, via the duality theorem of \cite{OC6}, to a Grothendieck topology $J$ on $\textrm{f.p.}{\mathbb T}\textrm{-mod}(\Set)^{\textrm{op}}$, and let $M$ be a finitely presentable $\mathbb T$-model. Then $M$ is a model of ${\mathbb T'}$ if and only if it is $J$-irreducible (i.e., every $J$-covering sieve on $M$ is maximal). 
\end{theorem}  

\begin{proofs}
First, we recall from \cite{OC6} that we can axiomatize ${\mathbb T}'$ as the quotient obtained from $\mathbb T$ by adding all the sequents $\sigma$ of the form $\phi \vdash_{\vec{x}} \mathbin{\mathop{\textrm{\huge $\vee$}}\limits_{i\in I}}(\exists \vec{y_{i}})\theta_{i}$, where, for any $i\in I$, $[\theta_{i}]:\{\vec{y_{i}}. \psi_{i}\}\to \{\vec{x}. \phi\}$ is an arrow in ${\cal C}_{\mathbb T}$, $\phi(\vec{x})$ and $\psi(\vec{y_{i}})$ are formulae presenting respectively ${\mathbb T}$-models $M_{\phi}$ and $M_{\psi_{i}}$, and the cosieve $S_{\sigma}$ on $M_{\phi}$ in the category $\textrm{f.p.} {\mathbb T}\textrm{-mod}(\Set)$ defined as follows is $J$-covering. For each $i\in I$, $[[\theta_{i}]]_{M_{\psi_{i}}}$ is the graph of a morphism $[[\vec{y_{i}}.\psi_{i}]]_{M_{\psi_{i}}}\to [[\vec{x}.\phi]]_{M_{\psi_{i}}}$; then the image of the generators of $M_{\psi_{i}}$ via this morphism is an element of $[[\vec{x}.\phi]]_{M_{\psi_{i}}}$ and this in turn determines, by definition of $M_{\phi}$, a unique arrow $s_{i}:M_{\phi} \to M_{\psi_{i}}$ in ${\mathbb T}\textrm{-mod}(\Set)$. We define $S_{\sigma}$ as the sieve in ${\textrm{f.p.} {\mathbb T}\textrm{-mod}(\Set)}^{\textrm{op}}$ on $M_{\phi}$ generated by the arrows $s_{i}$ as $i$ varies in $I$.

Conversely, the Grothendieck topology $J$ associated to the quotient ${\mathbb T}'$ can be defined as the collection of all the sieves of the form $S_{\sigma}$ where $\sigma$ is a sequent provable in ${\mathbb T}'$ of the form $\phi \vdash_{\vec{x}} \mathbin{\mathop{\textrm{\huge $\vee$}}\limits_{i\in I}}(\exists \vec{y_{i}})\theta_{i}$, where, for any $i\in I$, $[\theta_{i}]:\{\vec{y_{i}}. \psi_{i}\}\to \{\vec{x}. \phi\}$ is an arrow in ${\cal C}_{\mathbb T}$ and $\phi(\vec{x})$, $\psi(\vec{y_{i}})$ are formulae presenting respectively ${\mathbb T}$-models $M_{\phi}$ and $M_{\psi_{i}}$.
 
For any formula $\phi(\vec{x})$ which presents a $\mathbb T$-model, we denote by $\vec{\xi_{\phi}}$ its set of generators in the $\mathbb T$-model $M_{\phi}$ presented by it.

We shall establish the following fact, from which our thesis will follow at once: for any sequent $\sigma:=(\phi \vdash_{\vec{x}} \mathbin{\mathop{\textrm{\huge $\vee$}}\limits_{i\in I}}(\exists \vec{y_{i}})\theta_{i})$, $\sigma$ holds in $M_{\psi}$ if and only if every arrow $f:M_{\phi}\to M_{\psi}$ in ${\textrm{f.p.} {\mathbb T}\textrm{-mod}(\Set)}$ belongs to $S_{\sigma}$.  

Clearly, $\sigma$ holds in $M_{\psi}$ if and only if for any $\vec{a}\in [[\phi(\vec{x})]]_{M_{\psi}}$ there exists $i\in I$ and $\vec{b}\in [[\psi_{i}(\vec{y_{i}})]]_{M_{\psi}}$ such that $[[\theta_{i}]]_{M_{\psi}}(\vec{b})=\vec{a}$. By the universal property of the model $M_{\phi}$, the tuples $\vec{a}\in [[\phi(\vec{x})]]_{M_{\psi}}$ can be identified with the arrows $M_{\phi}\to M_{\psi}$ in ${\textrm{f.p.} {\mathbb T}\textrm{-mod}(\Set)}$ via the correspondence which assigns to $\vec{a}\in [[\phi(\vec{x})]]_{M_{\psi}}$ the unique $\mathbb T$-model homomorphism $f_{\vec{a}}:M_{\phi}\to M_{\psi}$ such that $f_{\vec{a}}(\vec{\xi_{\phi}})=\vec{a}$, and similarly for the tuples $\vec{b}\in [[\psi_{i}(\vec{y_{i}})]]_{M_{\psi}}$. We can easily see that the property that $[[\theta_{i}]]_{M_{\psi}}(\vec{b})=\vec{a}$ is exactly equivalent to the requirement that $f_{\vec{b}}\circ s_{i}=f_{\vec{a}}$; indeed, the equality $f_{\vec{b}}\circ s_{i}=f_{\vec{a}}$ holds if and only if $(f_{\vec{b}}\circ s_{i})(\vec{\xi_{\phi}})=f_{\vec{a}}(\vec{\xi_{\phi}})$ and, since $s_{i}(\vec{\xi_{\phi}})=[[\theta_{i}]]_{M_{\psi_{i}}}(\vec{\xi_{\psi_{i}}})$, we have that  
$(f_{\vec{b}}\circ s_{i})(\vec{\xi_{\phi}})=f_{\vec{b}}([[\theta_{i}]]_{M_{\psi_{i}}}(\vec{\xi_{\psi_{i}}}))=[[\theta_{i}]]_{M_{\psi}}(f_{\vec{b}}(\vec{\xi_{\phi}}))=[[\theta_{i}]]_{M_{\psi}}(\vec{b})$, where the central equality follows from the naturality of the interpretation of $\theta_{i}$. 

Now the thesis immediately follows from these remarks by observing that

\begin{enumerate}[(i)]
\item if every $J$-covering sieve on $M_{\psi}$ is maximal then any sequent $\sigma$ which is provable in ${\mathbb T}'$ is valid in $M_{\psi}$; indeed, if $\sigma$ is provable in ${\mathbb T}'$ then the cosieve $S_{\sigma}$ is $J$-covering on $M_{\phi}$ and hence for any arrow $f:M_{\phi}\to M_{\psi}$ in ${\textrm{f.p.} {\mathbb T}\textrm{-mod}(\Set)}$ the pullback cosieve $f^{\ast}(S_{\sigma})$ is $J$-covering on $M_{\psi}$ and hence maximal, i.e. $f$ belongs to $S_{\sigma}$.     

\item if $M_{\psi}$ is a model of ${\mathbb T}'$ then for any cosieve $S:=S_{\sigma}$ on $M_{\psi}$ in the category ${\textrm{f.p.} {\mathbb T}\textrm{-mod}(\Set)}$, $\sigma:=(\psi \vdash_{\vec{y}} \mathbin{\mathop{\textrm{\huge $\vee$}}\limits_{i\in I}}(\exists \vec{y_{i}})\theta_{i})$ is valid in $M_{\psi}$ and hence the identity on $M_{\psi}$ belongs to $S_{\sigma}$, i.e. $S_{\sigma}$ is maximal.
\end{enumerate}
\end{proofs}

The following list of examples exhibits instances of Galois theories in different mathematical contexts arising from the application of our general theory. 

\subsection{Discrete Galois theory} 

There are several interesting examples of discrete Galois theories, notably including the well-known ones given by the classical Galois theory of a finite Galois extension and the theory of universal coverings and the fundamental group. 

\subsubsection{Classical Galois theory}

Let $F\subset L$ be a finite-dimensional Galois extension of a field $F$, and let ${\cal L}_{F}^{L}$ be the category of intermediate extensions and homomorphisms between them. We have an equivalence of toposes $\Sh({{\cal L}_{F}^{L}}^{\textrm{op}}, J_{at})\simeq \Cont(Aut_{F}(L))$, where $Aut_{F}(L)$ is the discrete Galois group of $L$ (the category ${\cal L}_{F}^{L}$ and the object $L$ satisfy the hypotheses of Theorem \ref{discreteGalois}). The fundamental theorem of classical Galois theory is exactly equivalent to the assertion that all the arrows in the category ${\cal L}_{F}^{L}$ are strict monomorphisms and that the category ${\cal L}_{F}^{L}$ is atomically complete.   

\subsubsection{Coverings and the fundamental group}

Let $B$ be a path connected, locally path connected and semi locally simply connected space; then there is a universal covering $\tau:\tilde{B}\to B$. Let $\cal B$ be the category whose objects are the connected coverings $p:E\to B$ of $B$ and whose arrows $(p:E\to B)\to (q:E'\to B)$ are the continuous maps $z:E\to E'$ such that $q\circ z=p$. Let $b_{0}$ be a fixed point of $B$, and $\pi_{1}(B, b_{0})$ the fundamental group of $B$ at $b_{0}$.

The category ${\cal B}$ and the object $\tau:\tilde{B}\to B$ satisfy the hypotheses of Theorem \ref{discreteGalois}, and the automorphism group of the object $\tau$ in ${\cal B}$ is isomorphic to the fundamental group $\pi_{1}(B, b_{0})$; hence we have an equivalence of toposes $\Sh({\cal B}^{\textrm{op}}, J_{at})\simeq \Cont(\pi_{1}(B, b_{0}))$. The fundamental theorem of covering theory provides a bijective correspondence between the subgroups of $\pi_{1}(B, b_{0})$ and the connected coverings of $B$; in other words, it says that this equivalence of toposes yields a standard Galois theory, or, equivalently, that every arrow in $\cal B$ is a strict monomorphism and that $\cal B$ is atomically complete.

\subsubsection{Ultrahomogeneous finite groups}

A finite group $G$ is said to be \emph{ultrahomogeneous} if every isomorphism between subgroups of $G$ can be extended to an automorphism of $G$. Ultrahomogeneous finite groups have been completely classified, and a full list can be found in \cite{felgner}. For any such group $G$, the category ${\cal C}_{G}$ of groups which can be embedded in $G$ and injective homomorphisms between them, together with the object $G$, satisfies the hypotheses of Theorem \ref{discreteGalois}. We thus have an equivalence of toposes $\Sh({{\cal C}_{G}}^{\textrm{op}}, J_{at})\simeq \Cont(Aut(G))$, where $Aut(G)$ is the discrete automorphism group of $G$ in ${\cal C}_{G}$.

\subsection{Decidable objects and infinite sets}. Let $\mathbb S$ be the theory of decidable objects, that is the theory over a signature consisting of a unique sort and a binary relation symbol $D$ with no axioms except for the two coherent sequents asserting that $D$ is a (provable) complement to the equality relation, i.e. $(\top \vdash_{x, x} D(x,y) \vee (x=y))$ and $(D(x,y) \wedge x=y \vdash_{x, y} \bot)$. Notice that $\mathbb S$ is the injectivization of the empty theory over the signature entirely consisting of one sort. Clearly, the models of $\mathbb S$ in $\Set$ are precisely the sets, while the $\mathbb S$-model homomorphisms are precisely the injective functions between them. It is easy to prove that $\mathbb S$ is of presheaf type (for example, by using the criteria established in section \ref{prelim} - cf. Remark \ref{injectivization} - or see p. 908 \cite{El} for a different proof), and that its category of finitely presentable models can be identified with the category $\mathbb I$ of finite sets and injections. This category satisfies both the amalgamation and joint embedding properties, and the set $\mathbb N$ of natural numbers is clearly a $\mathbf{I}$-ultrahomogeneous model (in fact, any \emph{infinite} set is equally $\mathbf{I}$-ultrahomogeneous). Thus Theorem \ref{maincategorical} yields an equivalence $\Sh({\mathbb I}^{\textrm{op}}, J_{at})\simeq \Cont(Aut({\mathbb N}))$. It is easy to see that every morphism in $I$ is a strict monomorphism, from which it follows that the given equivalence restricts to an equivalence of ${\mathbb I}^{\textrm{op}}$ with a full subcategory of $\Cont_{t}(Aut({\mathbb N}))$. This in particular implies the (easy) Galois-type property that for any two subsets $S$ and $T$ of $\mathbb N$, if all the automorphisms of $\mathbb N$ which fix $T$ fix $S$ then $S\subseteq T$.  

We notice that the formula $\top(\vec{x})$ presents a $\mathbb S$-model, namely the singleton set. This implies, in view of the results of section \ref{universal}, that for any natural number $n$ the assignments $M \to R_{M}$ sending to any finite set $M$ a subset $R_{M}\subseteq M^{n}$ in such a way that for any injective function $f:M\to N$ between finite sets $M$ and $N$ and any $x\in M^{n}$, $f^{n}(x)\in R_{N}$ if and only if $x\in R_{M}$ are in bijective correspondence with the $Aut({\mathbb N})$-invariant subsets of ${\mathbb N}^{n}$ (via the map sending to any such assignment the subset of elements of ${\mathbb N}^{n}$ which are of the form $f^{n}(x)$ for some $x\in R_{M}$ and injective function $f:M\to {\mathbb N}$), and with the geometric formulae (up to $\mathbb S$-provable equivalence) over the signature of $\mathbb S$ (via the map sending any such geometric formula $\phi(\vec{x})$ in $n$ to its interpretation $[[\phi(\vec{x})]]_{M}$ in $M$).      

\subsection{Atomless Boolean algebras}

Let $\mathbf{Bool}$ be (a skeleton of) the category of Boolean algebras and homomorphisms between them, that is the category of set-based models of the algebraic theory $\mathbb T$ of Boolean algebras. This theory is clearly of presheaf type, and its finitely presentable models are precisely the finite Boolean algebras. The injectivization of $\mathbb T$ is also of presheaf type (this follows for instance from our Remark \ref{injectivization}, and also from the criterion in \cite{vickers}); this implies that its category of models, which coincides with the category $\mathbf{Bool}_{i}$ of Boolean algebras and embeddings between them, is equivalent to the ind-completion of its full subcategory  $\mathbf{Bool}^{f}_{i}$ on the finite Boolean algebras. By using Stone duality between the category of finite Boolean algebras and embeddings between them and the category of finite sets and surjections, it is immediate to prove that the category $\mathbf{Bool}^{f}_{i}$ has the property that every arrow of it is a strict monomorphism. From Fra\"iss\'e's construction in Model Theory, we know that the unique countable atomless Boolean algebra $B$ is a $\mathbf{Bool}^{f}_{i}$-universal and $\mathbf{Bool}^{f}_{i}$-ultrahomogeneous object in $\mathbf{Bool}_{i}$. We thus obtain a concrete Galois theory for finite Boolean algebras.  

\subsection{Dense linear orders without endpoints} Let $\mathbf{LOrd}$ be (a skeleton of) the category of linear orders and order-preserving injections between them. It is easy to see that $\mathbf{LOrd}$ is the ind-completion of the full subcategory $\mathbf{LOrd}^{f}$ of it consisting of the finite linear orders. It is immediate to verify that every arrow of $\mathbf{LOrd}^{f}$ is a strict monomorphism, and that the set $\mathbb Q$ of rational numbers, endowed with its usual ordering, is a $\mathbf{LOrd}^{f}$-universal and $\mathbf{LOrd}^{f}$-ultrahomogeneous object in $\mathbf{LOrd}$. We thus have a concrete Galois theory for finite linear orders.

\subsection{Universal locally finite groups}

Let $\cal C$ be (a skeleton of) the category of finite groups and injective homomorphisms between them. The ind-completion of $\cal C$ can be identified with the category of ${\cal L}$ locally finite groups (i.e., the groups such that any finitely generated subgroup of them is finite) and injective homomorphisms between them. 

As remarked in \cite{Hodges}, the fact that $\cal C$ satisfies the amalgamation property (and hence also the joint embedding property) can be proved using the permutation products of B. H. Neumann (see section 3 
of \cite{neumann}). Also, it is clear that $\cal C$ satisfies all the other hypotheses of Fra\"iss\'e's construction (cf. Theorem 7.1.2 \cite{Hodges}). The Fra\"iss\'e's of $\cal C$ is known as Philip Hall's 
universal locally finite group (cf. \cite{hall}); this group is countable, simple, and any two isomorphic finite subgroups are conjugate.

A locally finite group $G$ is said to be \emph{universal} if $(a)$ every finite group is 
embeddable in $G$ and $(b)$ every isomorphism between finite subgroups of $G$ can be extended to some inner automorphism of $G$. It follows at once that $G$, regarded as an object of $\cal L$, is $\cal C$-universal and $\cal C$-ultrahomogeneous. Universal locally finite groups are known to exist in all infinite cardinalities, and Fra\"iss\'e's theorem implies that there is exactly one countable universal locally finite group (up to isomorphism), namely Philip Hall's group.  

For any universal locally finite group $G$, Theorem \ref{maincategorical} thus provides an equivalence of toposes $\Sh({\cal C}^{\textrm{op}}, J_{at})\simeq \Cont(Aut(G))$, where $Aut(G)$ is endowed with the topology of pointwise convergence. Such an equivalence always yields a concrete Galois theory for the category $\cal C$ (that is, ${\cal C}^{\textrm{op}}$ embeds as a full subcategory of the topos $\Cont(Aut(G))$), as any arrow of $\cal C$ is a strict monomorphism (cf. Exercise 7H(a) \cite{adamek}). These Galois theories are non-standard; in other words, there are `imaginaries' for the theory of finite groups. To prove this, we use the criterion for proving the atomic completeness of a category provided by Theorem \ref{atomcompl}. Clearly, the category $\cal C$ has equalizers and arbitrary intersections of subobjects; we thus have to find a pair of arrows $h,k:c\to e$ with equalizer $m:d\to c$ and a pair of arrows $k,n:c\to e'$ with $l\circ m=n\circ m$ and such that for any arrow $s:e'\to e''$, $(s\circ l, s\circ n)$ does not belong to the equivalence relation $R^{e''}_{h,k}$ on $Hom_{\cal C}(c, e'')$ generated by the relation consisting of the pairs of the form $(t\circ h, t\circ k)$ for an arrow $t:e \to e''$. Let us take $c=e$ to be the additive group $\mathbb Z\slash 15\mathbb Z$, $h$ to be the identity on $c$ and $k$ to be the function $[x]\to [2x]$; clearly, these maps are arrows in $\cal C$, and the equalizer $m$ of $h$ and $k$ is the trivial group. Take $e'$ to be the product group $\mathbb Z\slash 15\mathbb Z \times \mathbb Z\slash 15\mathbb Z$, $l$ to be the arrow $[x]\to ([x], [3x])$ and $n$ to be the arrow $[y]\to ([5y], [y])$; obviously, $l\circ m=n\circ m$. Now, for any $a,b\in Hom_{\cal C}(c, e'')$, $(a,b)\in R^{e''}_{h,k}$ clearly implies that for any element $z$ of $c$ there exist natural numbers $p$ and $q$ such that $2^{p}a(z)=2^{q}b(z)$. From this it follows that for any arrow $s:e'\to e''$, if $(s\circ l, s\circ n)$ belongs to $R^{e''}_{h,k}$ then there exist natural numbers $p$ and $q$ such that $2^{p}(s\circ l)([1])=2^{q}(s\circ n)([1])$. But $2^{p}(s\circ l)([1])=s(2^{p}l([1]))$ and $2^{q}(s\circ n)([1])=s(2^{q}n([1]))$, whence, as $s$ is injective, $2^{p}l([1])=2^{q}n([1])$, which is not true since by definition of $l$ and $n$ we have $2^{p}l([1])=([2^{p}], [2^{p}\cdot 3])$ and $2^{q}l([1])=([2^{q}\cdot 5], [2^{q}])$.

\subsection{The random graph} 

Let $\cal C$ be (a skeleton of) the category of finite (irreflexive) graphs and injective homomorphisms between them. It is easy to verify that $\cal C$ satisfies the amalgamation and joint embedding properties, and that its ind-completion can be identified with the category of (irreflexive) graphs and injective homomorphisms between them. Also, one can easily prove that every arrow in $\cal C$ is a strict monomorphism. The \emph{random graph} (also known as the Rado graph) is a $\cal C$-universal and $\cal C$-ultrahomogeneous; indeed, it can be identified with the Fra\"iss\'e's limit of the class of finite graphs (cf. Theorem 7.1.2 and Theorem 7.4.4 \cite{Hodges}). We thus have a concrete Galois theory for finite graphs.

\subsection{Infinite Galois theory} Let $F\subset L$ be a Galois extension (not necessarily finite-dimensional), and let ${\cal L}_{F}^{L}$ be the category of finite intermediate extensions and field homomorphisms between them. The ind-completion $\Ind{{\cal L}_{F}^{L}}$ of ${\cal L}_{F}^{L}$ can be identified with the category of intermediate (separable) extensions of $F$, and $L$ is a ${\cal L}_{F}^{L}$-universal and ${\cal L}_{F}^{L}$-ultrahomogeneous object in $\Ind{{\cal L}_{F}^{L}}$. The fundamental theorem of classical Galois theory ensures that all the arrows in ${\cal L}_{F}^{L}$ are strict monomorphisms and that the category ${\cal L}_{F}^{L}$ is atomically complete (in fact, the theorem is equivalent to the conjunction of these two assertions). Notice that if $L$ is finite-dimensional then the resulting Galois theory is discrete. Of course, infinite Galois theory also falls into the framework of Grothendieck's Galois theory (see below).

\subsection{Grothendieck's Galois Theory}\label{GrothGalois} Let $\cal C$ be a Galois category (in the sense of \cite{grothendieckg}), with fibre functor $F:{\cal C}\to \Set$, and let ${\cal C}_{t}$ be the full subcategory of $\cal C$ formed by the objects $c$ of $\cal C$ such that $F(c)$ is non-empty and the action of $Aut(F)$ on $F(c)$ is transitive (equivalently, by the atomic objects of $\cal C$, i.e. the objects of $\cal C$ which are not isomorphic to $0$ and which do not admit any proper subobjects). Grothendieck's theory provides an equivalence ${\cal C}\simeq \Cont_{f}(Aut(F))$, where $\Cont_{f}(Aut(F))$ is the category of continuous actions of the profinite group $Aut(F)$ on a finite set. The open subgroups of $Aut(F)$ are exactly those of the form $U_{c, x}:=\{\alpha:F\cong F \textrm{ | } \alpha_{c}(x)=x\}$ for $c\in {\cal C}$ and $x\in F(c)$; from this it follows in particular that every transitive continuous action of $Aut(F)$ lies in $\Cont_{f}(Aut(F))$ (that is, its underlying set is finite). 

The equivalence ${\cal C}\simeq \Cont_{f}(Aut(F))$ thus restricts to an equivalence ${\cal C}_{t}\simeq \Cont_{t}(Aut(F))$, where $\Cont_{t}(Aut(F))$ is the category of continuous transitive actions of $Aut(F)$ over discrete sets. On the other hand, the equivalence ${\cal C}\simeq \Cont_{f}(Aut(F))$ can be extended to an equivalence of toposes by equipping each of the two categories with the coherent topology. By the Comparison Lemma, the category of sheaves on $\Cont_{f}(Aut(F))$ with respect to the atomic topology is equivalent to the topos $\Cont(Aut(F))$ of continuous actions of $Aut(F)$; indeed, as every continuous action can be written as a coproduct of transitive ones, and every transitive action lies in $\Cont_{f}(Aut(F))$, the subcategory $\Cont_{f}(Aut(F))$ of $\Cont(Aut(F))$ is dense with respect to the canonical topology on $\Cont(Aut(F))$, and the induced Grothendieck topology on $\Cont_{f}(Aut(F))$ coincides with the coherent one. Thus we have an equivalence of toposes
\[
\Sh({\cal C}, J_{coh})\simeq \Cont(Aut(F)),
\]    
where $J_{coh}$ is the coherent topology on $\cal C$.

Another application of the Comparison Lemma yields an equivalence
\[
\Sh({\cal C}, J_{coh})\simeq \Sh({\cal C}_{t}, J_{at}),
\]
where $J_{at}$ is the atomic topology on ${\cal C}_{t}$. We thus have an equivalence
\[
\Sh({\cal C}_{t}, J_{at}) \simeq \Cont(Aut(F)).
\]

These equivalences show that $F$ corresponds to a ${{\cal C}_{t}}^{\textrm{op}}$-universal and ${{\cal C}_{t}}^{\textrm{op}}$-ultrahomogeneous object of $\Ind{{{\cal C}_{t}}^{\textrm{op}}}$, that all the arrows of ${{\cal C}_{t}}^{\textrm{op}}$ are strict monomorphisms and that the category ${{\cal C}_{t}}^{\textrm{op}}$ is atomically complete. In particular, the canonical embedding of ${\cal C}_{t}$ into the topos $\Sh({\cal C}_{t}, J_{at})\simeq \Sh({\cal C}, J_{coh})$ can be identified with the free small coproduct-completion of ${\cal C}_{t}$ (cf. Proposition \ref{freecoproduct}), while the canonical embedding ${\cal C}_{t}\hookrightarrow {\cal C}$ realizes $\cal C$ as the free finite-coproduct completion of the category ${\cal C}_{t}$. 

Notice that, as $Aut(F)$ is a profinite group, the topos $\Cont(Aut(F))$ is coherent (see section D3.4 \cite{El}). In fact, the category $\cal C$ is a pretopos, as it is equivalent to the full subcategory of $\Cont(Aut(F))$ on its coherent objects; indeed, from the proof of Lemma D3.4.2 \cite{El} we know that the coherent objects of the topos $\Cont(Aut(F))$ are exactly the compact objects, and it is clear that an object of $\Cont(Aut(F))$ is compact if and only if it lies in $\Cont_{f}(Aut(F))$.  

The discussion above shows in particular that the Galois categories can be characterized as the pretoposes $\cal D$ in which every object can be written as a coproduct of atomic objects and with the property that the opposite ${\cal D}_{a}$ of their subcategory of atomic objects satisfies AP and JEP and there is a ${\cal D}_{a}$-universal and ${\cal D}_{a}$-ultrahomogeneous object in $\Ind{{\cal D}_{a}}$. Indeed, we observed above that any Galois category has this property, while the converse follows from the fact that if a category $\cal D$ satisfies this condition then, denoted by $u$ the ${\cal D}_{a}$-universal and ${\cal D}_{a}$-ultrahomogeneous object in $\Ind{{\cal D}_{a}}$, we have, by the Comparison Lemma and Theorem \ref{maincategorical}, equivalences 
\[
\Sh({\cal D}, J_{coh})\simeq \Sh({{\cal D}_{a}}^{\textrm{op}}, J_{at})\simeq \Cont(Aut_{{\cal D}_{a}}(u)),
\] 
and hence ${\cal D}$ is equivalent to $\Cont_{f}(Aut_{{\cal D}_{a}}(u))$ (recall that any pretopos can be recovered, up to equivalence, as the category of coherent objects of the topos of coherent sheaves on it).       

The infinitary version of Grothendieck's Galois theory obtained in \cite{noohi} admits a similar interpretation in the context of our framework; Noohi's Galois categories $\cal C$ are themselves toposes which can be represented as categories of atomic sheaves on the subcategory ${\cal C}_{t}$ of atomic objects of $\cal C$.

\section{The four ways to topological Galois representations}\label{four}

Our abstract approach to Galois theory is based on the intuition that the generating `kernel' of a Galois theory lies at the topos-theoretic level; more precisely, it is expressed by an equivalence between a topos of sheaves on an atomic site $\Sh({\cal C}^{\textrm{op}}, J_{at})$ and a topos $\Cont(G)$ of continuous actions of a topological group $G$. Starting from such an equivalence
\[
\Sh({\cal C}^{\textrm{op}}, J_{at})\simeq \Cont(G),
\]  
one can establish relationships between the category $\cal C$ and the group $G$ by considering appropriate topos-theoretic invariants in terms of the two different representations of the classifying topos of the given `Galois-type' theory; in other words, the given topos acts as a `bridge' (in the sense of \cite{OC10}) for transferring information between the two sides. 

In particular, it is natural to investigate under which conditions the given equivalence restricts to a categorical equivalence at the level of sites, leading to a `concrete' Galois theory (we did this in section \ref{concrete}). The fact that every concrete Galois theory can be obtained as a restriction of a topos-theoretic `Galois-type' equivalence and that, on the other hand, there are many equivalences of the latter kind which do not specialize to equivalences of sites (but are nonetheless relevant for other purposes, notably including the computation of cohomological invariants), shows that this `top-down', topos-theoretic viewpoint carries with itself a higher generality and a greater technical flexibility with respect to the categorical, site-level, analysis. 

Another illustration of the naturality of our topos-theoretic viewpoint is provided by the fact that, as a given topos can admit several different representations, one can arrive at establishing such Galois-type topos-theoretic equivalences in a variety of different ways. We shall briefly discuss these different ways below, organizing them in four main groups.  
         
\subsection{Representation theory of Grothendieck toposes}

The representation theory of Grothendieck toposes provides natural ways for generating Galois-type equivalences of toposes. For instance, as we saw in section \ref{discreteness}, discrete Galois theories can all be obtained as instances of Grothendieck's Comparison Lemma. Another natural source of Galois-type equivalences of the above kind is provided by the representation theory of atomic connected toposes in terms of localic groups established by Dubuc in \cite{dubuc3}; indeed, these representations yield, in the case of spatial localic groups, topological Galois representations fitting into our framework.

\subsection{Ultrahomogeneous structures}

Our theorems \ref{maincategorical} and \ref{maincategoricalD} establish a central connection between the theory of ultrahomogeneous structures and the subject of Galois representations. Ultrahomogeneous structures naturally arise in a great variety of different mathematical contexts (cf. section \ref{examples} for a list of notable examples). Their existence can be proved either directly through an explicit construction or through abstract logical arguments (see \cite{Hodges}). A general method for building countable ultrahomogeneous structures is provided by Fra\"iss\'e's construction in Model Theory (cf. chapter 7 of \cite{Hodges}), while the categorical generalization established in \cite{OC2} allows to construct ultrahomogeneous structures of arbitrary cardinality.

\subsection{Special models}

Recall from section \ref{repres} that a model $M$ of an atomic complete theory $\mathbb T$ is said to be \emph{special} if every $\mathbb T$-complete formula $\phi(\vec{x})$ is realized in $M$ and for any tuples $\vec{a}$ and $\vec{b}$ of elements of $M$ which satisfy the same $\mathbb T$-complete formulae there is an automorphism $f:M\to M$ of $M$ which sends $\vec{a}$ to $\vec{b}$. Theorem \ref{main} shows that, in presence of any special model for an atomic complete theory, we have a Galois-type representation for its classifying topos.

An interesting aspect of special models is that their existence can be proved in many contexts by using model-theoretic techniques. For instance, the unique (up to isomorphism) countable model of a countably categorical theory is always special. Also, it follows from the fact that every consistent finitary first-order theory has a $\omega$-big model (cf. chapter 10 of \cite{Hodges}) that every consistent finitary atomic complete theory has a special model. This implies that for any consistent (i.e., with at least a set-based model) finitary theory of presheaf type $\mathbb T$ satisfying the hypotheses of Theorem \ref{corollarypresheaf} there exists a $\textrm{f.p.}{\mathbb T}\textrm{-mod}(\Set)$-universal and $\textrm{f.p.}{\mathbb T}\textrm{-mod}(\Set)$-ultrahomogeneous model of $\mathbb T$.

Notice also that, since the concept of special model for an atomic complete theory is expressible through a topos-theoretic invariant (cf. Remark \ref{rem}(c)), one can obtain special models for a theory $\mathbb T$ starting from special models of any theory ${\mathbb T}'$ which is Morita-equivalent to it.  

Another result from chapter 10 of \cite{Hodges} is that any consistent finitary first-order theory has a $\lambda$-big model for any infinite cardinal $\lambda$; in particular, it has a model $M$ such that every model of $\mathbb T$ of cardinality less than $\lambda$ can be elementarily embedded in $M$ and any isomorphism between any two such models can be extended to an automorphism of $M$. Let ${\mathbb T}\textrm{-mod}_{e}(\Set)$ be the category of $\mathbb T$-models in $\Set$ and elementary embeddings between them. In order to obtain a Galois-type equivalence from these data by applying Theorem \ref{catd}, we would need the embedding ${\mathbb T}\textrm{-mod}_{\lambda}(\Set)\hookrightarrow {\mathbb T}\textrm{-mod}_{e}(\Set)$, where ${\mathbb T}\textrm{-mod}_{\lambda}(\Set)$ is the full subcategory of ${\mathbb T}\textrm{-mod}_{e}(\Set)$ on the $\mathbb T$-models of cardinality less than $\lambda$, together with the object $M$ of ${\mathbb T}\textrm{-mod}_{e}(\Set)$, to satisfy the hypotheses of the theorem. Notice that if $\lambda$ is chosen in such a way that all finitely generated $\mathbb T$-models have cardinality less than $\lambda$ (note that the cardinality of a finitely generated model of $\mathbb T$ is always $\lt \omega + card(L)$, where $card(L)$ is the cardinality of the signature of $\mathbb T$, cf. Theorem 1.2.3 \cite{Hodges}) and every $\mathbb T$-model can be written as a filtered colimit in ${\mathbb T}\textrm{-mod}_{e}(\Set)$ of finitely generated models (notice that this condition is always satisfied if $\mathbb T$ is the theory of homogeneous models of a theory of presheaf type, as such a theory is Boolean being atomic and hence all the $\mathbb T$-model homomorphisms are elementary embeddings) then the embedding ${\mathbb T}\textrm{-mod}_{f.g.}(\Set)\hookrightarrow {\mathbb T}\textrm{-mod}_{e}(\Set)$, where ${\mathbb T}\textrm{-mod}_{f.g.}(\Set)$ is the full subcategory of ${\mathbb T}\textrm{-mod}_{e}(\Set)$ on the finitely generated $\mathbb T$-models satisfies, together with the object $M$, the hypotheses of Theorem \ref{catd}, provided that the category ${\mathbb T}\textrm{-mod}_{f.g.}(\Set)$ satisfies the amalgamation and joint embedding properties and every arrow of it is a strict monomorphism.

\subsection{Galois categories} 

As remarked in \ref{GrothGalois}, to every Galois category we can naturally associate an equivalence of toposes which `materializes' the relationship between it and the associated Galois group. In fact, the theory of Galois categories provides a significant number of examples of standard Galois theories, as the subcategories of atomic objects of a Galois category are atomically complete.

\section{Conclusions and future work}\label{conclusions} 

The theory developed in this paper provides simple means for constructing `Galois-type' equivalences in a variety of different mathematical contexts. In fact, Galois-type theories (i.e., theories whose classifying topos can be represented as a topos of continuous actions of a topological group) are ubiquitous in Mathematics; indeed, they are maximal elements in the lattice of geometric theories over their signature, and every geometric theory over a countable signature whose Booleanization is consistent can be extended to a Galois-type theory (cf. \cite{OC2}, \cite{OC5}, \cite{OC6} and section \ref{intro} above). A natural direction of future investigation thus consists in trying to generate new meaningful applications of our general machinery across different mathematical domains; in particular, we intend to study Connes' topos of cyclic sets (\cite{connes1} and \cite{connes2}), exploring the possibility of building an associated Galois theory. Also, we plan to investigate the notion of random group, as introduced by M. Gromov in \cite{gromov}, in relation to the purpose of building Galois-type theories for various classes of discrete groups.         

On the theoretical side, we plan to investigate the possibility of extending our theory by replacing topological groups with topological groupoids, in order to obtain generalized Galois-type representations holding for arbitrary geometric theories with enough models; this will be achieved by generalizing Butz and Moeridjk's representation theorem \cite{BM} for Grothendieck toposes with enough points in the direction of allowing more flexibility in the construction of the groupoids presenting the toposes. We also intend to analyze the relationships between the localic groups of automorphisms of points of atomic toposes introduced in \cite{dubuc3} and the topological automorphism groups considered in the present work. Finally, we would like to explore possible connections between the theory developed in this paper and the categorical Galois theory of Borceux and Janelidze \cite{borjan}, which is based on a rather different framework.  

\vspace{0.5cm}

\textbf{Acknowledgements:} I am grateful to Jean-Pierre Serre for kindly providing me with two different proofs of the (standard) fact that in the category of finite groups and embeddings between them every arrow is a strict monomorphism, and to Laurent Lafforgue for useful discussions on the subject of `imaginaries' in the Galois theory for finite groups.

\end{document}